\documentclass[article,ij4uq]{ij4uq}     

\usepackage[hang]{footmisc}
\fancypagestyle{plain}{%
  \fancyhf{}
  \fancyhead[R]{\small {\it International Journal for Uncertainty Quantification}, x(x): \thepage--\pageref{LastPage} (\myyear\today)}
  \fancyfoot[R]{\small\bf\thepage }
  \fancyfoot[L]{\fottitle}
  }

\renewcommand{\dmyy}{17}
\renewcommand{\myyear}{2025}
\renewcommand{\today}{}

\numberwithin{equation}{section}

\usepackage{lipsum}
\usepackage{amsfonts}
\usepackage{graphicx}
\usepackage{epstopdf}
\usepackage{algorithmic}
\ifpdf
  \DeclareGraphicsExtensions{.eps,.pdf,.png,.jpg}
\else
  \DeclareGraphicsExtensions{.eps}
\fi
\usepackage{amsmath,amssymb, amsfonts}
\usepackage{graphicx}
\usepackage{booktabs}

\renewcommand{\vec}[1] {\ensuremath{\boldsymbol{#1}}}

\usepackage{amsopn}

\def\R{{\mathbb R}}
\def\U{{\mathcal U}}
\def\Z{{\mathcal Z}}
\def\t{{\vec{\theta}}}
\def\u{{\vec{u}}}

\def\z{{\vec{z}}}

\def\ztilde{{\tilde{\z}}}

\def\v{{\vec{v}}}

\def\d{{\delta}}
\def\A{{\vec{A}}}
\def\B{{\vec{B}}}
\def\E{{\vec{E}}}
\def\M{{\vec{M}}}
\def\K{{\vec{K}}}

\def\S{{\vec{S}}}

\def\L{{\vec{L}}}

\def\H{{\vec{H}}}
\def\P{{\vec{P}}}

\def\V{{\vec{V}}}

\def\I{{\vec{I}}}

\def\W{{\vec{W}}}
\def\T{{\vec{T}}}

\def\F{{\vec{F}}}

\usepackage{array}
\newcolumntype{P}[1]{>{\centering\arraybackslash}p{#1}}
\newcolumntype{M}[1]{>{\centering\arraybackslash}m{#1}}

\begin{document}

\volume{Volume x, Issue x, \myyear\today}
\title{Hyper-differential sensitivity analysis with respect to model discrepancy: Prior Distributions}
\titlehead{Discrepancy Prior Distributions}
\authorhead{Joseph Hart, Bart van Bloemen Waanders, Jixian Li, Timbwaoga A. J. Ouermi, and Chris R. Johnson}

\corrauthor[1]{Joseph Hart}
\corremail{joshart@sandia.gov}
\author[1]{Bart van Bloemen Waanders}
\author[2]{Jixian Li}
\author[2]{Timbwaoga A. J. Ouermi}
\author[2]{Chris R. Johnson}

\address[1]{Department of Scientific Machine Learning, Sandia National Laboratories}
\address[2]{Scientific Computing and Imaging Institute, University of Utah}

\abstract{
Hyper-differential sensitivity analysis with respect to model discrepancy was recently developed to enable uncertainty quantification for optimization problems. The approach consists of two primary steps: (i) Bayesian calibration of the discrepancy between high- and low-fidelity models, and (ii) propagating the model discrepancy uncertainty through the optimization problem. When high-fidelity model evaluations are limited, as is common in practice, the prior discrepancy distribution plays a crucial role in the uncertainty analysis. However, specification of this prior is challenging due to its mathematical complexity and many hyper-parameters. This article presents a novel approach to specify the prior distribution. Our approach consists of two parts: (1) an algorithmic initialization of the prior hyper-parameters that uses existing data to initialize a hyper-parameter estimate, and (2) a visualization framework to systematically explore properties of the prior and guide tuning of the hyper-parameters to ensure that the prior captures the appropriate range of uncertainty. We provide detailed mathematical analysis and a collection of numerical examples that elucidate properties of the prior that are crucial to ensure uncertainty quantification.}


\maketitle

\section{Introduction}
Optimization constrained by computational models, such as partial differential equations, arises in a wide range of science and engineering applications. Although numerical optimization has a rich theoretical and algorithmic foundation, many problems remain unsolved due to the complexity of the computational models involved. For instance, state-of-the-art optimization algorithms require derivative information to efficiently traverse high-dimensional spaces. However, many complex computational models are developed with a focus on forward prediction and consequently fail to provide access to derivatives. Furthermore, derivative-based optimization algorithms still require many model evaluations to converge, which is prohibitive when each evaluation requires significant computational resources. Consequently, a common approach to enable computationally tractable optimization is to work with a lower-fidelity model that can provide derivative information and be evaluated more efficiently. 

Hyper-differential sensitivity analysis with respect to model discrepancy was first developed in~\cite{hart_bvw_cmame} to leverage a small number of high-fidelity model evaluations (without derivatives) to enhance an optimization solution computed using a low-fidelity model. Because high-fidelity model evaluations are limited,~\cite{hart_bvw_mods} extended the approach to facilitate uncertainty quantification through a Bayesian approach. Specifically, the model discrepancy, defined as the difference in the high- and low-fidelity models, is estimated using Bayesian calibration. The resulting posterior distribution is propagated through the optimization problem to characterize the most likely optimization solutions given the high-fidelity model evaluations. This framework leverages tools from numerical linear algebra to achieve computational scalability for high-dimensional optimization problems involving nonlinear partial differential equations; the computational complexity scales with rank rather than the native dimension of the problem.

The aforementioned works provide a computational framework to enable uncertainty quantification of the optimization solution. However, because the number of high-fidelity model evaluations are frequently limited, the model discrepancy calibration is ill-posed and consequently the quality of the uncertainty estimate depends heavily on the prior. Specifying this prior is difficult because the mathematical object being calibrated is an operator that maps from the optimization variable space to the state variable space. Moreover, the prior depends on many hyper-parameters, involving both the state and optimization variables, which require considerable expertise and manual tuning to determine appropriate values. This article focuses on specification of the prior discrepancy distribution and provides a pragmatic framework to facilitate analysis for complex systems. 

Our contributions include:
\begin{enumerate}
\item[$\bullet$] Proving that the prior discrepancy distribution is well-defined in the discretization refinement limit.
\item[$\bullet$] Developing an algorithmic approach to initialize the prior hyper-parameters and thus reduce the difficultly of manual tuning of the prior.
\item[$\bullet$] Proposing a visualization framework to facilitate analysis of the prior and guide any hyper-parameter modifications required after the algorithmic initialization. 
\item[$\bullet$] Demonstrating properties of the prior and providing insight through a hierarchy of numerical examples.
\end{enumerate}

Hyper-differential sensitivity analysis with respect to model discrepancy (HDSA-MD) is related to multifidelity methods~\cite{multifidelity_review_peherstorfer}, and in particular, multifidelity optimization approaches such as~\cite{multifidelity_quasinewton_bryson,trmm_lewis,biegler_rom_opt,willcox_multifi_opt_2012}. However, HDSA-MD differs in that it does not require derivative information from the high-fidelity model and it does not require access to query the high-fidelity model. Rather, HDSA-MD works with high-fidelity data, which is provided from offline evaluations and can be used with as little as one high-fidelity solve. Because HDSA-MD places few requirements on the source of the high-fidelity data, it puts greater emphasis on the prior discrepancy distribution to enrich the high-fidelity data. Bayesian modeling of model discrepancy is a classical concept pioneered by Kennedy and O'Hagan~\cite{ohagan2001} and subsequently studied by many others~\cite{ohagan2002,ohagan2014,Higdon_2008,Arendt_2012,Maupin,Ling_2014,Gardner_2021}. Although our approach has some similarities, we differ from this literature in our focus on optimization problems, which leads to differences in the discrepancy representation and calibration. 

The article is organized as follows. Section~\ref{sec:background} introduces the problem formulation and provides relevant background information. The Bayesian discrepancy calibration problem is discussed in Section~\ref{sec:bayes_calibration} where we show that the prior discrepancy distribution is well defined. Section~\ref{sec:hyperparam_selection} presents our proposed approach for hyper-parameter specification through a combination of algorithm-based initialization and prior sample visualization. Numerical results demonstrate properties of the prior in Section~\ref{sec:numerics} and concluding remarks are made in Section~\ref{sec:conclusion}.

\section{Problem Formulation and Background}  \label{sec:background}

We consider optimization problems of the form
\begin{align}
\label{eqn:true_opt_prob}
& \min_{z \in \Z} J(S(z),z) 
\end{align}
where $z$ is an optimization variable in a Hilbert space $\Z$, $S:\Z \to \U$ is the solution operator for a differential equation whose state variable $u$ is in a Hilbert space $\U$, and $J:\U \times \Z \to \R$ is the objective function. 

Solving~\eqref{eqn:true_opt_prob} for high-fidelity models is frequently intractable due to the complexity of the high-fidelity model; many models require considerable computational expense per evaluation and do not provide access to derivative information. In such cases, a common alternative is to solve
\begin{align}
\label{eqn:approx_opt_prob}
& \min_{z \in \Z} J(\tilde{S}(z),z) 
\end{align}
where $\tilde{S}:\Z \to \U$ is the solution operator for a low-fidelity model. The goal of HDSA-MD is to approximate the solution of~\eqref{eqn:true_opt_prob}, and characterize uncertainty in this optimal solution approximation, by updating the solution of~\eqref{eqn:approx_opt_prob} using limited evaluations of the high-fidelity model.

The state space $\U$ is infinite-dimensional, and the optimization space $\Z$ can be finite- or infinite-dimensional. We discretize by letting $\{\eta_i\}_{i=1}^{n_u}$ be a basis for a finite-dimensional subspace $\U_h \subset \U$ and $\{\varphi_j\}_{j=1}^{n_z}$ be a basis for a finite-dimensional subspace $\Z_h \subset \Z$. For finite-dimensional optimization spaces, we have $\Z_h = \Z$. Our subsequent analysis focuses on $\U$ and $\Z$ being Sobolev spaces such as $H^1$, but we note that $\Z$ can be a finite-dimensional inner product space, such as Euclidean space. Our analysis is valid for such cases, with some simplifications that are highlighted in the article. The coordinate representations of $u$ and $z$ are denoted as $\u \in \R^{n_u}$ and $\z \in \R^{n_z}$. For simplicity, we abuse notation and use $J:\R^{n_u} \times \R^{n_z} \to \R$, $S:\R^{n_z} \to \R^{n_u}$, and $\tilde{S}:\R^{n_z} \to \R^{n_u}$ to denote the discretized objective function and solution operators. 

To implement the $\U$ and $\Z$ inner products in the coordinate space, we define $\M_\u \in \R^{n_u \times n_u}$ and $\M_\z \in \R^{n_z \times n_z}$ as $(\M_\u)_{i,j} = (\eta_i,\eta_j)_{\mathcal{U}}$ and $(\M_\z)_{i,j} = (\varphi_i,\varphi_j)_{\mathcal{Z}}$. The gradient of a scalar-valued function (or the Jacobian of a vector-valued function) with respect to $\u$ will be denoted as $\nabla_\u$. Similarly, Hessians with respect to $(\u,\u)$ and $(\u,\z)$ will be denoted as $\nabla_{\u,\u}$ and $\nabla_{\u,\z}$, respectively. There is a large body of literature devoted to solving optimization problems constrained by computational models. We refer the reader to~\cite{Vogel_99, Archer_01,Haber_01,Vogel_02,Biegler_03,Biros_05,Laird_05,Hintermuller_05,Hazra_06,Biegler_07,Borzi_07,Hinze_09,Biegler_11,frontier_in_pdeco} for an overview.

\subsection{Post-optimality sensitivity with respect to model discrepancy} 
To represent model discrepancy, the difference between $S(\z)$ and $\tilde{S}(\z)$, in the optimization problem, we consider 
\begin{align}
\label{eqn:dis_approx_opt_prob_pert_rs}
 \min_{\z \in \R^{n_z}} \hspace{1 mm} \mathcal{J}(\z,\t):=J(\tilde{S}(\z)+\d(\z,\t),\z) ,
\end{align}
where $\d:\R^{n_z} \times \R^{n_\theta} \to \R^{n_u}$ is an operator that we parameterize by $\t \in \R^{n_\theta}$ and calibrate so that $\d(\z,\t) \approx S(\z)-\tilde{S}(\z)$. To relate~\eqref{eqn:dis_approx_opt_prob_pert_rs} to the low-fidelity optimization problem, we parameterize $\d$ so that $\d(\z,\vec{0})=\vec{0}$ $\forall \z$. To facilitate our subsequent derivations, assume that $\mathcal J(\z,\t)$ is twice continuously differentiable with respect to $(\z,\t)$. Let $\ztilde$ be a local minimizer of the low-fidelity optimization problem, which satisfies the first order optimality condition
\begin{align}
\label{eqn:first_order_opt}
\nabla_\z \mathcal{J}(\ztilde,\vec{0}) = \vec{0}.
\end{align}
Applying the Implicit Function Theorem to~\eqref{eqn:first_order_opt} implies that there exists an operator $\F:
\R^{n_\theta} \to \R^{n_z}$ such that
\begin{align*}
\nabla_\z \mathcal{J}(\F(\t),\t) = \vec{0} \qquad \forall \t \in \mathcal N(\vec{0}),
\end{align*}
where $ \mathcal N(\vec{0})$ denotes a neighborhood of $\vec{0} \in \R^{n_\theta}$. Assuming that the Hessian $\nabla_{\z,\z} \mathcal{J}(\F(\t),\t)$ is positive definite $\forall \t \in \mathcal N(\vec{0})$, we can interpret $\F$ as a mapping from the model discrepancy to the optimal solution. Furthermore, the Implicit Function Theorem provides the Jacobian of $\F$,
\begin{eqnarray}
\label{eqn:dis_sen_op}
\nabla_\t \F(\vec{0}) = - \H^{-1} \B \in \R^{n_z \times n_\theta},
\end{eqnarray}
where $\H=\nabla_{\z,\z} \mathcal{J}(\ztilde,\vec{0})$ and $\B=\nabla_{\z,\t} \mathcal{J}(\ztilde,\vec{0})$ are second order derivatives of $\mathcal{J}$, evaluated at $(\ztilde,\vec{0})$. We refer to $\nabla_\t \F$ as the post-optimality sensitivity operator, which has been the subject of different studies~\cite{shapiro_SIAM_review,fiacco,Griesse_part_1,Griesse_part_2,griesse2}.

Assume that $\ztilde$ has been computed, that $S$ and $\tilde{S}$ have been evaluated $N$ (a small number) times to provide a dataset $\{ \z_\ell,S(\z_\ell)-\tilde{S}(\z_\ell)\}_{\ell=1}^N$, and that $\z_1=\ztilde$ \footnote{The assumption that $\z_1=\ztilde$ may be relaxed. However, it is convenient for subsequent analysis and is typically desirable in practice.}. Hyper-differential sensitivity analysis with respect to model discrepancy is summarized as the two-step procedure:
\begin{itemize}
\item[] (discrepancy calibration) determine $\t \in \R^{n_\theta}$ such that 
$$\d(\z_\ell,\t) \approx S(\z_\ell)-\tilde{S}(\z_\ell), \qquad \ell=1,2,\dots,N,$$

\item[] (optimal solution updating) approximate the high-fidelity optimization solution as
$$\ztilde - \H^{-1} \B \t.$$

\end{itemize}

Since the number of high-fidelity model evaluations $N$ is typically small, the calibration step is ill-posed in that there may be many $\t$'s such that the calibrated $\d$ matches the data. To incorporate prior information, and quantify uncertainty due to the data limitations, we pose the calibration step in a Bayesian framework. The result is a probabilistic description of $\t$\textemdash the Bayesian posterior distribution. Then the updating step consists of approximately propagating the $\t$ posterior through the optimality system to determine a probability distribution for the optimization solution. It has been demonstrated that this approach is computationally efficient and can provide valuable optimal solution posteriors with as little as one high-fidelity model evaluation~\cite{hart_bvw_cmame,hart_bvw_mods}. 

\begin{figure}[h]
\centering
  \includegraphics[width=0.75\textwidth]{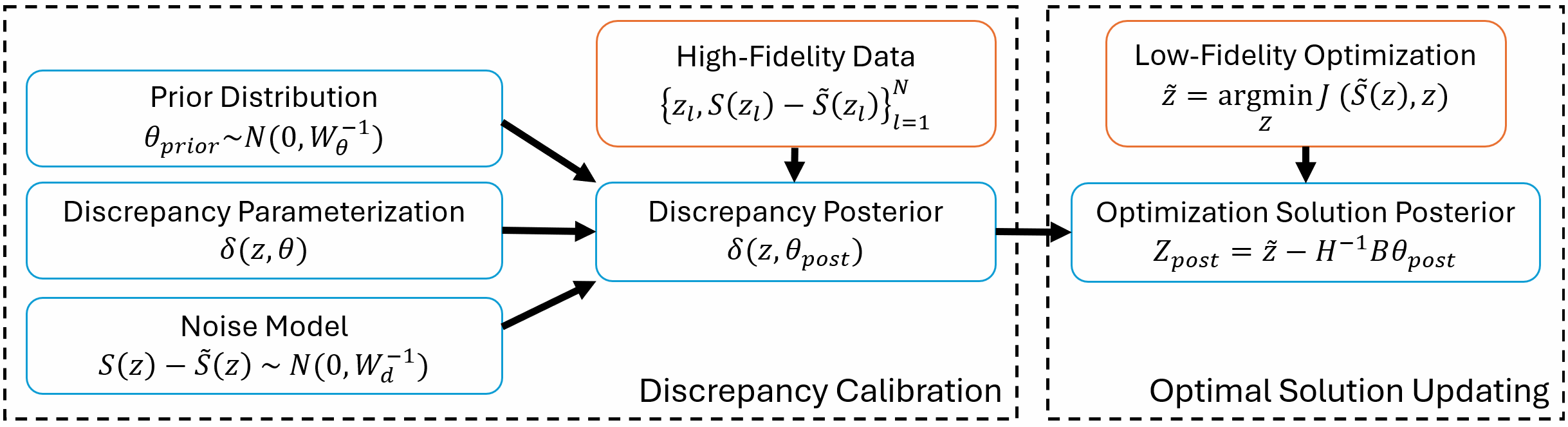}
    \caption{Depiction of the HDSA-MD framework.}
  \label{fig:hdsa_md_process}
\end{figure}

Figure~\ref{fig:hdsa_md_process} depicts the process of calibration and updating that combines low-fidelity optimization and $N$ high-fidelity forward solves to produce an optimization solution posterior distribution. The components of Figure~\ref{fig:hdsa_md_process} (the blue boxes) are interrelated in that the prior distribution for $\t$ depends on the parametric form of $\d(\z,\t)$, and the $\t$ posterior depends on the form of the prior and noise models. Furthermore, to achieve computational efficiency, it is critical that $\d(\z,\t)$ is parameterized in a way that facilitates both efficient posterior sampling and computation of the matrix vector products $\H^{-1} \B \t$ with the posterior samples. Since $\d$ is an operator mapping from the optimization variable space to the state space, the dimension of $\t$ will be $\mathcal O(n_u n_z).$ For problems where the optimization variable dimension $n_z$ is large, it is crucial to have a closed form expression for the $\t$ posterior samples to enable computation in the extremely high-dimensional parameter space $\R^{n_\theta}$.

This article focuses on the specification and properties of the prior distribution. We refer the reader to~\cite{hart_bvw_cmame,hart_bvw_mods} for details on other aspects of the HDSA-MD framework.
\subsection{Model discrepancy parameterization}
As a prerequisite to analyzing the prior distribution of $\t$, we recall the form of the discrepancy parameterization. Observe that the post-optimality sensitivity operator~\eqref{eqn:dis_sen_op} depends on $\H=\nabla_{\z,\z} \mathcal{J}(\ztilde,\vec{0})$ and $\B=\nabla_{\z,\t} \mathcal{J}(\ztilde,\vec{0})$. Since we assumed that $\d(\z,\vec{0})=\vec{0}$ $\forall \z$, $\H$ does not depend on the parameterization of $\d$. Furthermore, we see that $\B$ depends on the mixed $\z,\t$ second derivative of $\d$. Hence, it is sufficient to consider a parameterization that is bilinear in $(\z,\t)$, or equivalently, an affine function of $\z$ with coefficients parameterized by $\t$. Starting from a generic affine operator mapping between the Hilbert spaces $\Z_h$ and $\U_h$, we arrive at the coordinate representation
 \begin{align*}
 \d(\z,\t) = \u_0(\t) + \mathbf{L}(\t) \M_\z \z,
 \end{align*}
with
\begin{eqnarray}
\label{eqn:u0_L}
 \u_0(\t) = \t_0, \qquad \mathbf{L}(\t) =  \left(
 \begin{array}{c}
 \t_1^T \\
 \t_2^T \\
 \vdots \\
 \t_{n_u}^T
 \end{array} \right),
 \qquad  \text{and} \qquad
 \t = \left(
\begin{array}{c}
\t_0 \\
\t_1 \\
\vdots \\
\t_{n_u}
\end{array}
\right),
\end{eqnarray} 
where $\t_0 \in \R^{n_u}$ and $\t_i \in \R^{n_z}$, $i=1,2,\dots,n_u$. We can equivalently express $\d$ as
\begin{eqnarray}
\label{eqn:delta_kron}
\d(\z,\t) =
\left( \begin{array}{cc}
\I_{n_u} & \I_{n_u} \otimes \z^T \M_\z
\end{array} \right) \t
\end{eqnarray}
where $\I_{n_u} \in \R^{n_u \times n_u}$ is the identity matrix. This Kronecker product representation~\eqref{eqn:delta_kron} is crucial to enable computationally scalable linear algebra as it allows computation in $\R^{n_\theta}$ to be done via computation in $\R^{n_u}$ and $\R^{n_z}$.

To complete our model discrepancy parameterization, we define an inner product on the affine operator space as the sum of a state space inner product on $\u_0(\t)$ and a linear operator space inner product on $\L(\t)\M_\z$. This leads to endowing $\R^{n_\theta}$ with a $\M_\t$ weighted inner product between two vectors $\t, \vec{\vartheta} \in \R^{n_\theta}$ as 
 \begin{align*}
 (\t,\vec{\vartheta})_{\M_\t} = (\u_0(\t),\u_0(\vec{\vartheta}))_{\M_\u} + (\mathbf{L}(\t) \M_\z ,\mathbf{L}(\vec{\vartheta}) \M_\z)_{\mathcal L(\M_\z,\M_\u)},
 \end{align*}
 where the notation $(\u_1,\u_2)_{\M_\u} = \u_2^T \M_\u \u_1$ is the $\M_\u$ weighted inner product on $\R^{n_u}$ and $(\A_1,\A_2)_{\mathcal L(\M_\z,\M_\u)} = \text{Tr}_{\M_\z}(\A_2^* \A_1)=\text{Tr}_2(\M_\z^{-1} \A_2^T \M_\u \A_1)$ is the Schatten 2-norm on the space of linear operators whose domain and range inner products are weighted by $\M_\z$ and $\M_\u$, respectively. $\A^*$ denotes the adjoint of $\A$, $\text{Tr}_{\M_\z}$ and $\text{Tr}_{2}$ denotes the matrix trace in the $\M_\z$ and $\ell^2$ inner products, respectively. 
 
 It follows that the inner product on $\R^{n_\theta}$ is defined by
$ (\t,\vec{\vartheta})_{\M_\t}  = \vec{\vartheta}^T \M_\t \t$, where
\begin{eqnarray*}
\M_\t =  \left( \begin{array}{cc}
\M_\u & \vec{0}  \\
\vec{0}  & \M_\u \otimes \M_\z
\end{array} \right)  \in \R^{n_\theta \times n_\theta}.
\end{eqnarray*}

\section{Bayesian discrepancy calibration} \label{sec:bayes_calibration}

The discrepancy calibration step seeks to find $\t$ such that $\d(\z_\ell,\t) \approx S(\z_\ell)-\tilde{S}(\z_\ell)$, $\ell=1,2,\dots,N$. For problems with large $n_z \ge N$, there exist many $\t$'s that interpolate the data, as stated in Theorem~\ref{thm:fit_solutions}.
\begin{theorem}
\label{thm:fit_solutions}
Given $N < n_z+1$ data pairs $\{\z_\ell,\vec{d}_\ell\}_{\ell=1}^N$, where $\{\z_\ell\}_{\ell=1}^N$ are linearly independent, there is a vector $\t^\star \in \R^{n_\theta}$ and a subspace $\Theta_{\text{interp}} \subset \R^{n_\theta}$ of dimension $(n_z+1-N)n_u$ such that $\d(\z_\ell,\t^\star+\t)=\vec{d}_\ell$, $\ell=1,2,\dots,N$, for any $\t \in \Theta_{\text{interp}}$.
\end{theorem}
A proof is given in~\ref{proof_thm:fit_solutions}. 

Theorem~\ref{thm:fit_solutions} is an unsurprising result that highlights the necessity of Bayesian calibration to weight the many potential parameters using prior knowledge. To illustrate Theorem~\ref{thm:fit_solutions}, consider an example where $S,\tilde{S}:\R^2 \to \R$ generates data $\vec{d}_1=S(\z_1)-\tilde{S}(\z_1)=1$ and $\vec{d}_2=S(\z_2)-\tilde{S}(\z_2)=2$ for inputs $\z_1=(0,0)$ and $\z_2=(1,0)$. Then 
\begin{align*}
\d(\z,\t) = \theta_0 + \theta_1 z_1 + \theta_2 z_2,
\end{align*}
with parameters $\t=(1,1,\theta_2)$, interpolates the data $\{\z_\ell,\vec{d}_\ell\}_{\ell=1}^2$ for any $\theta_2 \in \R$. The problem is that interpolation does not constrain the variation of $\d$ in directions that are not informed by the data. Specifically, the derivative of $\d$ in the direction $z_2$ is uniformed and hence may be arbitrarily large or small. The prior is therefore critically important to inform the range of variation. Characteristic scales of the problem will be encoded in prior hyper-parameters, which will be analyzed in Section~\ref{sec:hyperparam_selection} where we provide a detailed approach to hyper-parameter specification. 

\subsection{Noise model}
The noise model in Bayesian inverse problems is typically specified based on characteristics of the data collection process. For instance, the noise in observations may be characterized by a known probability distribution. However, in our context, the data is simulation output.  Although we can fit the data perfectly (see Theorem~\ref{thm:fit_solutions}), we are fitting an affine operator $\d(\z,\t)$ to approximate $S(\z)-\tilde{S}(\z)$, which is a nonlinear operator, in general. The noise model reflects our belief about how well the affine map should fit the data. To enable a closed-form expression for posterior samples, we assume an additive Gaussian noise model with mean $\vec{0}$ and covariance $\W_{\vec{d}}^{-1} = \alpha_{\vec{d}} \I_N \otimes \M_\u^{-1}$, where $\alpha_{\vec{d}} > 0$. In other words, the inverse of the noise covariance defines a data misfit norm in the state space. We interpret $\alpha_{\vec{d}}$ as a threshold for how well the calibrated discrepancy should fit the data. As a default, we suggest taking $\alpha_{\vec{d}} = (0.001 C_\d)^2$, where $C_\d \in (0,\infty)$ is defined as the average magnitude of the discrepancy data. This can be interpreted as a noise standard deviation that is $0.1\%$ of the average data magnitude, thus ensuring that we fit the data well. 

\subsection{Prior distribution}
Bayesian priors on function spaces has been the topic of considerable research~\cite{ghattas_infinite_dim_bayes_1,ghattas_infinite_dim_bayes_2,stuart_inv_prob}. However, development for priors defined on spaces of operators is less mature and has only received recent attention with the rise of operator learning~\cite{stuart_learning_lin_op}. 
To leverage existing theory on the specification of function space priors, and to achieve computational efficiency via closed-form expressions for prior and posterior samples, we assume a Gaussian prior for $\t$. To simplify the exposition, we assume that the data is centered by computing the average magnitude of the discrepancy data $\{\vec{d}_\ell\}_{\ell=1}^N$, denote it as $\overline{d} \in \R$, and subtracting $\overline{d}$ element wise from each $\vec{d}_\ell$, $\ell=1,2,\dots,N$. We abuse notation and use $\{\vec{d}_\ell\}_{\ell=1}^N$ to denote the shifted data. 

Since the discrepancy data is centered, we assume a Gaussian prior for $\t$ with mean zero and covariance matrix $\W_\t^{-1} \ \in \R^{n_\theta \times n_\theta}$, which we will specify next. Recall that a mean zero Gaussian random vector has a probability density function proportional to $\exp\left( - \frac{1}{2} \t^T \W_\t \t \right)$. Hence, the inverse covariance matrix, called the precision matrix, defines a weighted inner product such that small values of $\t^T \W_\t \t$ implies high probability. From this perspective, the precision matrix induces a norm penalizing deviation from prior belief.

To define $\W_\t^{-1}$, we introduce mean-zero Gaussian priors for the state and optimization variables, whose covariances are denoted as $\W_\u^{-1} \in \R^{n_u \times n_u}$ and $\W_\z^{-1} \in \R^{n_z \times n_z}$, respectively. We will discuss the specification and interpretation of $\W_\u^{-1}$ and $\W_\z^{-1}$ later. To define a precision matrix for $\t$, we consider that $\d(\z,\t)$ is an affine operator and hence is characterized in terms of a point evaluation and slope (the linear operator $\L(\t) \M_\z$). Since the post-optimality sensitivity analysis is local about $\ztilde$, we localize the weighted norm via
\begin{align}
\label{eqn:discrepancy_norm}
\vert \vert \d(\ztilde,\t) \vert \vert_{\W_\u}^2 + \vert \vert \mathbf{L}(\t) \M_\z \vert \vert_{\mathcal L(\M_\z \W_\z^{-1} \M_\z,\W_\u)}^2 ,
\end{align}
where the Schatten 2-norm of $\mathbf{L}(\t) \M_\z$ has its domain and range weighted by $\M_\z \W_\z^{-1} \M_\z$ and $\W_\u$, respectively. The range's inner product arises naturally as the precision matrix for the state space prior. The domain's inner product may seem unintuitive. This choice of inner product ensures that prior discrepancy samples have adequate smoothness properties, which is necessary to ensure that the prior is well defined in the mesh refinement limit.

Thanks to our choices of the model discrepancy parameterization and inner products, we arrive at the expression
\begin{align*}
\vert \vert \d(\ztilde,\t) \vert \vert_{\W_\u}^2 + \vert \vert \mathbf{L}(\t) \M_\z \vert \vert_{\mathcal L(\M_\z \W_\z^{-1} \M_\z,\W_\u)}^2 = \t^T \W_\t \t,
\end{align*}
where
\begin{eqnarray}
\label{eqn:W_theta}
\W_\t =  \left( \begin{array}{cc}
\W_\u & \W_\u \otimes \ztilde^T \M_\z \\
\W_\u \otimes \M_\z \ztilde & \W_\u \otimes (\W_\z+\M_\z \ztilde  \ztilde^T \M_\z)
\end{array} \right)  \in \R^{n_\theta \times n_\theta}.
\end{eqnarray} 
The Kronecker products in~\eqref{eqn:W_theta} enable computation whose complexity scales like $\mathcal O(n_u + n_z)$ rather than $\mathcal O(n_u n_z)$.

It follows that the covariance matrix is given by
\begin{eqnarray*}
\W_\t^{-1} =  \left( \begin{array}{cc}
(1+\ztilde^T \M_\z \W_\z^{-1} \M_\z \ztilde)\W_\u^{-1} & \W_\u^{-1} \otimes -\ztilde^T \M_\z \W_\z^{-1}  \\
\W_\u^{-1} \otimes - \W_\z^{-1} \M_\z \ztilde & \W_\u^{-1} \otimes  \W_\z^{-1}
\end{array} \right)  \in \R^{n_\theta \times n_\theta}.
\end{eqnarray*} 

 \subsection{Trace class covariance operators}
 
 Since $\U$ is infinite-dimensional and is discretized with the finite-dimensional subspace $\U_h$, $\W_\u^{-1}$ should be defined as the discretization of a trace class covariance operator defined on $\U$. That is, in the mesh refinement limit, $\W_\u^{-1}$ converges to a trace class covariance operator, and consequently, $\text{Tr}_{\M_\u}(\W_\u^{-1})$, the $\M_\u$ weighted trace, converges to a real number. Similarly, when $\Z$ is infinite-dimensional, $\W_\z^{-1}$ should be the discretization of a trace class covariance operator so that $\text{Tr}_{\M_\z}(\W_\z^{-1})$ converges to a real number in the limit as the discretization is refined. In practice, such trace class covariance operators are typically defined as the squared inverse of Laplacian-like differential operators~\cite{stuart_inv_prob}. To this end, let $\E_\u \in \R^{n_u \times n_u}$ and $\E_\z \in \R^{n_z \times n_z}$ denote the discretization of Laplacian-like differential operators and define the covariances 
 \begin{align} \label{eqn:Wu_Wz_def}
 \W_\u^{-1} = \alpha_\u \E_\u^{-1} \M_\u \E_\u^{-1} \qquad \text{and} \qquad \W_\z^{-1} = \alpha_\z \E_\z^{-1} \M_\z \E_\z^{-1},
 \end{align}
  where $\alpha_\u,\alpha_\z >0$ are variance hyper-parameters. Defining the covariances in this way ensures that they are trace class and that samples may be computed efficiently. Theorem~\ref{thm:trace_class_cov} demonstrates that $\W_\t^{-1}$ is trace class in the discretization refinement limit, i.e., it is the discretization of a trace class covariance operator.
 
 \begin{theorem}
 \label{thm:trace_class_cov}
 If
 \begin{align*}
\lim\limits_{\U_h \to \U} \text{Tr}_{\M_\u}(\W_\u^{-1}) \qquad \text{and} \qquad \lim\limits_{\Z_h \to \Z} \text{Tr}_{\M_\z}(\W_\z^{-1})
\end{align*}
exist and are finite, where the limit $\U_h \to \U, \Z_h \to \Z$ indicates refinement of the discretization, and
 \begin{align*}
\lim\limits_{\U_h \to \U, \Z_h \to \Z} \ztilde = \tilde{z} \in \Z ,
\end{align*}
 that is, the optimization solution has a well defined discretization refinement limit, then
\begin{align*}
\lim\limits_{\U_h \to \U, \Z_h \to \Z} \text{Tr}_{\M_\t}(\W_\t^{-1})
\end{align*}
exists and is finite.
 \end{theorem}
 A proof is given in~\ref{proof_thm:trace_class_cov}.

The specific form of the operators $\E_\u$ and $\E_\z$ depends on the problem, most notably, if $\u$ and $\z$ are stationary or transient. We consider the definition of $\E_\u$ for these two scenarios and note that $\E_\z$ may be defined in an analogous fashion. If $\Z$ is finite-dimensional, then $\E_\z=\M_\z$ is a reasonable default which leads to $\W_\z=\alpha_\z^{-1} \M_\z$. We focus on a scalar-valued state noting that analogous derivations may be done component-wise for a vector-valued state. 

\subsubsection*{Stationary problems}
For stationary problems, the state space $\U$ contains functions $f:\Omega \to \R$, where $\Omega \subset \R^s$, $s=1,2,$ or $3$. In this case, a Laplacian-like operator is given by $\E_\u = \beta_\u \K_\u + \M_\u$, where $\beta_\u >0$ and $\K_\u \in \R^{n_u \times n_u}$ is the stiffness matrix\footnote{The $(i,j)$ entry of the stiffness matrix given by the inner product of the gradients of the $i^{th}$ and $j^{th}$ basis functions.}. This corresponds to the discretization of the elliptic operator $-\beta_\u \Delta + I$, where $\Delta$ and $I$ are the Laplacian and identity operators, respectively. 

\subsubsection*{Transient problems}
For transient problems, the state space $\U$ contains functions $f:\Omega \times [0,T] \to \R$, where $\Omega \subset \R^s$, $s=1,2,$ or $3$, and $T>0$ is the final time. Assume that time is discretized at $n_t$ nodes $0=t_1 < t_2 < \cdots < t_{n_t}=T$ and that space is discretized with $n_s$ basis functions. Consequently, the state coordinate vector is $\u \in \R^{n_u}$, $n_u = n_t n_s$, with an ordering $\u=\begin{pmatrix} \u_1^T & \u_2^T & \cdots & \u_{n_t}^T \end{pmatrix}^T$, where $\u_i \in \R^{n_s}$, $i=1,2,\dots,n_t$, are the spatial coordinates at each time step. It follows that the discretized inner product is weighted by $\M_\u = \M_t \otimes \M_s$, where $\M_t \in \R^{n_t \times n_t}$ is the mass matrix corresponding to piecewise linear basis functions (i.e., hat functions) defined on the time nodes and $\M_s$ is the mass matrix corresponding to the spatial basis functions defined on $\Omega$. 

The state covariance matrix is defined by $\W_\u^{-1} = \alpha_\u \W_t^{-1} \otimes \W_s^{-1}$, where $\W_t^{-1}$ and $\W_s^{-1}$ are covariance matrices corresponding to the temporal and spatial domains, respectively. As in the stationary case, we model $\W_s^{-1}= \E_s^{-1} \M_s \E_s^{-1}$, where $\E_s = \beta_\u \K_s + \M_s$, $\K_s$ is the stiffness matrix corresponding to the spatial discretization. To define $\W_t^{-1}$, we consider the elliptic operator $\E_t =\beta_t \K_t + \M_t$ on the temporal domain, where $\beta_t>0$ and $\K_t$ is the stiffness matrix corresponding to the time discretization. We defined $\W_s^{-1}$ as square of $\E_s^{-1}$ to ensure that the covariance is trace class for spatial dimensions $s=2$ and $s=3$. However, squaring is not required in the temporal domain as $\E_t^{-1}$ is a trace class operator since it is defined on an interval~\cite{stuart_inv_prob}. In many applications, the model discrepancy is small, or identically zero, at time $t=0$, and its magnitude increases as a function of time due to the compounding of errors. To model this, we consider a time dependent variance weighting $\vec{\alpha}_t \in \R^{n_t}$ and define $\W_t^{-1} = \mathbf{D}(\vec{\alpha}_t)^\frac{1}{2} \E_t^{-1}\mathbf{D}(\vec{\alpha}_t)^\frac{1}{2}$, where $\mathbf{D}(\vec{\alpha}_t) \in \R^{n_t \times n_t}$ is the diagonal matrix whose diagonal is $\vec{\alpha}_t$. This provides the modeling flexibility to capture temporal trends in the discrepancy while preserving the computational efficiency of the space-time Kronecker product structure and elliptic solves. Note that the space-time Kronecker product $\W_\u^{-1}$ can still be written in the form of~\eqref{eqn:Wu_Wz_def} with $\E_\u$ corresponding to a Kronecker product of time and space operators.

\subsection{Computing prior discrepancy samples}
We seek to determine $\T_\t \in \R^{n_\theta \times n_\theta}$ such that $\W_\t^{-1} = \T_\t \T_\t^T$, and compute prior samples as $\T_\t \vec{\omega}_\t$, where $\vec{\omega}_\t \in \R^{n_\theta}$ is sampled from a standard normal distribution on $\R^{n_\theta}$. To this end, let $\W_\u^{-1} = \T_\u \T_\u^T$ and $\W_\z^{-1} = \T_\z \T_\z^T$ be factorizations and assume that matrix-vector products may be computed with $\T_\u$ and $\T_\z$. For instance, matrix-vector products with $\T_\u = \sqrt{\alpha_\u} \E_\u^{-1} \M_\u^{\frac{1}{2}}$ may be computed using Krylov subspace methods that require matrix-vector products $\E_\u \v$ and $\M_\u \v$, respectively. It follows that we may decompose $\W_\t^{-1} = \T_\t \T_\t^T$ with
\begin{eqnarray*}
\T_\t =  \left( \begin{array}{cc}
\T_\u & -\T_\u \otimes \ztilde^T \M_\z \T_\z \\
\vec{0} & \T_\u \otimes \T_\z
\end{array} \right)  \in \R^{n_\theta \times n_\theta}.
\end{eqnarray*} 

To analyze samples from the prior discrepancy, let $\vec{\omega}_\t \in \R^{n_\theta}$ be sampled from a standard normal distribution on $\R^{n_\theta}$ and consider 
\begin{align}
\label{eqn:prior_discrepancy_samples}
\d(\z,\T_\t \vec{\omega}_\t) & = 
\left( \begin{array}{cc}
\I_{n_u} & \I_{n_u} \otimes \z^T \M_\z
\end{array} \right) \T_\t \vec{\omega}_\t \\
& = 
\left( \begin{array}{cc}
\T_\u & \T_\u \otimes (\z-\ztilde)^T \M_\z \T_\z
\end{array} \right) \vec{\omega}_\t \nonumber \\
& = \T_\u \u_0(\vec{\omega}_\t) + \T_\u \L(\vec{\omega}_\t) \T_\z^T \M_\z (\z-\ztilde), \nonumber
\end{align}
where $\u_0$ and $\L$ are defined in~\eqref{eqn:u0_L}. Plugging in $\T_\u = \sqrt{\alpha_\u} \E_\u^{-1} \M_\u^{\frac{1}{2}}$ and $\T_\z = \sqrt{\alpha_\z} \E_\z^{-1} \M_\z^{\frac{1}{2}}$, and exploiting properties of Gaussian distributions, we can equivalently express the samples as
\begin{align}
\label{eqn:prior_discrepancy_samples_2}
\d(\z,\T_\t \vec{\omega}_\t) = \sqrt{\alpha_\u} \E_\u^{-1} \M_\u^{\frac{1}{2}} \vec{\omega}_\u^{(1)} + \sqrt{\alpha_\u} \sqrt{\alpha_\z} \sqrt{s_\z} \E_\u^{-1} \M_\u^{\frac{1}{2}} \vec{\omega}_\u^{(2)} 
\end{align}
where $s_\z = (\z-\ztilde)^T \M_\z \E_\z^{-1} \M_\z \E_\z^{-1} \M_\z (\z-\ztilde)$ and $\vec{\omega}_\u^{(1)},\vec{\omega}_\u^{(2)}$ are independent samples from a standard normal distribution on $\R^{n_u}$. In what follows, we use~\eqref{eqn:prior_discrepancy_samples} and~\eqref{eqn:prior_discrepancy_samples_2} to analyze the prior model discrepancy and guide the specification of its hyper-parameters.

 \section{hyper-parameter specification} \label{sec:hyperparam_selection}
There are four prior hyper-parameters, $\alpha_\u,\beta_\u,\alpha_\z$, and $\beta_\z$, for stationary problems, and two additional hyper-parameters, $\vec{\alpha}_t$ and $\beta_t$, for transient problems. If $Z$ is finite-dimensional, we can assume $\beta_\z=0$. For generality, we will assume that $Z$ is infinite-dimensional and include $\beta_\z$ in our subsequent analysis. Although these hyper-parameters have physical interpretations to guide their specification, determining an appropriate value for them is nontrivial since they are embedded in blocks of a larger covariance that is defined on the space of operators. In this section, we present an approach that uses the data $\{ \z_\ell,\vec{d}_\ell \}_{\ell=1}^N$ to initialize the hyper-parameters and introduce a prior sample visualization framework. Visualization of the prior discrepancy is crucial to verify the hyper-parameter values generated by the algorithm-based initialization and adjust them as necessary. We suggest this algorithm-based initialization and prior sample visualization approach rather than a hierarchical Bayesian approach because the calibration data is limited. 

\subsection{Smoothness hyper-parameters}
The hyper-parameters, $\beta_\z,\beta_\u$, and $\beta_t$ (for transient problems) scale the Laplacian in the optimization variable and state prior covariances, respectively, and consequently determine the smoothness of samples from the prior discrepancy. Since we are assuming a Gaussian prior, these smoothness hyper-parameters are related to the correlation length, defined as the distance at which we observe an empirical correlation of $0.1$ in a Gaussian random field~\cite{lindgren11,matern_cov_villa}. The values of $\beta_\u$ and $\beta_t$ should be specified to reflect the smoothness characteristics of the $S(\z)-\tilde{S}(\z)$. In some applications, $\beta_\u$ and $\beta_t$ may be selected based on physical knowledge such as knowing the length scales that are resolved by the high- and low-fidelity models. However, this information is already present in $\{\vec{d}_\ell\}_{\ell=1}^N$. Accordingly, we estimate the spatial and temporal correlation lengths of $\vec{d}_\ell=S(\z_\ell)-\tilde{S}(\z_\ell)$, $\ell=1,2,\dots,N$, as described in~\ref{appendix_corr_length}. For stationary problems, $\beta_t$ is ignored. For transient problems, we separate $\{ \vec{d}_\ell \}_{\ell=1}^N$ into space and time components to calculate the spatial and temporal correlation lengths, $\kappa_\u$ and $\kappa_t$, separately. Specifically, each time snapshot is considered separately to compute $\kappa_\u^i$, $i=1,2,\dots,n_tN$, and define $\kappa_\u=\frac{1}{n_tN} \sum_{i=1}^{n_tN} \kappa_\u^i.$ Similarly, the time series for each spatial node is considered separately to compute $\kappa_t^j$, $j=1,2,\dots,n_sN$, and define $\kappa_t=\frac{1}{n_sN} \sum_{j=1}^{n_sN} \kappa_t^j.$ The smoothness hyper-parameters are specified as $\beta_\u = \kappa_\u^2/C_s$ and $\beta_t = \kappa_t^2/4$~\cite{matern_cov_villa}, where $C_s$ depends on the spatial dimension $s$ of functions in $\U$. Specifically, $C_1=12$, $C_2=8$, and $C_3=4$. We compute $\beta_\z$ in an analogous manner using the data $\z_\ell$, $\ell=1,2,\dots,N$. Note that $\beta_t = \kappa_t^2/4 \ne  \kappa_t^2/C_1$ because the temporal operator is not squared.

\subsection{Temporal variance weighting hyper-parameter} \label{ssec:time_weighting}
The temporal trend of the discrepancy magnitude is controlled by the hyper-parameter $\vec{\alpha}_t \in \R^{n_t}$. Note that if $\vec{\alpha}_t \in \R^{n_t}$ is scaled by a constant $C$, then $\alpha_\u$ may be scaled by $C^{-1}$ to render the same prior covariance. Hence, we are not concerned about the norm of the vector $\vec{\alpha}_t$ but rather the magnitude of its entries relative to one another. Note that all entries must be strictly positive since they specify variance. 

Defining $\vec{\alpha}_t = \vec{e}$, the vector of $1$'s, will treat all time nodes equally (up to some boundary effects of the elliptic operator). However, the discrepancy data $\{ \vec{d}_\ell \}_{\ell=1}^N$ provides temporal trends to inform the specification of $\vec{\alpha}_t$. To this end, we compute the spatial norm (squared) of $\{ \vec{d}_\ell \}_{\ell=1}^N \subset \R^{n_t n_s}$ for each time step, yielding a dataset $\{ \hat{\vec{d}}_\ell \}_{\ell=1}^N \subset \R^{n_t}$, where the $i^{th}$ entry of $\hat{\vec{d}}_\ell$ is $\| \vec{d}_{\ell,i} \|_{\M_s}^2 = \vec{d}_{\ell,i}^T \M_s \vec{d}_{\ell,i}$, $\vec{d}_{\ell,i} \in \R^{n_s}$ denotes the time node $i$ snapshot of $\vec{d}_\ell$. Note that we compute the norm squared since the variance is a squared quantity. Averaging the data and normalizing the average by its maximum (the vector infinity norm $\| \cdot \|_{\infty}$) yields 
\begin{align*}
\hat{\alpha}_t = \frac{\frac{1}{N} \sum\limits_{\ell=1}^N  \hat{\vec{d}}_\ell}{ \left\| \frac{1}{N} \sum\limits_{\ell=1}^N  \hat{\vec{d}}_\ell \right\|_{\infty}} \in \R^{n_t},
\end{align*}
which weights each time step with a value in $[0,1]$. To ensure sufficient expressivity of the prior, we inflate the variance by adding a constant $\epsilon_t >0$ to define $\vec{\alpha}_t = \hat{\vec{\alpha}}_t + \epsilon_t \vec{e}$. We suggest $\epsilon_t = 0.01$ as a default meaning that the variance is inflated by $1\%$ of the maximum variance magnitude. 

\subsection{State variance hyper-parameter}
The prior discrepancy, evaluated at $\z=\ztilde$, is normally distributed with mean zero and covariance $\W_\u^{-1}$, see~\eqref{eqn:prior_discrepancy_samples}. Hence, the variance hyper-parameter $\alpha_\u$, which scales $\W_\u^{-1}$, determines the magnitude of the prior model discrepancy at $\z=\ztilde$. We initialize $\alpha_\u$ so that the expectation of the prior discrepancy norm, evaluated at $\ztilde$, equals the observed magnitude of the discrepancy data. Specifically, consider a prior sample $\T_\t \vec{\omega}_\t$, where $\vec{\omega}_\t$ is sampled from a standard normal distribution on $\R^{n_\theta}$. The discrepancy norm squared is 
\begin{align*}
\| \d(\ztilde,\T_\t \vec{\omega}_\t) \|_{\M_\u}^2 = \u_0(\vec{\omega}_\t)^T \T_\u^T \M_\u \T_\u \u_0(\vec{\omega}_\t),
\end{align*}
which is distributed as a weighted sum of chi-squared random variables, where the weights are defined by the eigenvalues of $\T_\u^T \M_\u \T_\u$. It follows that
\begin{align*}
\mathbb{E}_{\vec{\omega_\t}} \left[ \| \d(\ztilde,\T_\t \vec{\omega}_\t) \|_{\M_\u}^2 \right] & = \mathbb{E}_{\vec{\omega_\t}} \left[ \u_0(\vec{\omega}_\t)^T \T_\u^T \M_\u \T_\u \u_0(\vec{\omega}_\t) \right] \\
& = \alpha_\u \text{Tr}_2( \M_\u^{\frac{1}{2}} \E_\u^{-1} \M_\u \E_\u^{-1} \M_\u^{\frac{1}{2}} ). \\
& = \alpha_\u \text{Tr}_{\M_\u}( \E_\u^{-1} \M_\u \E_\u^{-1}  )
\end{align*}
Hence, we initialize the state variance hyper-parameter as
\begin{align*}
\alpha_\u = \frac{ \| \vec{d}_1 \|_{\M_\u}^2 }{\text{Tr}_{\M_\u}( \E_\u^{-1} \M_\u \E_\u^{-1}  )},
\end{align*}
so that $\mathbb{E}_{\vec{\omega_\t}} \left[ \| \d(\ztilde,\T_\t \vec{\omega}_\t) \|_{\M_\u}^2 \right]=\| \vec{d}_1 \|_{\M_\u}^2$. Note that the initialization of $\alpha_\u$ depends only on $\vec{d}_1$ and not $\{\vec{d}_\ell \}_{\ell=1}^N$ since the analysis is based on evaluating the discrepancy at $\z=\ztilde=\z_1$.

\subsection{Optimization variable variance hyper-parameter}

The variance hyper-parameter $\alpha_\z$ determines the rate at which the prior discrepancy varies with respect to changes in $\z$, i.e., it controls the magnitude of the Jacobian of the discrepancy with respect to $\z$. To specify $\alpha_\z$, consider a perturbation $\Delta \z \in \R^{n_z}$ and examine
\begin{align*}
\d(\ztilde + \Delta \z,\vec{T}_\t \vec{\omega}_\t) - \d(\ztilde ,\vec{T}_\t \vec{\omega}_\t) = \T_\u \L(\vec{\omega}_\t) \T_\z^T \M_\z \Delta \z ,
\end{align*}
which follows from~\eqref{eqn:prior_discrepancy_samples}. Theorem~\ref{thm:exp_delta_diff} provides a convenient expression for the expectation (over the prior distribution) of the discrepancy difference, normalized by the expected discrepancy at the low-fidelity optimization solution $\ztilde$.
\begin{theorem} \label{thm:exp_delta_diff}
If the covariances are defined as in~\eqref{eqn:Wu_Wz_def}, then
\begin{align} \label{eqn:exp_delta_diff_1_pert}
\frac{\mathbb{E}_{\vec{\omega_\t}} \left[ \| \d(\ztilde + \Delta \z,\L_\t \vec{\omega}_\t) - \d(\ztilde ,\L_\t \vec{\omega}_\t) \|_{\M_\u}^2 \right] }{\mathbb{E}_{\vec{\omega_\t}} \left[ \| \d(\ztilde ,\L_\t \vec{\omega}_\t) \|_{\M_\u}^2 \right] } = \alpha_\z \Delta \z^T \M_z \E_\z^{-1} \M_\z \E_\z^{-1} \M_\z \Delta \z .
\end{align}
\end{theorem}
A proof is given in~\ref{proof_thm:exp_delta_diff}.

To determine $\alpha_\z$ that is appropriate for a range of $\Delta \z$'s, we consider samples from a mean-zero Gaussian distribution with covariance $\E_\z^{-1} \M_\z \E_\z^{-1}$ and rescaled to have norm $\| \ztilde \|_{\M_\z}$. The samples take the form
\begin{align*}
\Delta \z = \| \ztilde \|_{\M_\z} \frac{\E_\z^{-1} \M_\z^\frac{1}{2} \vec{\omega}_\z}{ \| \E_\z^{-1} \M_\z^\frac{1}{2} \vec{\omega}_\z \|_{\M_\z}},
\end{align*}
where $\vec{\omega}_\z$ is sampled from standard normal distribution on $\R^{n_z}$. We would like the prior model discrepancy to satisfy 
\begin{align}
\label{eqn:gamma_criteria}
& \mathbb{E}_{\Delta \z} \left [\mathbb{E}_{\vec{\omega_\t}} \left[ \| \d(\ztilde + \Delta \z,\L_\t \vec{\omega}_\t) - \d(\ztilde ,\L_\t \vec{\omega}_\t) \|_{\M_\u}^2 \right]  \right] \nonumber \\= 
&   \mathbb{E}_{\Delta \z} \left [ \| (S(\ztilde+\Delta \z)-\tilde{S}(\ztilde+\Delta \z)) - (S(\ztilde)-\tilde{S}(\ztilde) ) \|_{\M_\u}^2 \right].
\end{align}
That is, the average (over $\Delta \z$) difference of $\d$, averaged over the prior, should equal the average (over $\Delta \z$) difference of $S-\tilde{S}$. Letting
\begin{align}
\label{eqn:gamma_def}
\gamma^2 =  \mathbb{E}_{\Delta \z} \left [ \| (S(\ztilde+\Delta \z)-\tilde{S}(\ztilde+\Delta \z)) - (S(\ztilde)-\tilde{S}(\ztilde) ) \|_{\M_\u}^2 \right]
\end{align}
denote the desired discrepancy difference from~\eqref{eqn:gamma_criteria}, recalling that $\mathbb{E}_{\vec{\omega_\t}} \left[ \| \d(\ztilde ,\L_\t \vec{\omega}_\t) \|_{\M_\u}^2 \right] =\| \vec{d}_1\|_{\M_\u}^2$ by our specification of $\alpha_\u$, and manipulating~\eqref{eqn:exp_delta_diff_1_pert}, yields the expression to initialize $\alpha_\z$ as
\begin{align}
\label{eqn:alpha_z}
\alpha_\z  = \frac{\gamma^2}{\| \vec{d}_1\|_{\M_\u}^2} \frac{1}{\mathbb{E}_{\Delta \z} \left[ \Delta \z^T \M_z \E_\z^{-1} \M_\z \E_\z^{-1} \M_\z \Delta \z \right]}.
\end{align}
Next we discuss estimating the numerator and denominator of~\eqref{eqn:alpha_z}, thus providing an approximate initialization of $\alpha_\z$. Exact computation of~\eqref{eqn:alpha_z} is not possible since $\gamma^2$ involves the high-fidelity model, which can be accessed only sparingly. 

\subsubsection*{Estimation of $\gamma^2$}
We leverage evaluations of $\tilde{S}$ to capture characteristic scales that are relevant to provide an estimate of $\gamma^2$. For notational simplicity, let $\Delta S(\Delta \z) = S(\ztilde+\Delta \z)-S(\ztilde)$ and $\Delta \tilde{S}(\Delta \z) =  \tilde{S}(\ztilde+\Delta \z) -  \tilde{S}(\ztilde)$ and observe that manipulations of~\eqref{eqn:gamma_def} yields
\begin{align*}
\gamma^2 & = \mathbb{E}_{\Delta \z} \left [( \Delta S(\Delta \z) - \Delta \tilde{S}(\Delta \z) )^T \M_\u ( \Delta S(\Delta \z) - \Delta \tilde{S}(\Delta \z) ) \right] \\
& = \mathbb{E}_{\Delta \z} \left [ \| \Delta S(\Delta \z) \|_{\M_\u}^2 + \| \Delta \tilde{S}(\Delta \z) \|_{\M_\u}^2 - 2 \Delta S(\Delta \z)^T \M_\u \Delta \tilde{S}(\Delta \z) \right] .
\end{align*}
To facilitate a practical estimator for $\gamma^2$, we assume that $\mathbb E_{\Delta \z} \left[ \| \Delta S(\Delta \z) \|_{\M_\u}^2 \right] = \mathbb E_{\Delta \z} \left[ \| \Delta \tilde{S}(\Delta \z) \|_{\M_\u}^2 \right]$, that is, the high- and low-fidelity models have the same magnitude of sensitivity with respect to $\z$, on average. Then we have
\begin{align}
\label{eqn:gamma_eqn}
\gamma^2  = 2 \mathbb{E}_{\Delta \z} \left [  \| \Delta \tilde{S}(\Delta \z) \|_{\M_\u}^2 - \cos(\zeta)\| \Delta S(\Delta \z) \|_{\M_\u} \| \Delta \tilde{S}(\Delta \z) \|_{\M_\u}  \right] ,
\end{align}
where $\zeta$ is the angle (in the $\M_\u$ inner product) between $\Delta S(\Delta \z)$ and $\Delta \tilde{S}(\Delta \z)$. We cannot use~\eqref{eqn:gamma_eqn} to compute $\gamma^2$ because $\zeta$ is an unknown that depends on the high-fidelity model $S$ and varies with the input $\Delta \z$. However, $\cos(\zeta) \in (-1,1)$ with $\cos(\zeta)=1$ and $\cos(\zeta)=0$ corresponding to the cases when $\Delta S(\Delta \z)$ and $\Delta \tilde{S}(\Delta \z)$ are aligned and orthogonal, respectively. Taking an intermediate value of $\cos(\zeta)=\frac{1}{2}$, and assuming that 
$$\mathbb E_{\Delta \z} \left[ \| \Delta S(\Delta \z) \|_{\M_\u} \|  \Delta \tilde{S}(\Delta \z) \|_{\M_\u}  \right] =\mathbb E_{\Delta \z} \left[ \| \Delta S(\Delta \z) \|_{\M_\u} \right] \mathbb E_{\Delta \z} \left[ \| \Delta \tilde{S}(\Delta \z) \|_{\M_\u} \right],$$ yields the approximation 
\begin{align}
\label{eqn:gamma_approx}
\gamma^2 \approx \mathbb{E}_{\Delta \z} \left [ \| \Delta \tilde{S}(\Delta \z) \|_{\M_\u}^2 \right].
\end{align}
Then $\gamma^2$ may be approximated via~\eqref{eqn:gamma_approx} using only low-fidelity model evaluations.

\subsubsection*{Estimation of $\mathbb{E}_{\Delta \z} \left[ \Delta \z^T \M_z \E_\z^{-1} \M_\z \E_\z^{-1} \M_\z \Delta \z \right]$}

To estimate the denominator of~\eqref{eqn:alpha_z}, Theorem~\ref{thm:exp_delta_z} demonstrates that the expectation may be expressed in a convenient form by leveraging the spectral decomposition of $\E_\z$.
\begin{theorem} \label{thm:exp_delta_z}
Let 
\begin{align}
\label{eqn:E_z_eig}
\E_\z = \M_\z \vec{V} \vec{\Lambda} \vec{V}^T \M_\z
\end{align}
denote the generalized eigenvalues decomposition of $\E_\z$ in the $\M_\z$ inner product, where the eigenvectors satisfies $\vec{V}^T \M_\z \vec{V}=\vec{I}$, $\vec{V} \in \R^{n_z \times n_z}$. Then
\begin{align*}
\mathbb E_{\Delta \z} \left[ \Delta \z^T \M_z \E_\z^{-1} \M_\z \E_\z^{-1} \M_\z \Delta \z \right] = \| \ztilde \|_{\M_\z}^2 \mathbb E_{\vec{\omega}_\z} \left[ \frac{ \vec{\omega}_\z^T \vec{\Lambda}^{-4} \vec{\omega}_\z} { \vec{\omega}_\z^T  \vec{\Lambda}^{-2} \vec{\omega}_\z} \right],
\end{align*}
where $\vec{\omega}_\z$ follows a standard normal distribution on $\R^{n_z}$.
\end{theorem}
A proof is given in~\ref{proof_thm:exp_delta_z}

Recall that $\E_\z = \beta_\z \K_\z + \M_\z$, so 
\begin{align*}
\beta_\z \v_i^T \K_\z \v_i + \v_i^T \M_\z \v_i = \lambda_i \v_i^T \M_\z \v_i \implies \lambda_i = 1+\beta_\z \v_i^T \K_\z \v_i  \ge 1.
\end{align*}
Furthermore, the eigenvalues are known to increase quadratically, i.e., $\lambda_i = \mathcal O(i^2)$, since $\K_\z$ is the discretization of the Laplacian. This implies that the integrand of
\begin{align}
\label{eqn:expectation_of_ratio}
\mathbb E_{\vec{\omega}_\z} \left[ \frac{ \vec{\omega}_\z^T \vec{\Lambda}^{-4} \vec{\omega}_\z} { \vec{\omega}_\z^T  \vec{\Lambda}^{-2} \vec{\omega}_\z} \right]
\end{align}
is strictly less than $1$ and can be approximated by truncating the eigenvalue contributions since the numerator and denominator scale like $\lambda_i^{-4}=\mathcal O(i^{-8})$ and $\lambda_i^{-2}=\mathcal O(i^{-4})$. Hence, we compute~\eqref{eqn:expectation_of_ratio} via Monte Carlo estimation, which requires evaluating the leading eigenvalues of $\E_\z^{-1}$ and generating standard normal random vectors whose dimension corresponds to the number of eigenvalues that are retained. For some problems, for instance those on a rectangular domain, the eigenvalues admit an analytic expression. Consequently, the Monte Carlo estimation may be computed at a negligible cost. 

\begin{figure}[h]
\centering
  \includegraphics[width=0.35\textwidth]{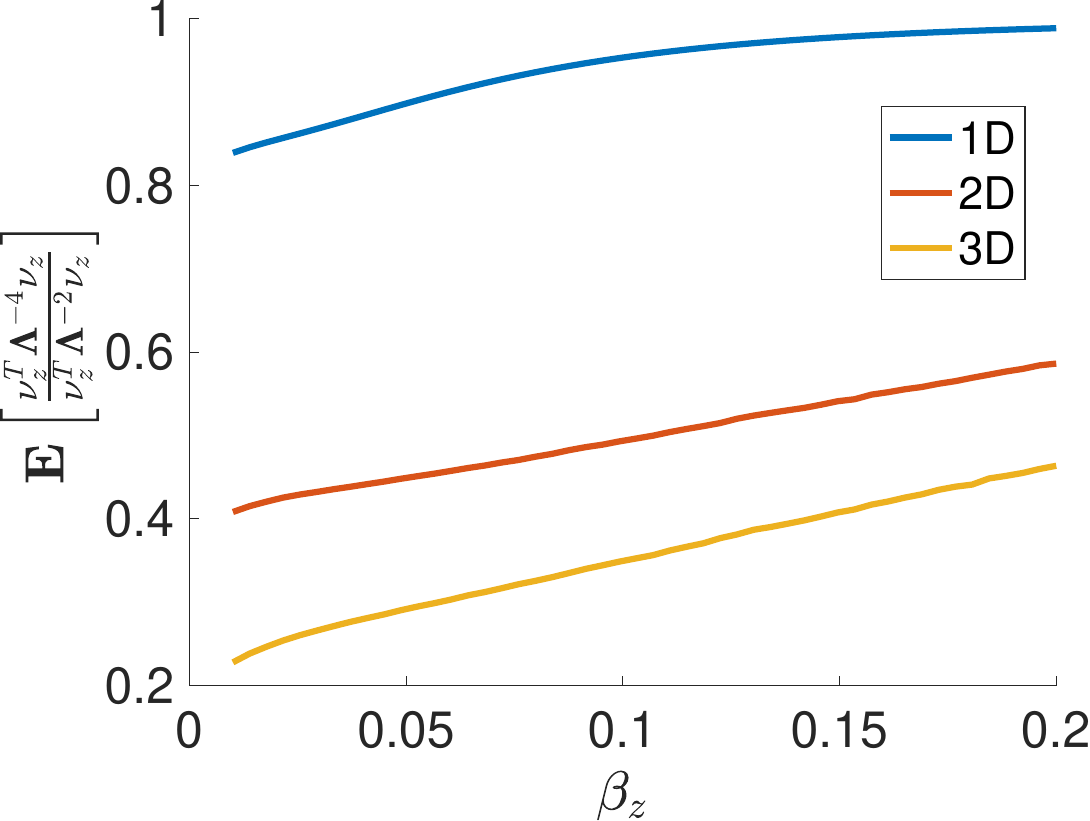}
    \caption{Value of~\eqref{eqn:expectation_of_ratio} for various smoothness hyper-parameters $\beta_\z$ and spatial dimensions $s=1,2,3$.}
  \label{fig:mc_eig_computation}
\end{figure}

Figure~\ref{fig:mc_eig_computation} displays~\eqref{eqn:expectation_of_ratio} as a function of $\beta_\z$ for elliptic operators defined on the unit cube $[0,1]^s$ in spatial dimensions $s=1,2,$ and $3$. We observe that~\eqref{eqn:expectation_of_ratio} assumes values within the interval $(.2,1)$, that it varies smoothly and monotonically with respect to $\beta_\z$, and that it takes smaller values for higher spatial dimensions. These advantageous properties ensure that we can accurately estimate~\eqref{eqn:expectation_of_ratio} at a modest computational cost.

\subsection{The role of visualization to verify or adjust the hyper-parameters} \label{ssec:viz_framework}
The algorithms presented above serve to initialize the prior hyper-parameters. However, the algorithms rely on many assumptions and approximations. Typically, $\alpha_\u$ and $\alpha_t$ can be initialized well since they do not require strong assumptions. On the other hand, $\beta_\u, \beta_\z,$ and $\beta_t$ rely on empirical correlation length estimation, which assumes that the data is sampled from a Gaussian random field. The most challenging hyper-parameter to initialize is $\alpha_\z$, which involves assumptions on the relationship between the high- and low-fidelity models and a crude estimation of $\gamma^2$. 

Misspecification of the prior hyper-parameters will lead to inappropriate prior distributions that are either overly restrictive (i.e., only permitting small discrepancies) or highly uninformed (i.e., permitting such large discrepancies that the uncertainty analysis is not useful). Visualizing prior discrepancy samples is crucial to ensure that the hyper-parameters produce a physically meaningful prior and consequently a useful quantification of uncertainty. However, visualization of the discrepancy is challenging since it is an operator mapping from $\Z$ to $\U$, which in general, are both function spaces. In the subsections below, we present an approach to curate prior sample data and visualize it in a way that guides verification and/or modification of the hyper-parameters determined by the initialization algorithms. 

\subsubsection*{Prior sample data}
 Let $\{\t_i\}_{i=1}^Q$ denote samples from the prior parameter distribution. For each $\t_i$, $\d(\z,\t_i)$ is an operator mapping from $\Z$ to $\U$ and hence can be visualized only as input/output pairs since each $\z$ and $\d(\z,\t_i)$ are functions. We center our analysis around the low-fidelity optimization solution $\ztilde$ and select optimization variable perturbations $\{ \Delta \z_k\}_{k=1}^P$ to form a dataset
\begin{eqnarray} \label{eqn:prior_sample_dataset}
\{ \Delta \z_k, \d(\ztilde,\t_i), \d(\ztilde + \Delta \z_k,\t_i) \}_{i=1,k=1}^{Q,P}.
\end{eqnarray}

The perturbations $\{ \Delta \z_k \}_{k=1}^P$ are determined by seeking directions that yield the largest change in the prior discrepancy $\d$. Specifically, we compute $\Delta \z_k$'s that maximize the expected value of the Jacobian of $\d$ with respect to $\z$, i.e., $\mathbb{E}_\t \left[ \nabla_\z \d(\cdot,\t) \right]$. In~\ref{appendix_z_pert}, we demonstrate that the $\Delta \z_k$'s should be taken as the leading eigenvectors of $\E_\z^{-1}$ (computed in the $\M_\z$ inner product). Observing from~\eqref{eqn:prior_discrepancy_samples_2} that the magnitude of $\d(\ztilde + \Delta \z_k,\t_i)-\d(\ztilde,\t_i)$ is proportional to the eigenvalue associated with the eigenvector $\Delta \z_k$, and recalling that the leading eigenvalue of $\E_\z^{-1}$ equals $1$, we determine the number of perturbations $P$ by computing all eigenpairs whose eigenvalue is $\ge 0.1$, i.e., within one order of magnitude of the leading eigenvalue. 

 Since an adequate number of samples $\{\t_i\}_{i=1}^Q$ is required to characterize the prior distribution and an adequate number of perturbations $P$  is required to characterize the sensitivity of $\d$, we frequently have $Q = \mathcal O(10^2)$ and $P=\mathcal O(10^2)$, or larger, which implies that the prior discrepancy dataset~\eqref{eqn:prior_sample_dataset} consists of $\mathcal O(Q \cdot P)=\mathcal O(10^4)$ functions, each of which vary in space and/or time. Visualization of such data is challenging since only a small number of functions may be visualized simultaneously. An interactive framework is proposed to identify extreme cases in the dataset and render visualizations of select samples and perturbations. 

\subsubsection*{Interactive visualization framework}
To summarize the data in a way that $\mathcal O(10^4)$ samples may be visualized, we define two metrics and project the data onto a scatter plot where each sample is represented by the metric values in $\R^2$. We select metrics that reflect the interpretation of the hyper-parameters and hence enable visualizations that guide their specification. The temporal hyper-parameters $\vec{\alpha}_t$ and $\beta_t$ can be visualized by overlaying time series plots, as demonstrated Section~\ref{ssec:transient_numerical_example} (see Figure~\ref{fig:trans_var}). We focus on two sets of visualizations relevant for two subsets of hyper-parameters, $\{\alpha_\u,\beta_\u\}$ and $\{\alpha_\z,\beta_\z\}$. Visualization for $\{\alpha_\z,\beta_\z\}$ is demonstrated in Section~\ref{ssec:num_2d} (see Figure~\ref{fig:interactive_viz_snapshot}).

Recall from~\eqref{eqn:prior_discrepancy_samples_2} that $\d(\ztilde,\t_i)$ is directly proportional to $\sqrt{\alpha_\u}$ and that its spatial smoothness properties are determined by $\beta_\u$. Inspired by this interpretation of $\alpha_\u$ and $\beta_\u$, we use the maximum absolute value and correlation length of $\d(\ztilde,\t_i)$ (which is a function defined on space and/or time) as metrics to project the data $\{ \d(\ztilde,\t_i) \}_{i=1}^Q$ onto a 2D scatter plot, which allows for an assessment if the magnitude and smoothness of the prior samples $\{ \d(\ztilde,\t_i) \}_{i=1}^Q$ are commensurate to characteristics of $S(\ztilde)-\tilde{S}(\ztilde)$. Visualizing the fields $\d(\ztilde,\t_i)$ for the extreme and average cases in the scatter plot reveals if $\alpha_\u$ and $\beta_\u$ have been specified correctly. Generally, $\beta_\u$ should be analyzed first since $\alpha_\u$ depends on it. Furthermore, there are potentially invalid assumptions made in our initialization of $\beta_\u$, whereas the initialization of $\alpha_\u$ did not make assumptions on the data. If the visualizations identify misspecification of the hyper-parameters, the interpretation provided by~\eqref{eqn:prior_discrepancy_samples_2} guides how the values should be modified. In particular, if the samples $\{ \d(\ztilde,\t_i) \}_{i=1}^Q$ are too smooth (or non-smooth) then $\beta_\u$ should be decreased (or increased) and if their magnitudes are too large (or small) then $\alpha_\u$ should be decreased (or increased).

We assume that $\{\alpha_\u,\beta_\u\}$ have been specified properly based on of $\{ \d(\ztilde,\t_i) \}_{i=1}^Q$, and now consider analysis of $\{\alpha_\z,\beta_\z\}$. Again recalling~\eqref{eqn:prior_discrepancy_samples_2}, observe that
 \begin{align} \label{eqn:delta_diff_pert}
 \d(\ztilde+\Delta \z_k,\t_i)-\d(\ztilde,\t_i) = \sqrt{\alpha_\u} \sqrt{\alpha_\z} \sqrt{s_{\Delta \z_k}} \E_\u^{-1} \M_\u^{\frac{1}{2}} \vec{\omega}_\u^{(2)},
 \end{align}
 where $s_{\Delta \z_k} = \Delta \z_k^T \M_\z \E_\z^{-1} \M_\z \E_\z^{-1} \M_\z \Delta \z_k=\| \E_\z^{-1} \M_\z \Delta \z \|_{\M_\z}^2$. We may interpret $\alpha_\z$ and $\beta_\z$ in~\eqref{eqn:delta_diff_pert} as controlling the rate of change when $\d$ is perturbed in the direction $\Delta \z_k$. Note that $\beta_\z$ is implicit in~\eqref{eqn:delta_diff_pert} since $\beta_\z$ appears in $\E_\z$, which enters~\eqref{eqn:delta_diff_pert} via $s_{\Delta \z_k}$.
 
 Notice that $\beta_\z$ weights perturbations based on their smoothness properties. In particular, the eigenvectors of $\E_\z^{-1}$ behave like sinusoidal functions and the eigenvalues are inversely proportional to the wave number. Hence, smooth (i.e., small wavenumber or low-frequency) eigenvectors correspond to larger eigenvalues. Furthermore, observe that $\d(\ztilde+\Delta \z_k,\t_i) -\d(\ztilde,\t_i)$ is proportional to $\sqrt{s_{\Delta \z_k}}$, which is proportional to the eigenvalue of $\E_\z^{-1}$ associated with the eigenvector $\Delta \z_k$. Hence, we interpret $\beta_\z$ in terms of its influence on the spectrum of $\E_\z^{-1}$, which determines the response of $\d$ to high-frequency perturbations $\Delta \z_k$. Large values of $\beta_\z$ will result in $\E_\z^{-1}$ having rapidly decaying eigenvalues so that $\d(\ztilde+\Delta \z_k,\t_i) \approx \d(\ztilde,\t_i)$ for all but a small number of leading low-frequency perturbations. Conversely, taking small values of $\beta_\z$ will result in $\E_\z^{-1}$ having a ``flat" spectrum (many eigenvalues close to $1$) and consequently $\d(\ztilde+\Delta \z_k,\t_i) -\d(\ztilde,\t_i)$ will be large for high-frequency perturbations $\Delta \z_k$. Visualizing the discrepancy response to low- and high-frequently perturbations provides the insight to verify or adjust $\beta_\z$. Assuming that $\beta_\z$ is specified appropriately, observe from~\eqref{eqn:delta_diff_pert} that $\sqrt{\alpha_\z}$ scales the magnitude of $ \d(\ztilde+\Delta \z_k,\t_i)-\d(\ztilde,\t_i) $. Hence, $\alpha_\z$ is verified or modified by linear scaling of $\sqrt{\alpha_\z}$ based upon the magnitude of the discrepancy differences. 
 
Given these observations about $\alpha_\z$ and $\beta_\z$, we use the maximum absolute value of $\d(\ztilde+\Delta \z_k,\t_i) -\d(\ztilde,\t_i)$ and the correlation length of $\Delta \z_k$ as metrics to project $\{ \Delta \z_k,\d(\ztilde+\Delta \z_k,\t_i) -\d(\ztilde,\t_i)\}_{i=1,k=1}^{Q,P}$ onto a 2D scatter plot. We select extreme points in the 2D scatter plot and visualize $\Delta \z_k$ and $\d(\ztilde + \Delta \z_k,\t_i)-\d(\ztilde,\t_i)$ to guide the verification or modification of $\alpha_\z$ and $\beta_\z$, as demonstrated in Figure~\ref{fig:interactive_viz_snapshot}.

\subsubsection*{Summary}

\begin{figure}[h]
\centering
  \includegraphics[width=0.64\textwidth]{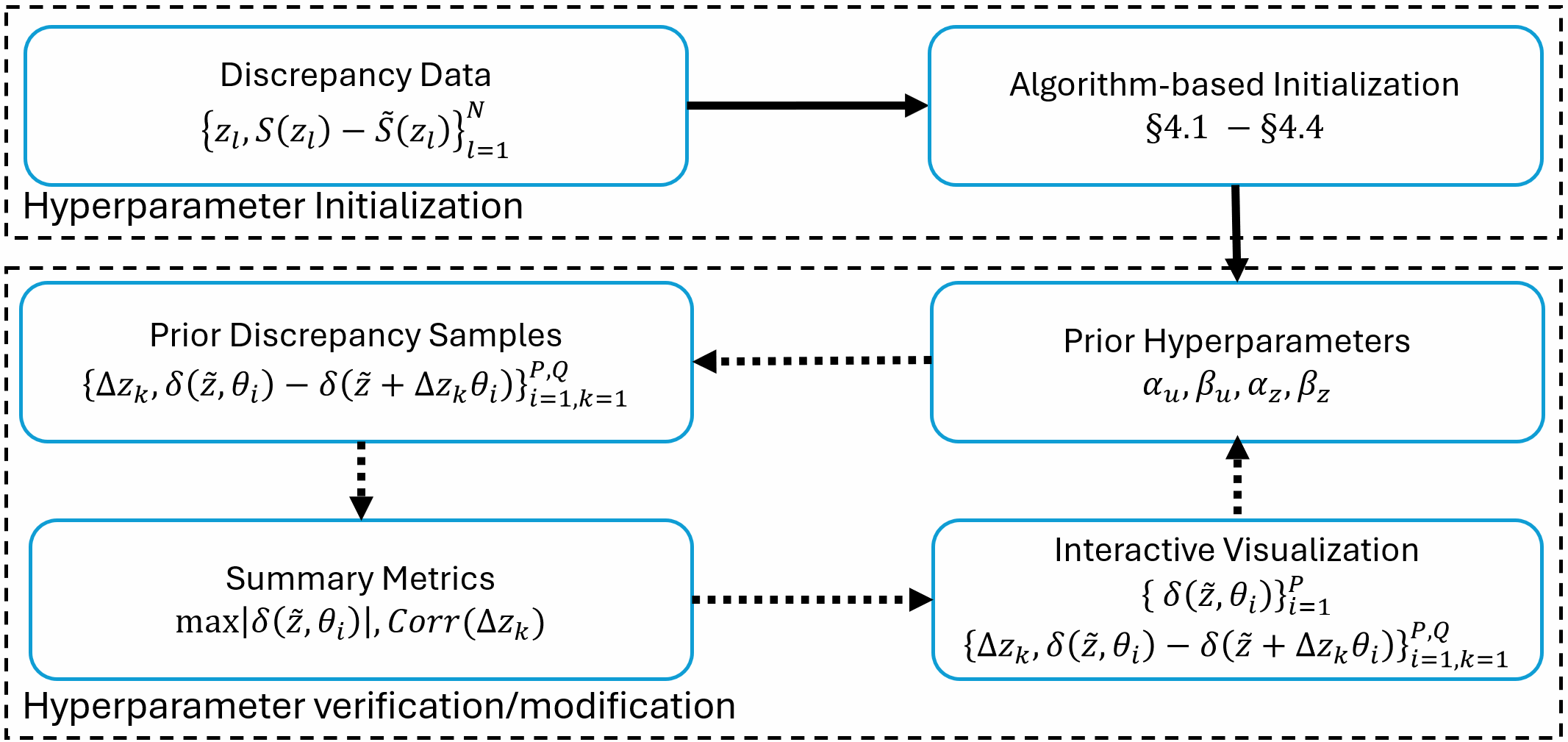}
    \caption{Depiction of the hyper-parameter specification process. The upper box indicates the algorithm-based initialization, which is performed once. The lower box indicates the visualization framework, which may be repeated multiple times if necessary to achieve appropriate hyper-parameters.}
  \label{fig:hyperparam_spec_process}
\end{figure}

We summarize our hyper-parameter specification process in Figure~\ref{fig:hyperparam_spec_process} and emphasize the important role of both algorithm-based initialization of the hyper-parameters to get ``close enough" and visualization techniques that synthesize complex prior sample data in 2D scatter plots to enable interactive rendering of visualizations for select samples and perturbations. Our use of algorithm-based initialization in the top part of Figure~\ref{fig:hyperparam_spec_process} ensures that the ``hyper-parameter verification/modification" loop in the lower part of Figure~\ref{fig:hyperparam_spec_process} does not require many iterations. Without a good initialization, the hyper-parameter verification/modification loop is laborious and time consuming.

\section{Numerical results} \label{sec:numerics}

This section illustrates our hyper-parameter specification framework and how the prior discrepancy depends on the hyper-parameters. We begin with a stationary problem in one spatial dimension for ease of explanation and visualization. Subsequently, we generalize to a transient problem in one spatial dimension to demonstrate the effect of the temporal variance hyper-parameter. Lastly, we consider a stationary problem in two spatial dimensions, where the hyper-parameter initialization algorithms produce a poor prior. We demonstrate the role of visualization to identify the issue and correct the initialization by adjusting the hyper-parameters.
  
\subsection{Stationary 1D problem} \label{ssec:stationary_1D_example}
  
In this subsection, we illustrate properties of our proposed hyper-parameter initialization algorithm and how the prior discrepancy distribution is affected by the hyper-parameters. For these illustrations, we focus on a stationary problem in $s=1$ spatial dimension and consider a high-fidelity model given by the advection-diffusion equation
\begin{align*}
& -\frac{d^2 u}{dx^2} + \frac{du}{dx} = z \qquad & \text{on } (0,1) \\
& -\frac{du}{dx}(0) = 2 u(0) \qquad \text{and} \qquad \frac{du}{dx}(1) = 2 u(1),
\end{align*}
which constrains an optimization problem to determine the optimal source function $z:[0,1] \to \R$ to achieve a desired target state $T(x)=50 - 60 (x - 0.5)^2$. We let $S(z)$ denote the solution operator for the high-fidelity advection-diffusion model and consider a low-fidelity model given by the diffusion equation
\begin{align*}
& -\frac{d^2 u}{dx^2} = z \qquad & \text{on } (0,1) \\
& -\frac{du}{dx}(0) = 2 u(0) \qquad \text{and} \qquad \frac{du}{dx}(1) = 2 u(1) 
\end{align*}
whose solution operator is denoted as $\tilde{S}(z)$. We seek to determine hyper-parameters $\alpha_\u, \beta_\u, \alpha_\z$, and $\beta_\z$, appearing in the covariance $\W_\t^{-1}$, so that the prior discrepancy has characteristics commensurate to $S(z)-\tilde{S}(z)$.

\begin{figure}[h]
\centering
  \includegraphics[width=0.234\textwidth]{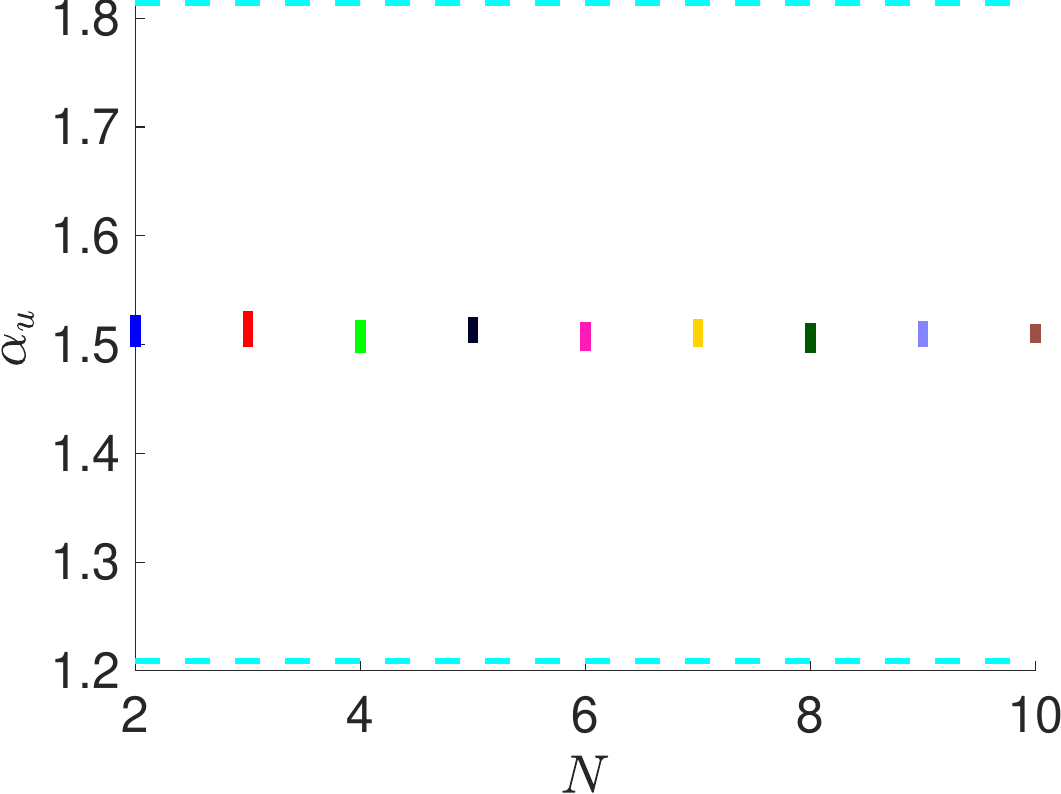}
    \includegraphics[width=0.234\textwidth]{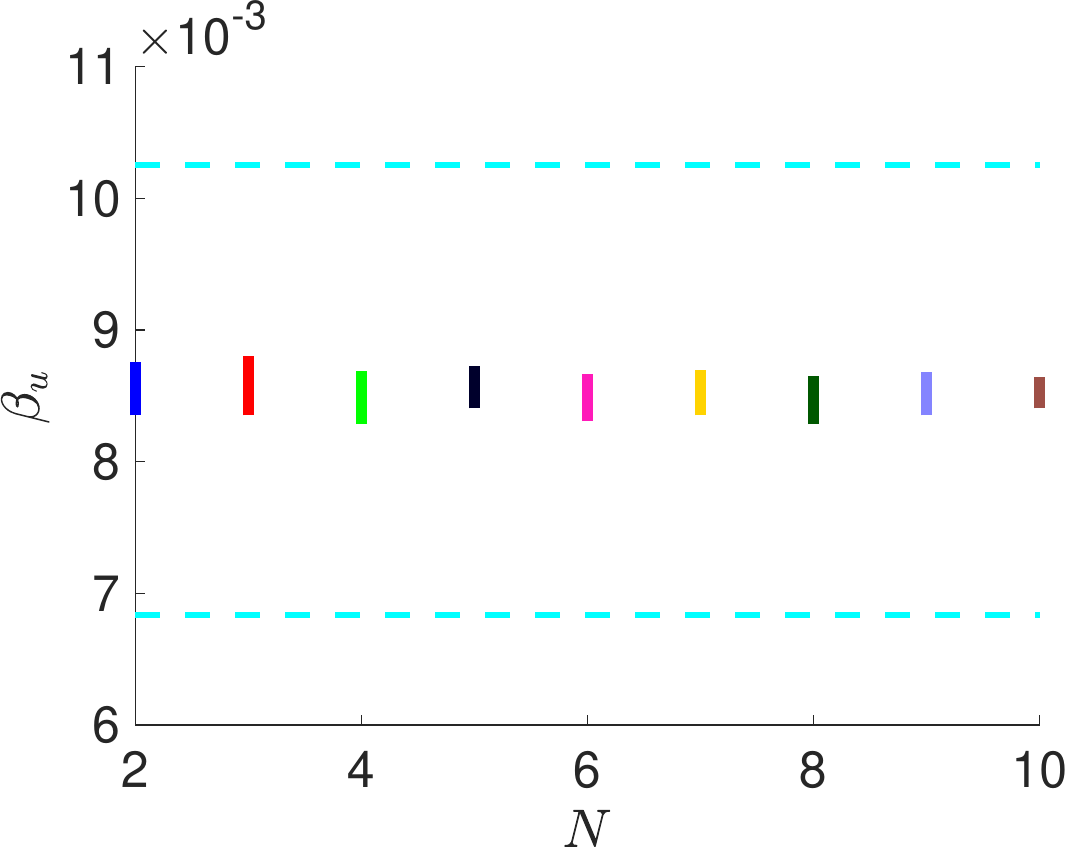}
  \includegraphics[width=0.234\textwidth]{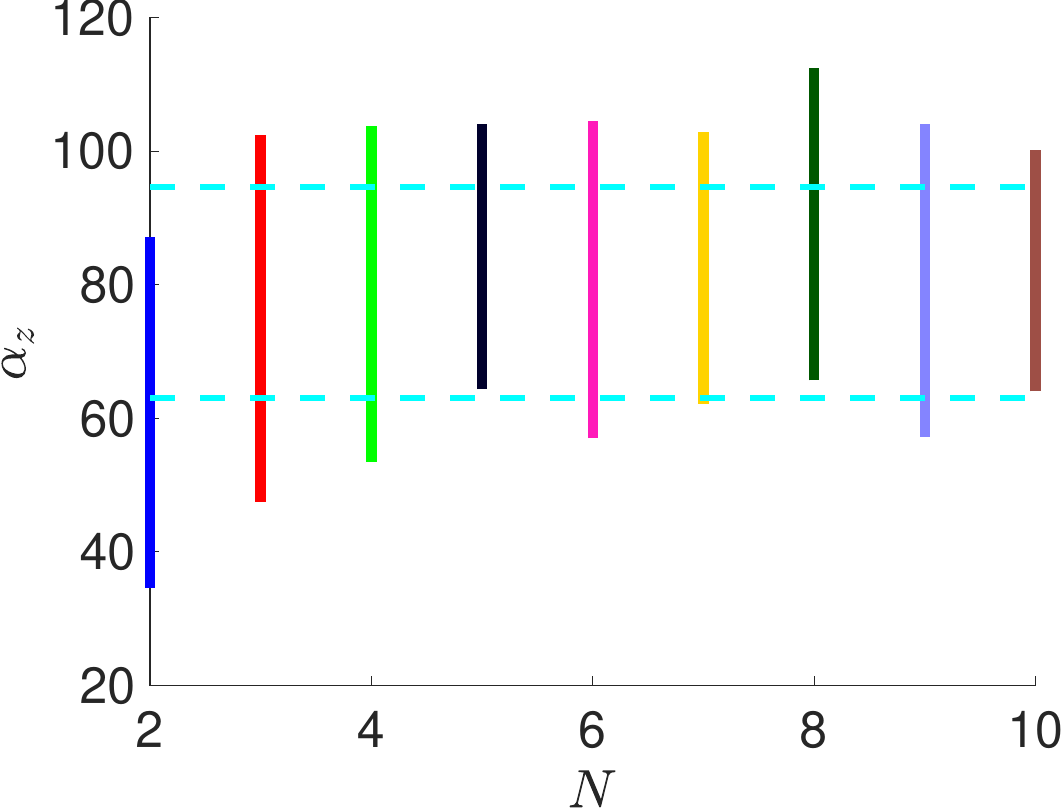}
    \includegraphics[width=0.234\textwidth]{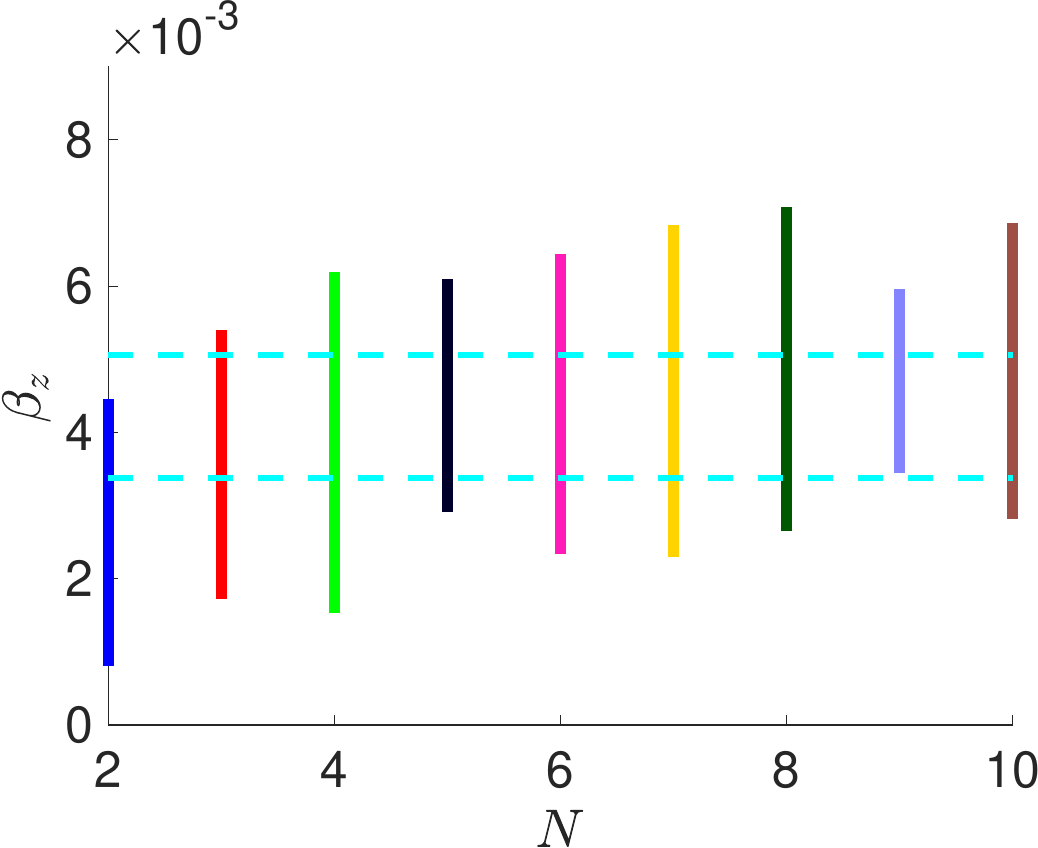}
    \caption{The panels from left to right correspond to: $\alpha_\u, \beta_\u, \alpha_\z$, and $\beta_\z$. The number of datapoints $N$ is varied on the horizontal axis and the range of hyper-parameters generated over all 50 datasets (generated uniquely for each $N$) is depicted by the vertical lines. The horizontal lines denote $\pm 20\%$ of the mean hyper-parameter value from all datasets.}
  \label{fig:hyperparam_variability}
\end{figure}

Our first illustration considers robustness of the hyper-parameter initialization algorithms with respect to variability in the data $\{ \z_\ell ,\vec{d}_\ell \}_{\ell=1}^N$. To this end, we vary the number of datapoints $N$, and for each fixed $N$, generate $50$ unique datasets. These $50$ datasets are generated by sampling the source terms $\{ \z_\ell \}_{\ell=2}^N$ from a Gaussian random field ($\z_1=\ztilde$ is fixed as the result of solving the low-fidelity optimization problem). The corresponding discrepancy data $\{ \vec{d} \}_{\ell=2}^N$ is recomputed for each dataset by evaluating the high- and low-fidelity models. Figure~\ref{fig:hyperparam_variability} displays the resulting hyparameters. Each panel corresponds to a hyper-parameter ($\alpha_\u, \beta_\u, \alpha_\z$, and $\beta_\z$) and displays the hyper-parameter value as a function of $N$ on the horizontal axis, with variability over the $50$ datasets per $N$ displayed by the range of the vertical lines. For perspective on scales, we average over all datasets and plot horizontal lines to denote $\pm 20\%$ variation from the mean hyper-parameter estimate.
  
  We draw two conclusions from Figure~\ref{fig:hyperparam_variability}:
  \begin{enumerate}
  \item[$\bullet$] The state hyper-parameters $\alpha_\u$ and $\beta_\u$ can be estimated more robustly than the optimization variable hyper-parameters $\alpha_\z$ and $\beta_\z$.
  \item[$\bullet$] The primary source of optimization variable hyper-parameter variability is the empirical correlation length calculations (which depend on the data) that determine $\beta_\z$. The variability of $\alpha_\z$ is primarily due to its dependence on $\beta_\z$, although some variability is also due to estimation of $\gamma^2$.
  \end{enumerate}
  
\begin{figure}
  \centering
  \begin{tabular}{M{.04\textwidth}M{.25\textwidth}M{.25\textwidth}M{.25\textwidth}}

 \hspace{1 mm}  & 
    \includegraphics[width=0.23\textwidth]{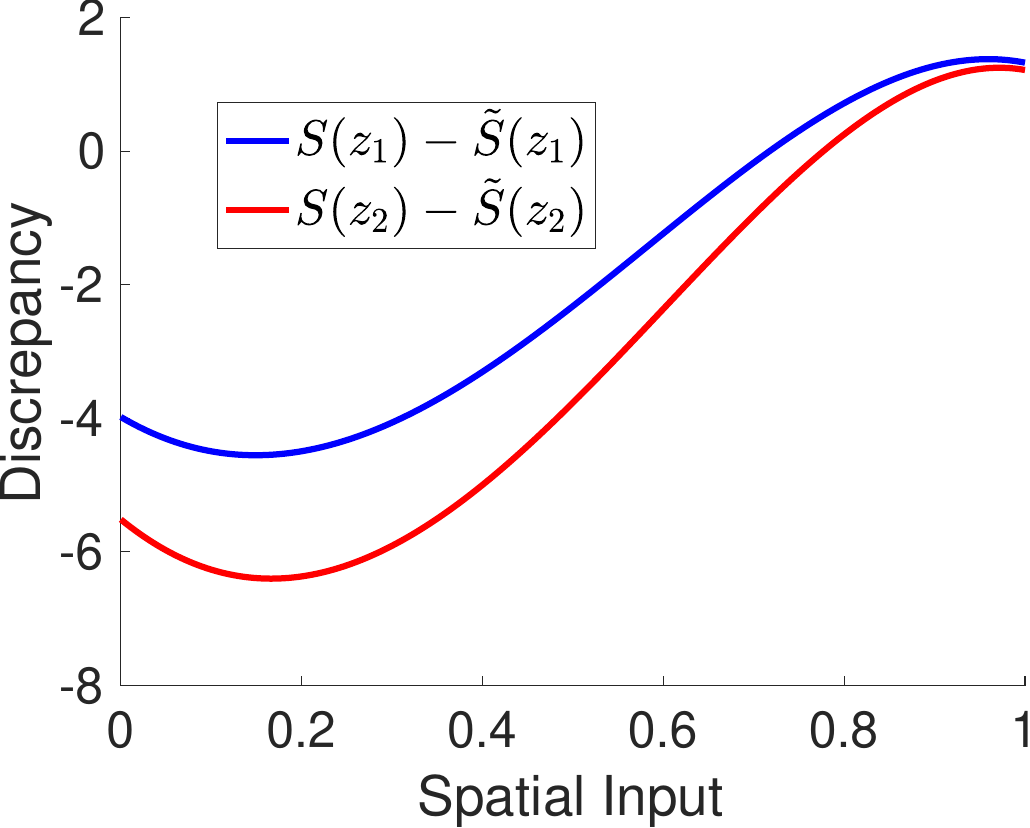} &
  \includegraphics[width=0.23\textwidth]{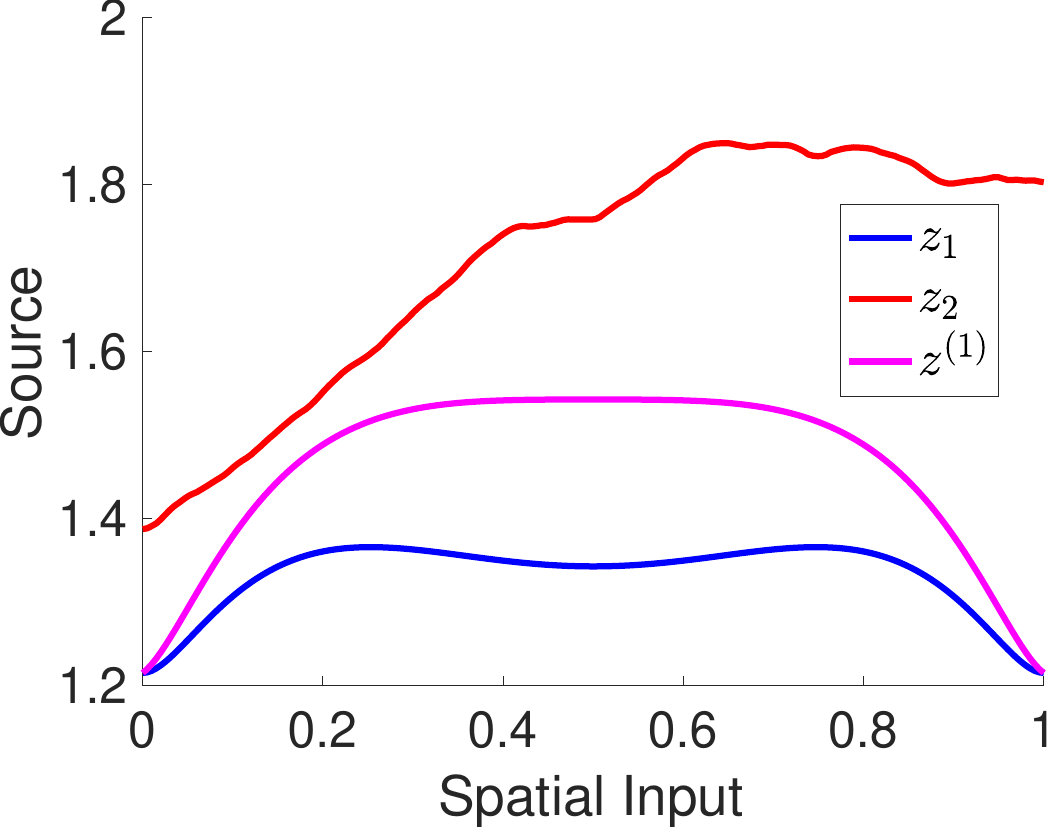} &
    \includegraphics[width=0.23\textwidth]{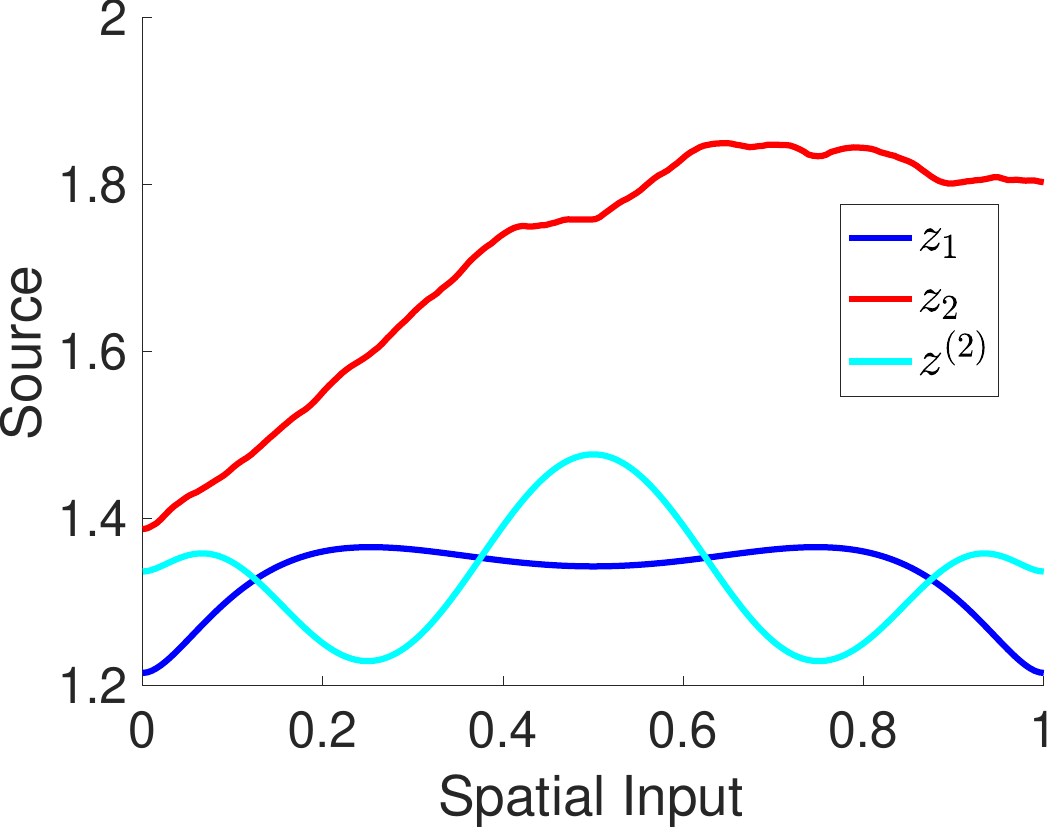} \\
     \hspace{1 mm}  & $\hspace{2mm} \d(\z_1,\t_i)$ & $\hspace{5mm} \d(\z^{(1)}-\z_1,\t_i)$ & $\hspace{5mm} \d(\z^{(2)}-\z_1,\t_i)$ \\
         \rotatebox{90}{Nominal} &
            \includegraphics[width=0.23\textwidth]{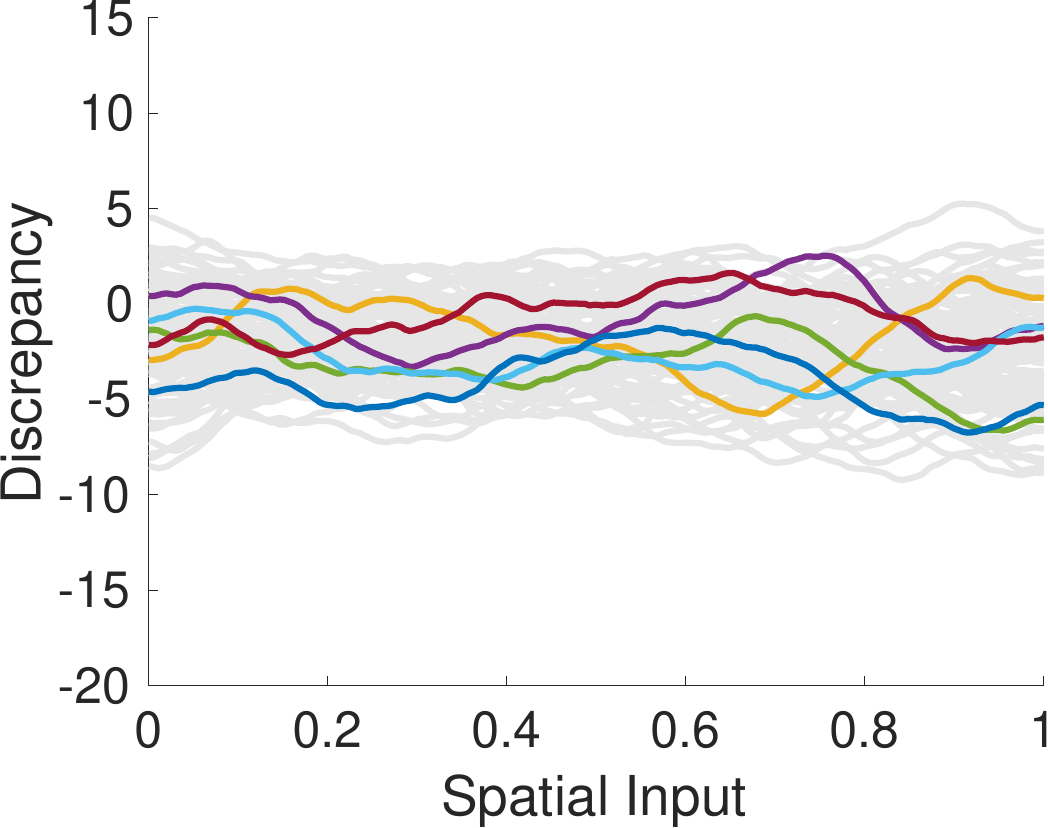} &
   \includegraphics[width=0.23\textwidth]{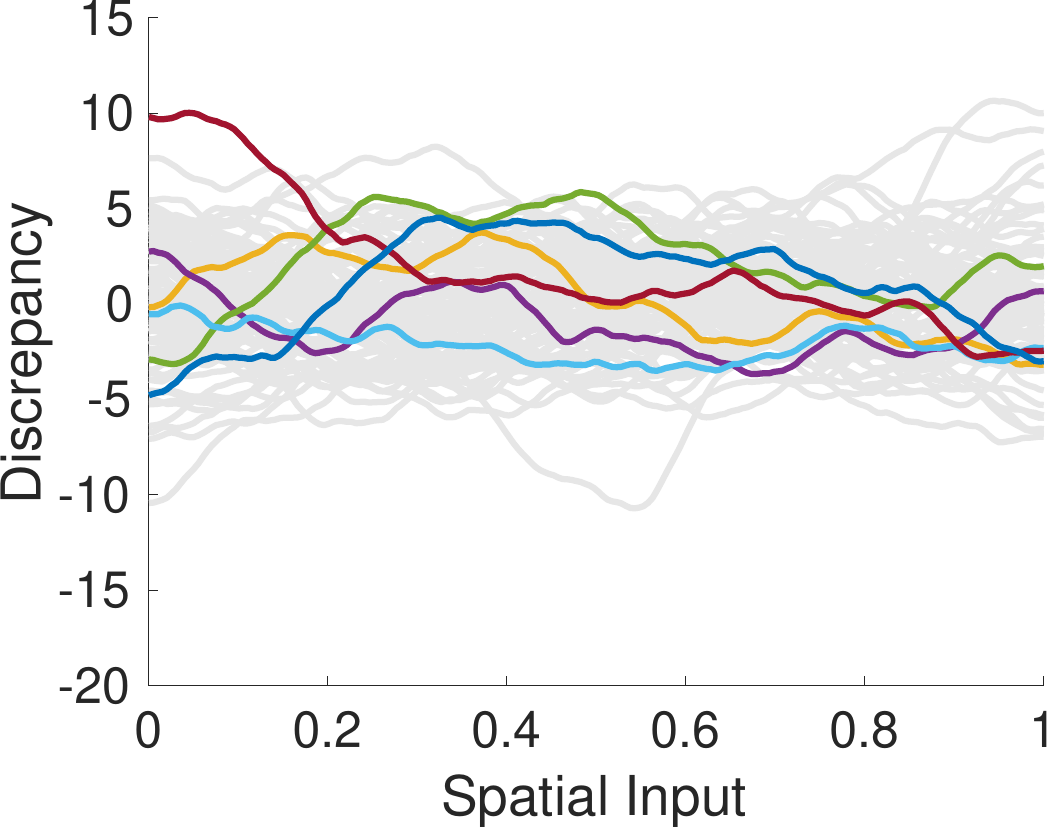} &
   \includegraphics[width=0.23\textwidth]{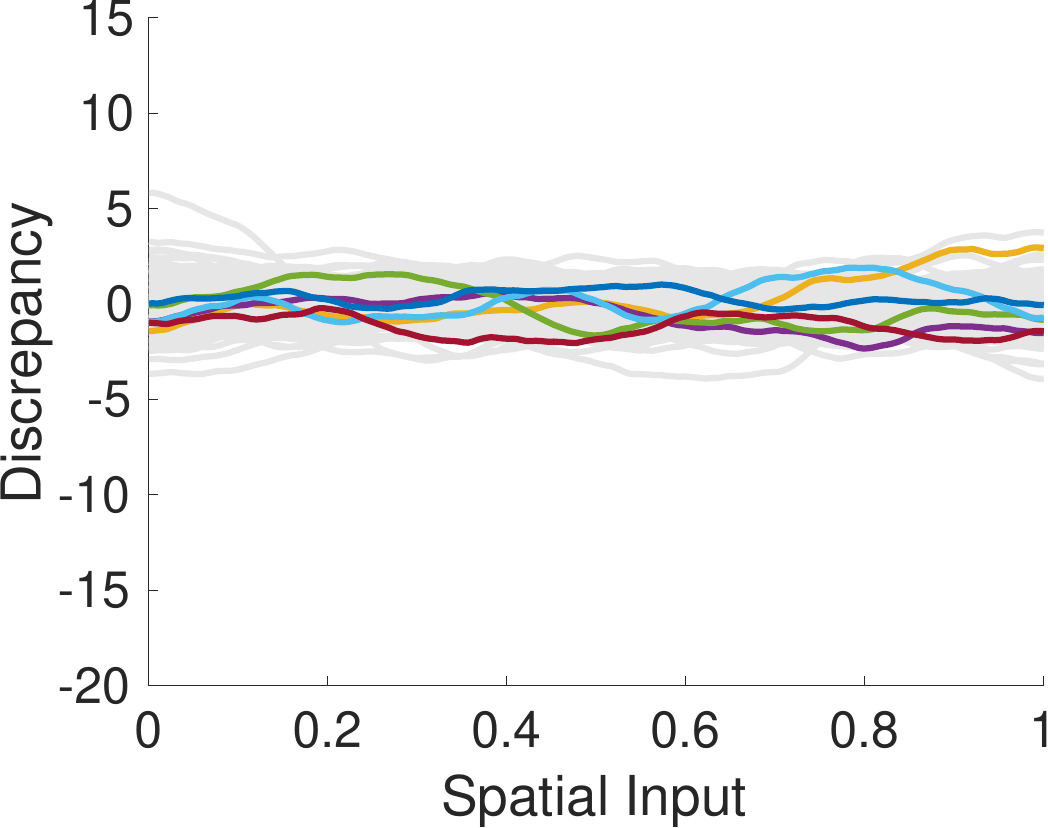} \\
            $\frac{\alpha_\u}{10}$ &
            \includegraphics[width=0.23\textwidth]{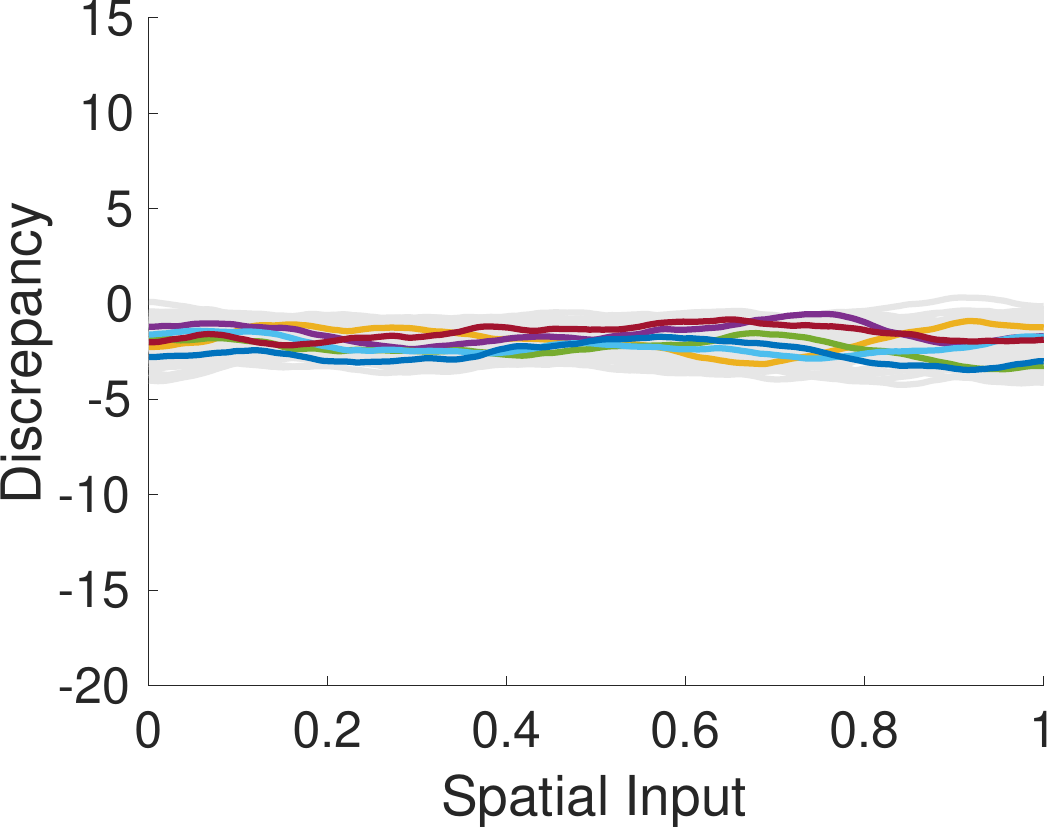} &
   \includegraphics[width=0.23\textwidth]{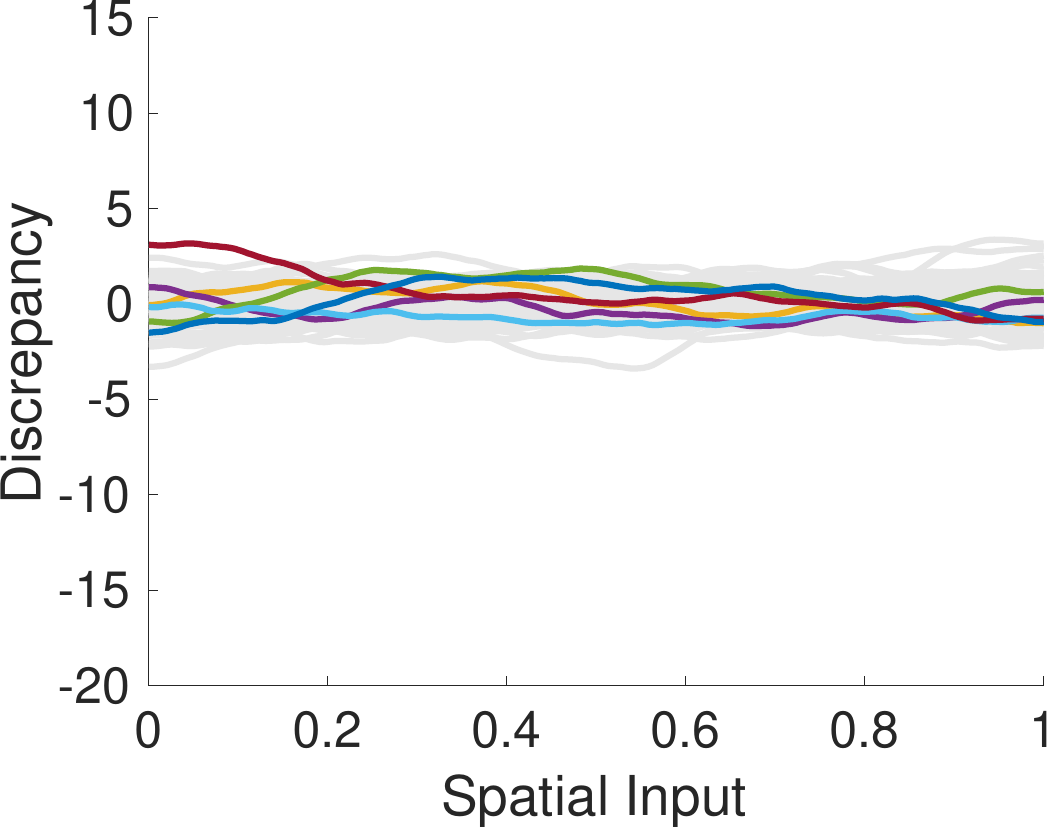} &
   \includegraphics[width=0.23\textwidth]{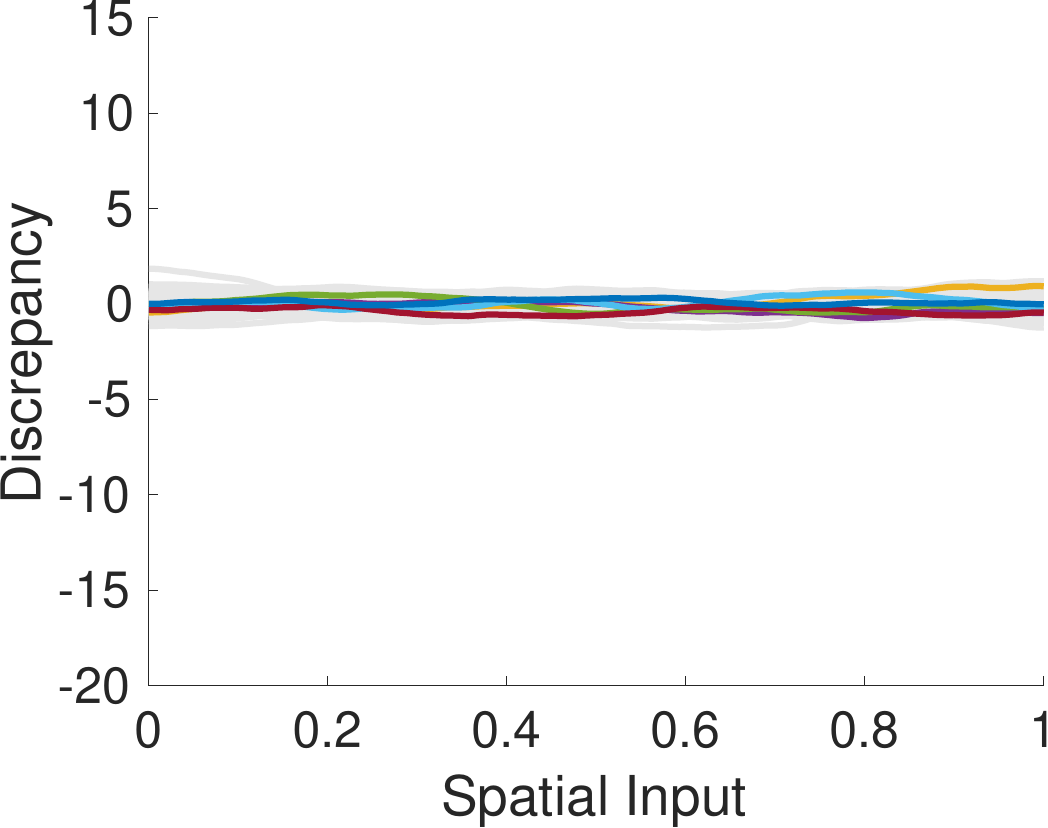} \\
               $\frac{\beta_\u}{10}$ &
            \includegraphics[width=0.23\textwidth]{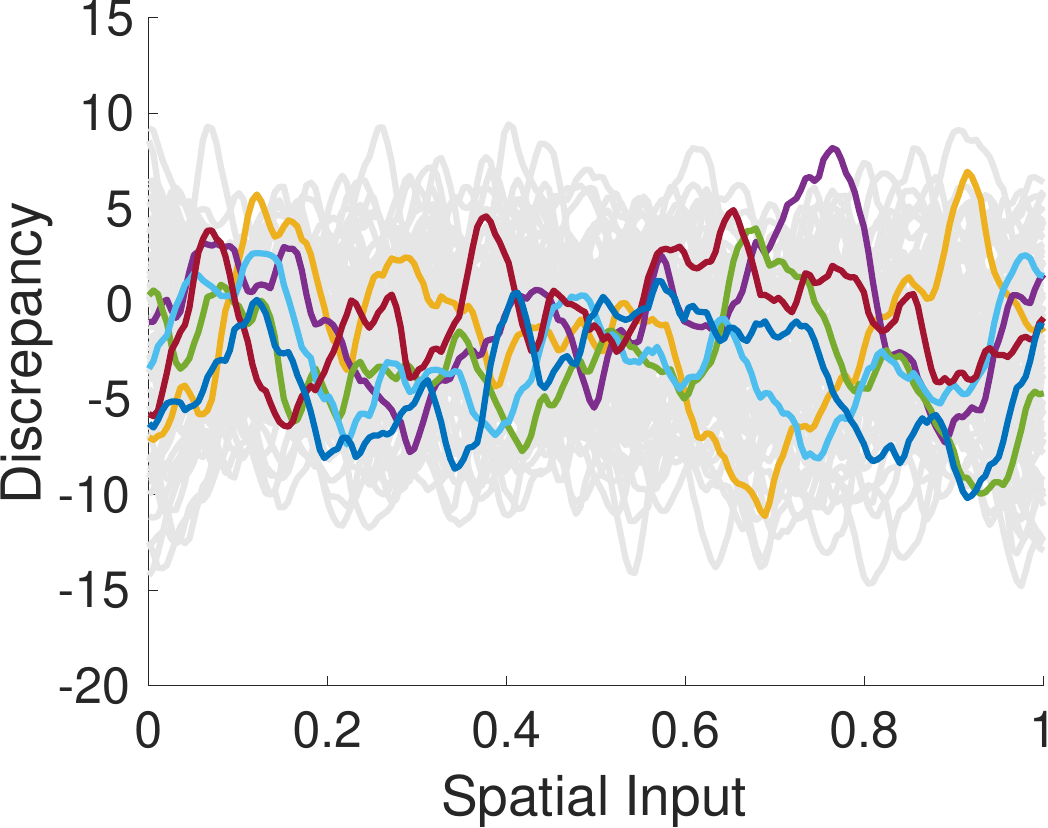} &
   \includegraphics[width=0.23\textwidth]{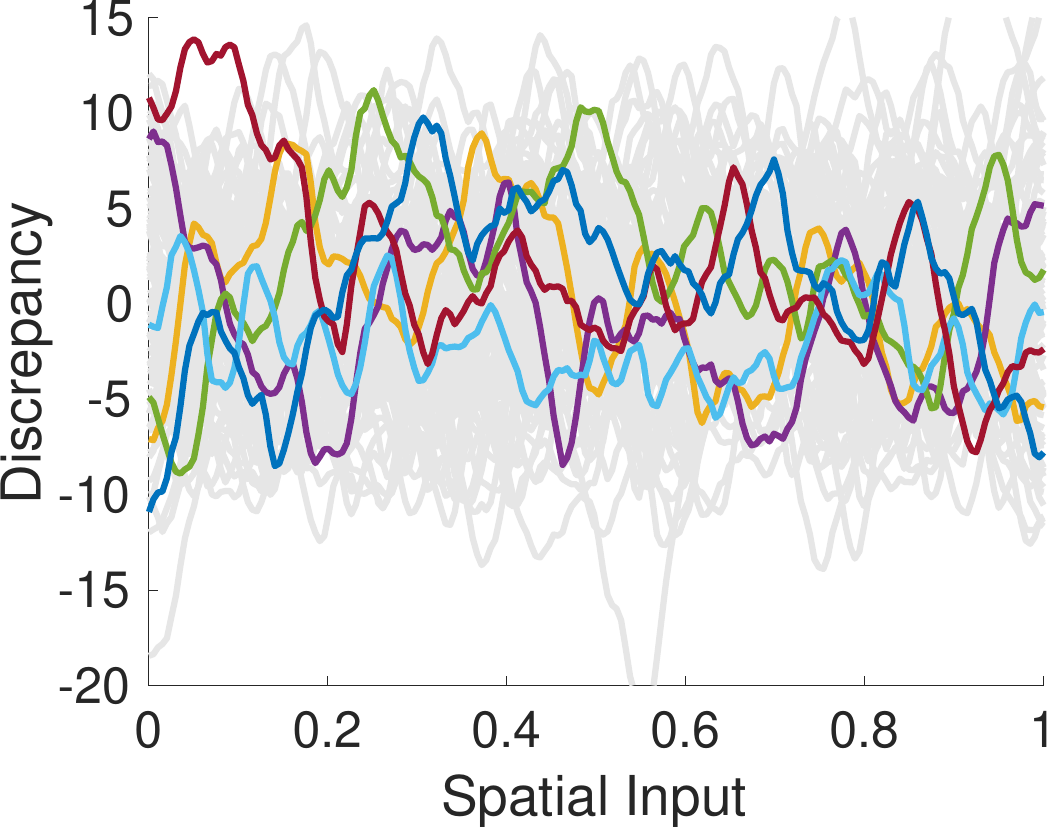} &
   \includegraphics[width=0.23\textwidth]{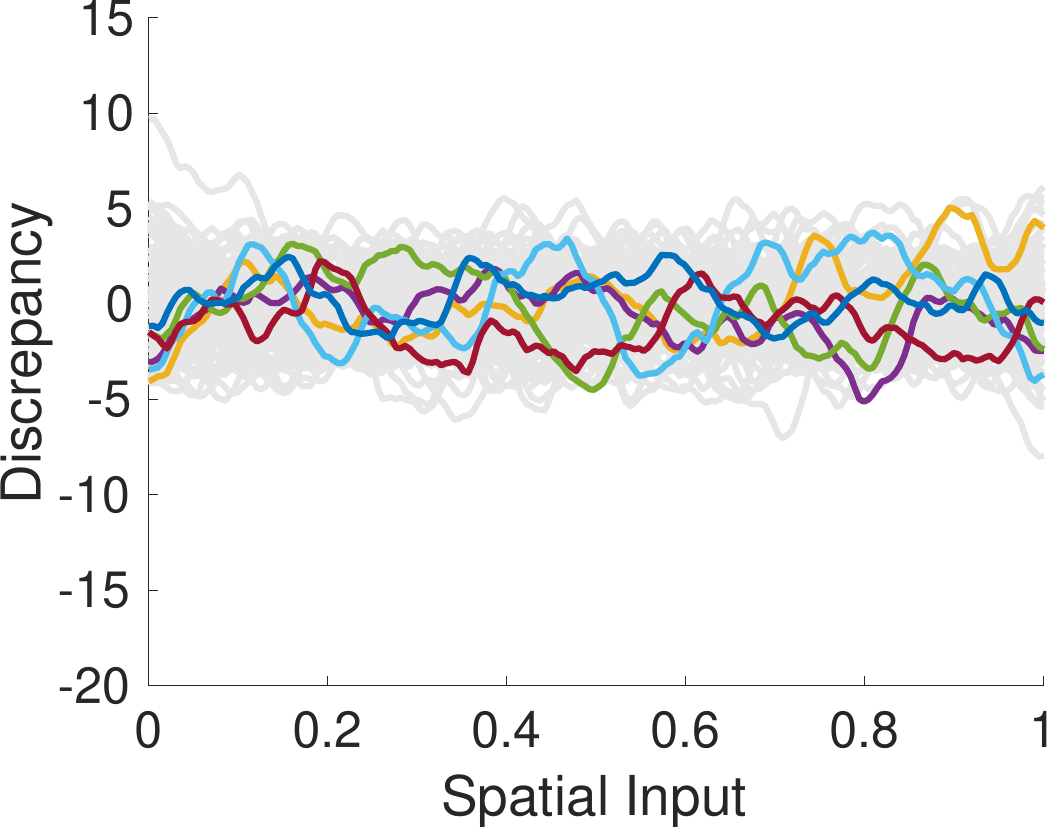} \\
                  $\frac{\alpha_\z}{10}$ &
            \includegraphics[width=0.23\textwidth]{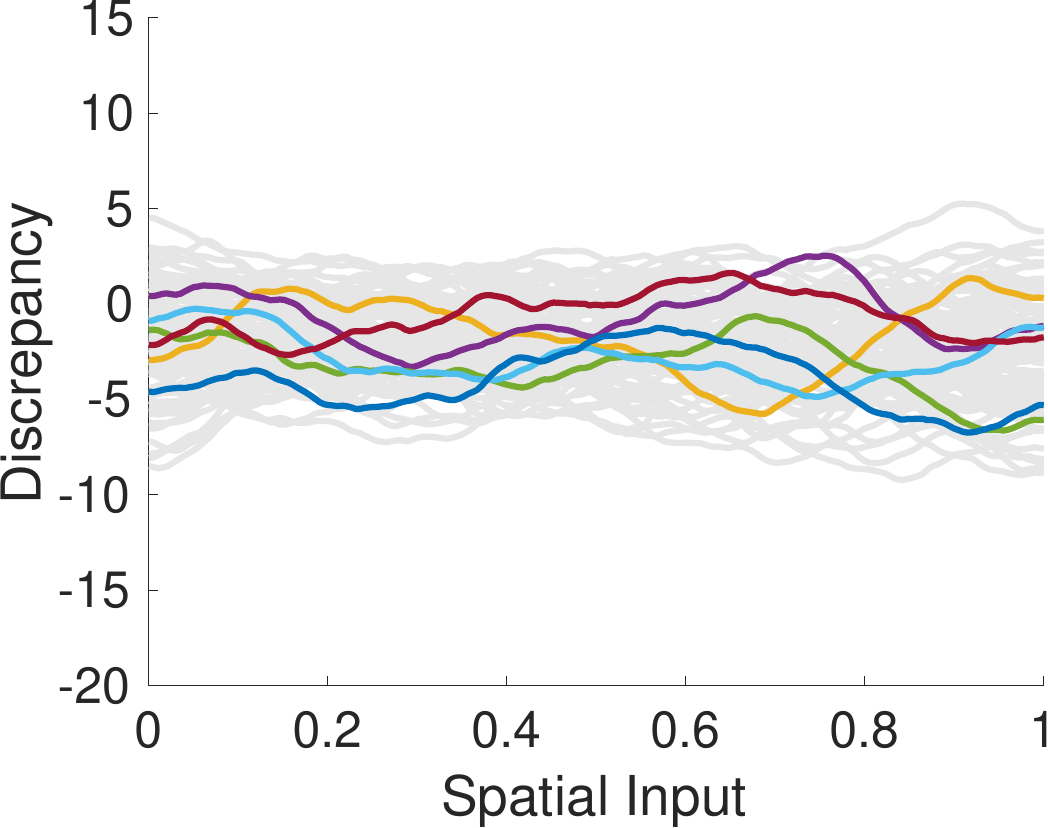} &
   \includegraphics[width=0.23\textwidth]{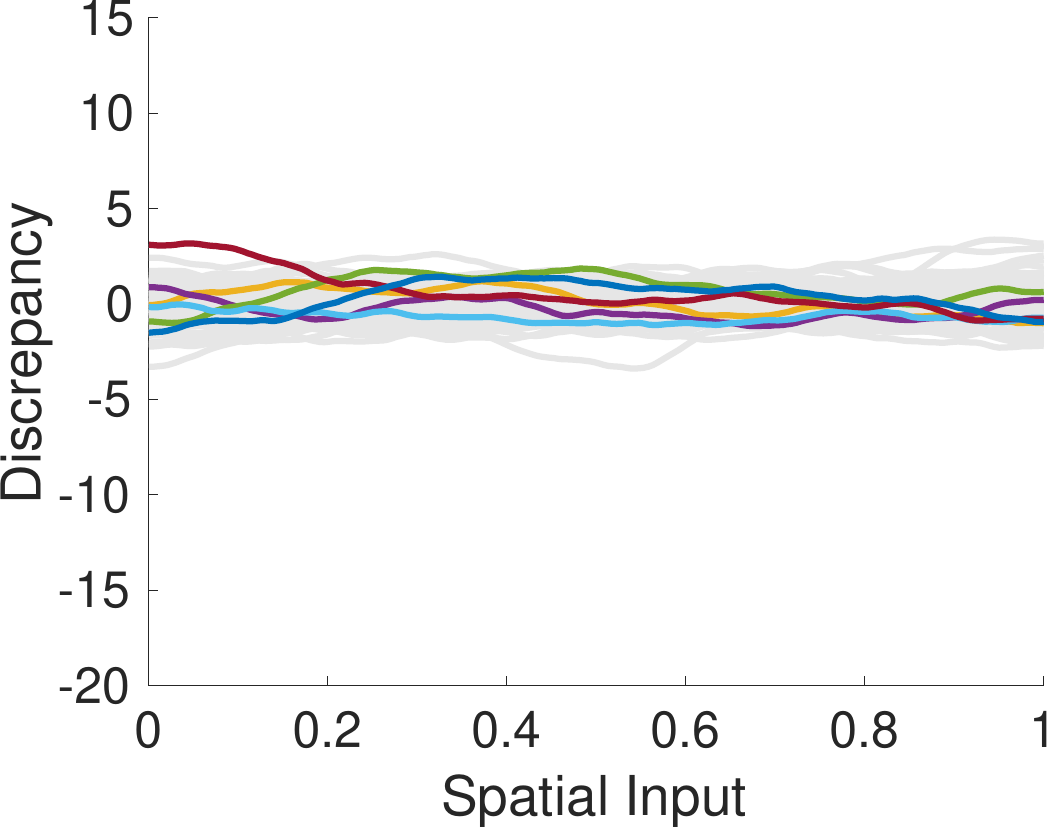} &
   \includegraphics[width=0.23\textwidth]{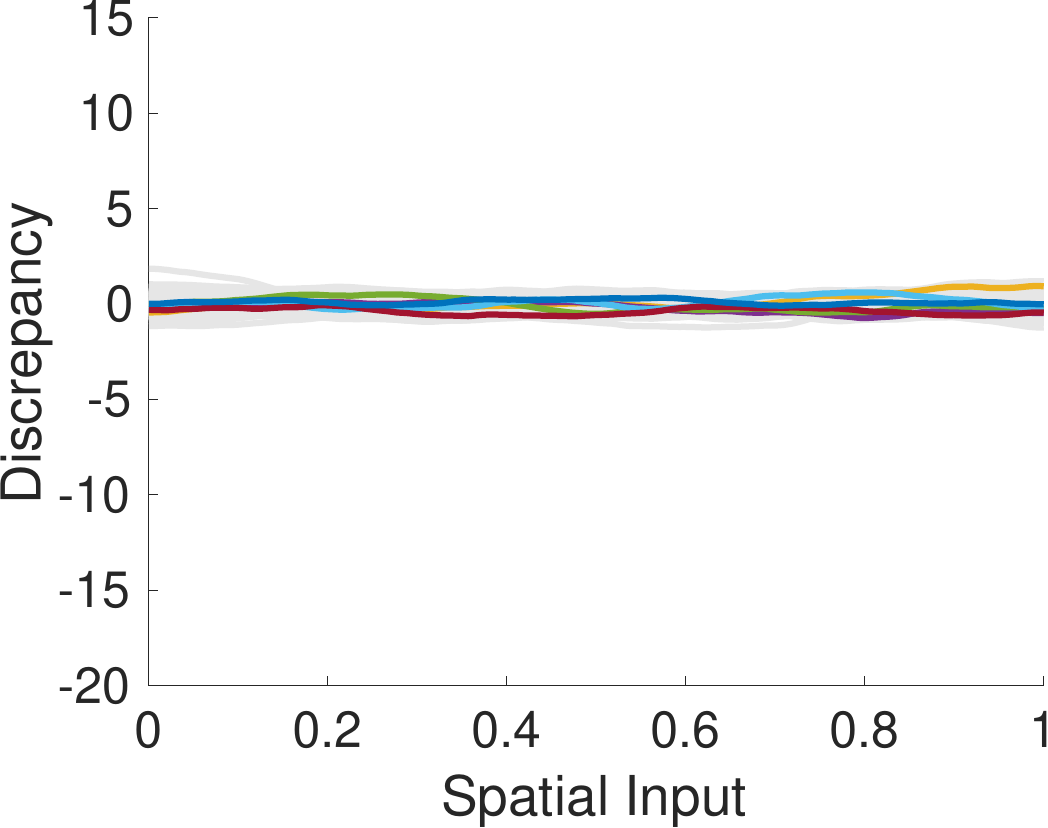} \\
                     $\frac{\beta_\z}{10}$ &
            \includegraphics[width=0.23\textwidth]{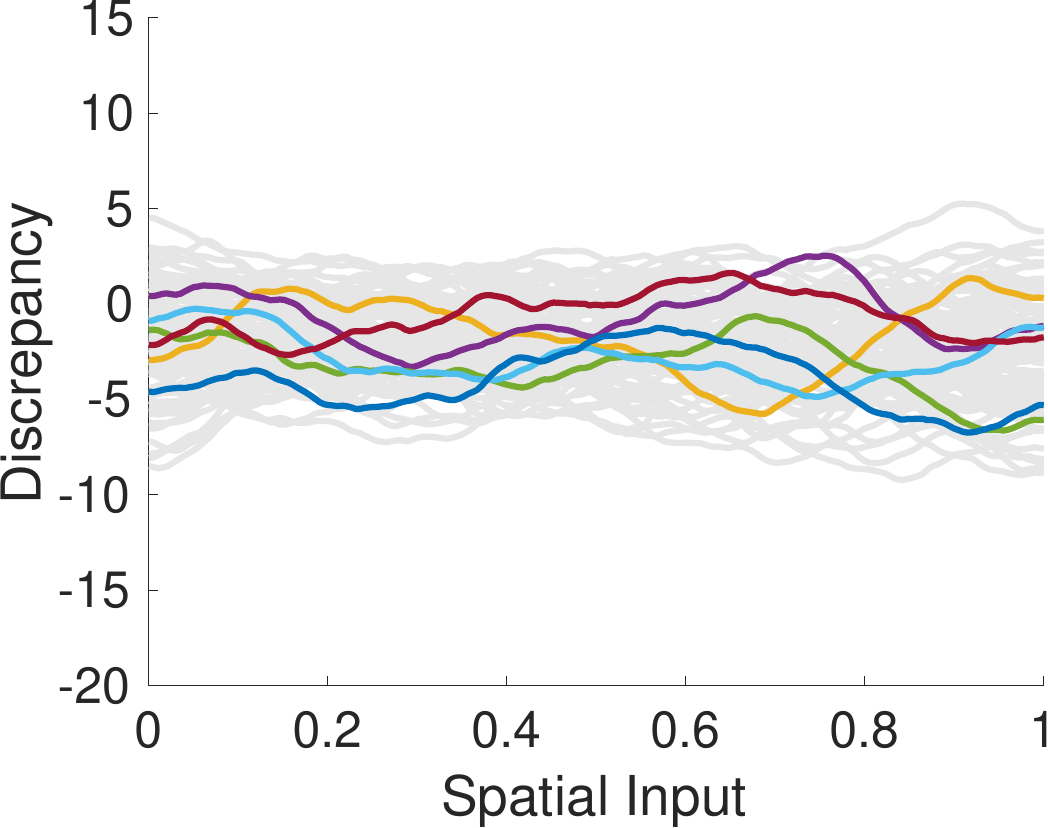} &
   \includegraphics[width=0.23\textwidth]{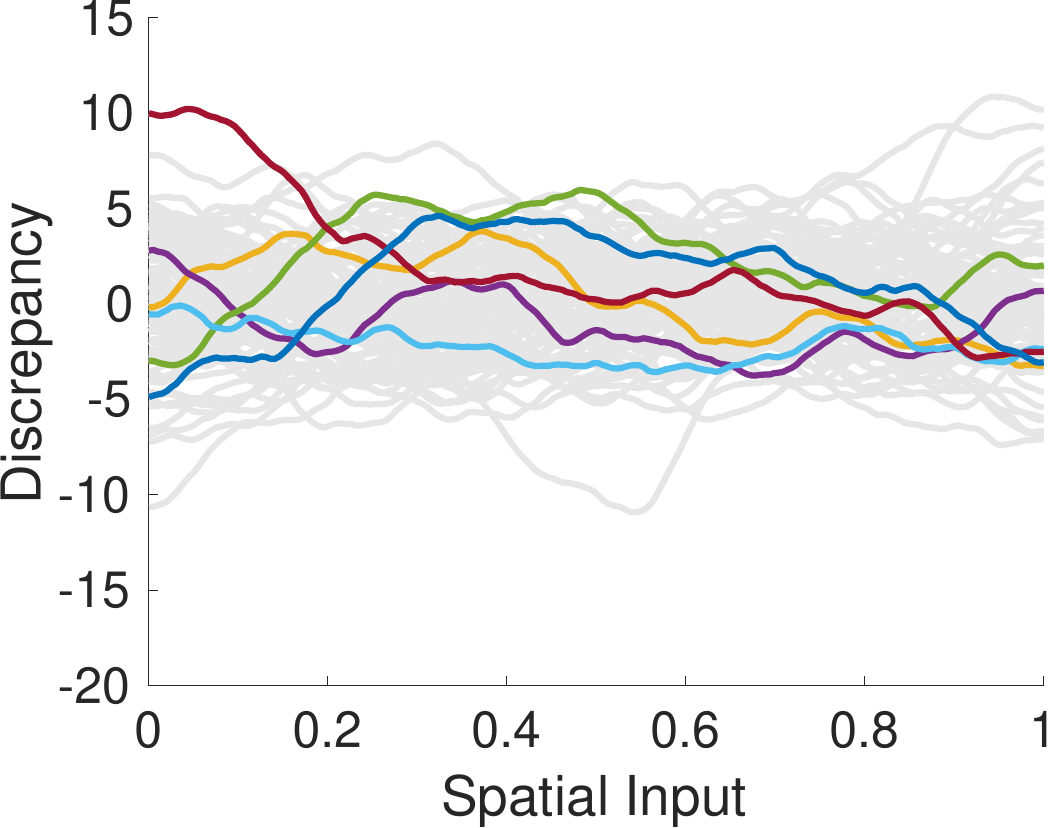} &
   \includegraphics[width=0.23\textwidth]{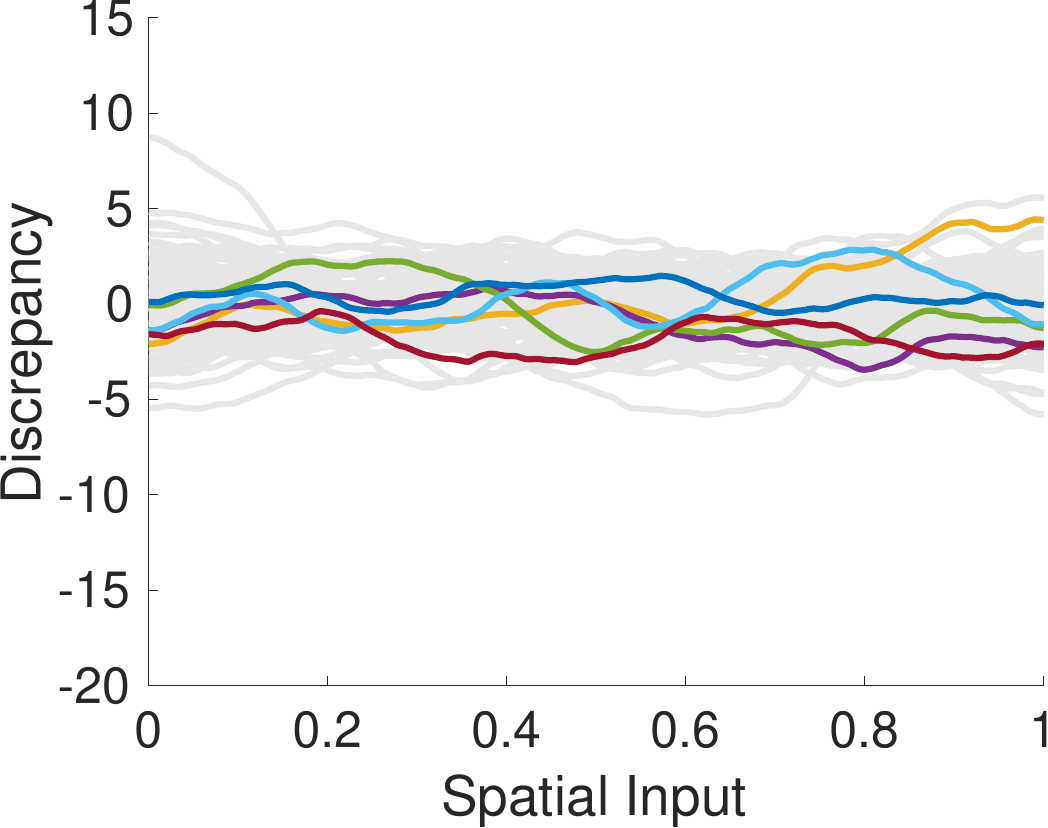} \\
  \end{tabular}
    \caption{
   \textbf{First row:} the discrepancy data $\{ \vec{d}_\ell \}_{\ell=1}^2$ (left) and the corresponding optimization variable data $\{ \z_\ell \}_{\ell=1}^2$ (center and right). The optimization variable inputs $\z^{(1)}$ and $\z^{(2)}$ are shown in the center and right columns, respectively. \textbf{Second row:} samples of $\d(\cdot,\t_i)$ evaluated at $\z_1$, $\z^{(1)}-\z_1$, and $\z^{(2)}-\z_1$, from left to right. Each curve is plotted in a gray scale with select samples plotted in color to highlight their characteristics. \textbf{Third through sixth rows:} analogous plots to the second row with each row corresponding to a different hyper-parameters divided by 10. }
  \label{fig:prior_samples}
\end{figure}

  To provide intuition on the effect of the prior hyper-parameters, we plot samples for various perturbations of the hyper-parameters. Specifically, we use a single dataset with $N=2$ to initialize the hyper-parameters. The discrepancy data $\{ \vec{d}_\ell \}_{\ell=1}^2$ is depicted in the left panel of the top row in Figure~\ref{fig:prior_samples}. The corresponding optimization variable data $\{ \z_\ell \}_{\ell=1}^2$ is shown in the center and right panels of the top row in Figure~\ref{fig:prior_samples}. To facilitate visualization of the prior discrepancy, we compute prior samples $\{\t_i\}_{i=1}^Q$ and evaluate $\d(\z,\t_i)$ for $\z \in \{\ztilde,\z^{(1)},\z^{(2)}\}$, where $\z^{(1)}$ and $\z^{(2)}$ are chosen for illustrative purposes. The top row (center and right panels) of Figure~\ref{fig:prior_samples} display $\z^{(1)}$ and $\z^{(2)}$.
  
In the second through sixth rows of Figure~\ref{fig:prior_samples}, we show samples of $\d(\z_1,\t_i)$, $\d(\z^{(1)}-\z_1,\t_i)=\d(\z^{(1)},\t_i)-\d(\z_1,\t_i)$, and $\d(\z^{(2)}-\z_1,\t_i)=\d(\z^{(2)},\t_i)-\d(\z_1,\t_i)$, in the left, center, and right columns, respectively. Hence, the left column demonstrates the prior discrepancy evaluated at the low-fidelity optimization solution $\ztilde$ and the center and right columns demonstrate how the prior discrepancy varies with respect to $\z$ in the directions of $\z^{(1)}$ and $\z^{(2)}$. The second row corresponds to samples using the nominal hyper-parameter values initialized using the dataset $\{\z_\ell , \vec{d}_\ell \}_{\ell=1}^2$, and each row below it corresponds to prior samples with one of the hyper-parameters modified. In each case, we reduce the hyper-parameter value by an order of magnitude to emphasize its effect. As demonstrated in Figure~\ref{fig:hyperparam_variability}, varying by an order of magnitude is very large, but it was chosen for emphasis in the visualization. We make several observations from Figure~\ref{fig:prior_samples}:
 \begin{enumerate}
 \item[$\bullet$] \textbf{From second row:} The magnitude of $\d(\z^{(1)}-\z_1,\t_i)=\d(\z^{(1)},\t_i)-\d(\z_1,\t_i)$ is greater than the magnitude of $\d(\z^{(2)}-\z_1,\t_i)=\d(\z^{(2)},\t_i)-\d(\z_1,\t_i)$. This is due to $\z^{(2)}-\z_1$ being more oscillatory that $\z^{(1)}-\z_1$, and hence $\z^{(2)}-\z_1$ is damped by the elliptic operator $\E_\z^{-1}$, yielding a smaller magnitude change in the discrepancy. 
 \item[$\bullet$] \textbf{From third row:} The uncertainty in the discrepancy is proportional to $\alpha_\u$, hence reducing $\alpha_\u$ by an order of magnitude reduces the range of magnitudes in the samples.
 \item[$\bullet$] \textbf{From fourth row:} Decreasing $\beta_\u$ results in discrepancy samples that are less smooth, and consequently some samples have extrema with larger magnitudes than in the nominal hyper-parameter case.
 \item[$\bullet$] \textbf{From fifth row:} Decreasing $\alpha_\z$ results in $\d(\z^{(1)},\t_i)-\d(\z_1,\t_i)$ and $\d(\z^{(2)},\t_i)-\d(\z_1,\t_i)$ having a small magnitude, i.e., the discrepancy has less sensitivity to changes in $\z$.
 \item[$\bullet$] \textbf{From sixth row:} Decreasing $\beta_\z$ results in $\d(\z^{(2)},\t_i)-\d(\z_1,\t_i)$ having a slightly larger magnitude since its oscillatory characteristic is not damped as significantly. This is most easily seen by examining the largest magnitude samples in the leftmost (spatial input $=0$) region of the rightmost panel. For instance, the largest magnitude curve attains a maximum value of approximately $9$, whereas the maximum value is approximately $6$ in the nominal hyper-parameter case.
    \end{enumerate}
  
\begin{figure}[h]
\centering
  \begin{tabular}{M{.32\textwidth}M{.32\textwidth}M{.32\textwidth}}
  Nominal & $\frac{\alpha_\u}{10}$ &  $\frac{\beta_\z}{10}$ \\
  \includegraphics[width=0.32\textwidth]{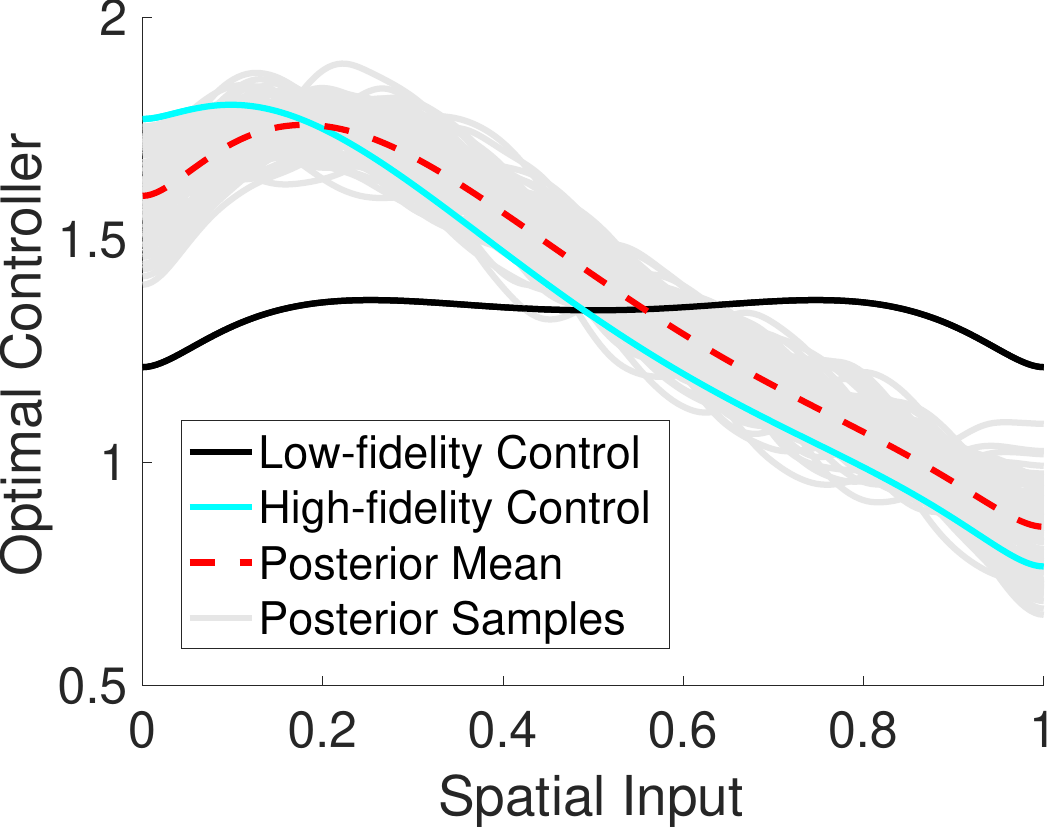} &
    \includegraphics[width=0.32\textwidth]{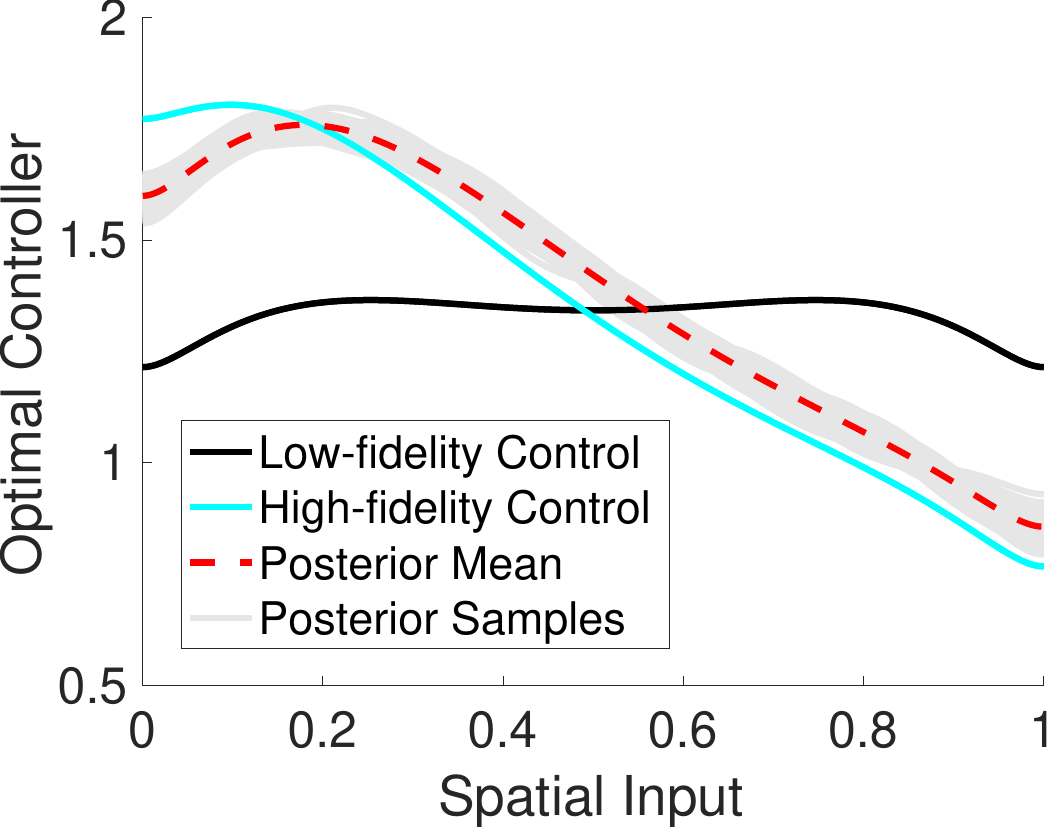}&
  \includegraphics[width=0.32\textwidth]{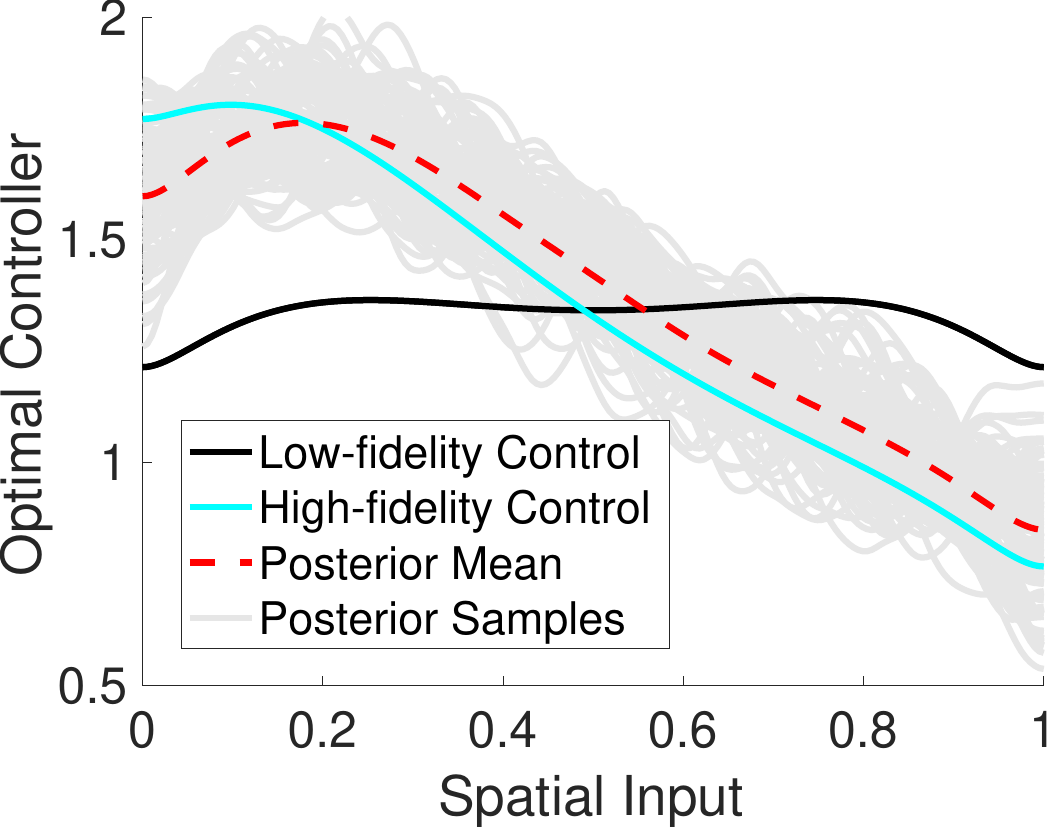}\\
  \end{tabular}
    \caption{Comparison of the posterior optimal controller distribution using the nominal hyper-parameters, modified hyper-parameters with $\alpha_\u$ divided by $10$, and modified hyper-parameters with $\beta_\z$ divided by $10$.}
  \label{fig:opt_sol_post}
\end{figure}

  To demonstrate the importance of the prior in the HDSA-MD framework, we follow the approach presented in~\cite{hart_bvw_mods} to compute the posterior model discrepancy and propagate it through the optimization problem to produce the posterior optimal solution. That is, we compute a probability distribution of the optimal controller, which reflects the uncertainty due to using only two high-fidelity model evaluations. Figure~\ref{fig:opt_sol_post} displays the posterior for three sets of prior hyper-parameters. The left panel uses the nominal hyper-parameters corresponding to the second row of Figure~\ref{fig:prior_samples} and the center and right panels use perturbed hyper-parameters corresponding to the third and sixth rows of Figure~\ref{fig:prior_samples}, respectively. In all three panels, the posterior is compared with the high- and low-fidelity optimal controllers. In practice, the high-fidelity optimal controller cannot be computed, but we include it as a ground truth for comparison here to highlight that an overly restrictive prior (when we divide $\alpha_\u$ by $10$) results in a posterior for which the high-fidelity optimal controller has approximately zero probability. The conclusions drawn from Figure~\ref{fig:opt_sol_post} motivate our study of prior model specification. In the subsections that follow, we focus on other facets of the prior model specification with an implicit awareness of the prior's influence on the posterior optimal solution.

\subsection{Transient 1D problem} \label{ssec:transient_numerical_example}
In this subsection, we demonstrate the effect of the temporal variance weighting hyper-parameter $\vec{\alpha}_t$. Consider a transient generalization of the advection-diffusion equation from the previous subsection. The high-fidelity model is 
\begin{align*}
& \frac{du}{dt} -\frac{d^2 u}{dx^2} + \frac{du}{dx} = z \qquad & \text{on } (0,1) \\
& \frac{du}{dx}(0) =\frac{du}{dx}(1) = 0 \\
& u(0,x)= 0,
\end{align*}
which constrains an optimization problem to determine an optimal source function $z:[0,1] \to \R$. The low-fidelity model is given by the diffusion equation 
\begin{align*}
&  \frac{du}{dt} -\frac{d^2 u}{dx^2} = z \qquad & \text{on } (0,1) \\
& \frac{du}{dx}(0) =\frac{du}{dx}(1) = 0 \\
& u(0,x)= 0,
\end{align*}
which omits the effect of advection. For this example, we assume that data $\{\z_\ell, \vec{d}_\ell \}_{\ell=1}^2$ is available to the hyper-parameter initialization algorithms.

\begin{figure}[h]
\centering
  \includegraphics[width=0.35\textwidth]{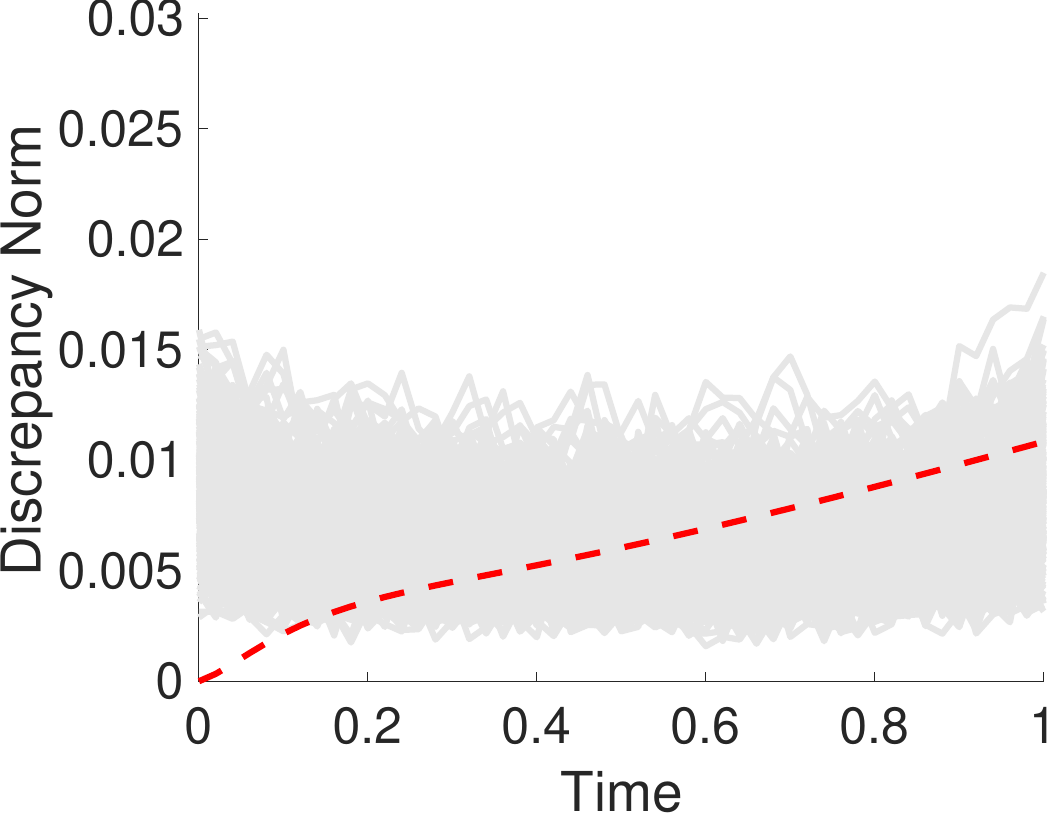}
    \includegraphics[width=0.35\textwidth]{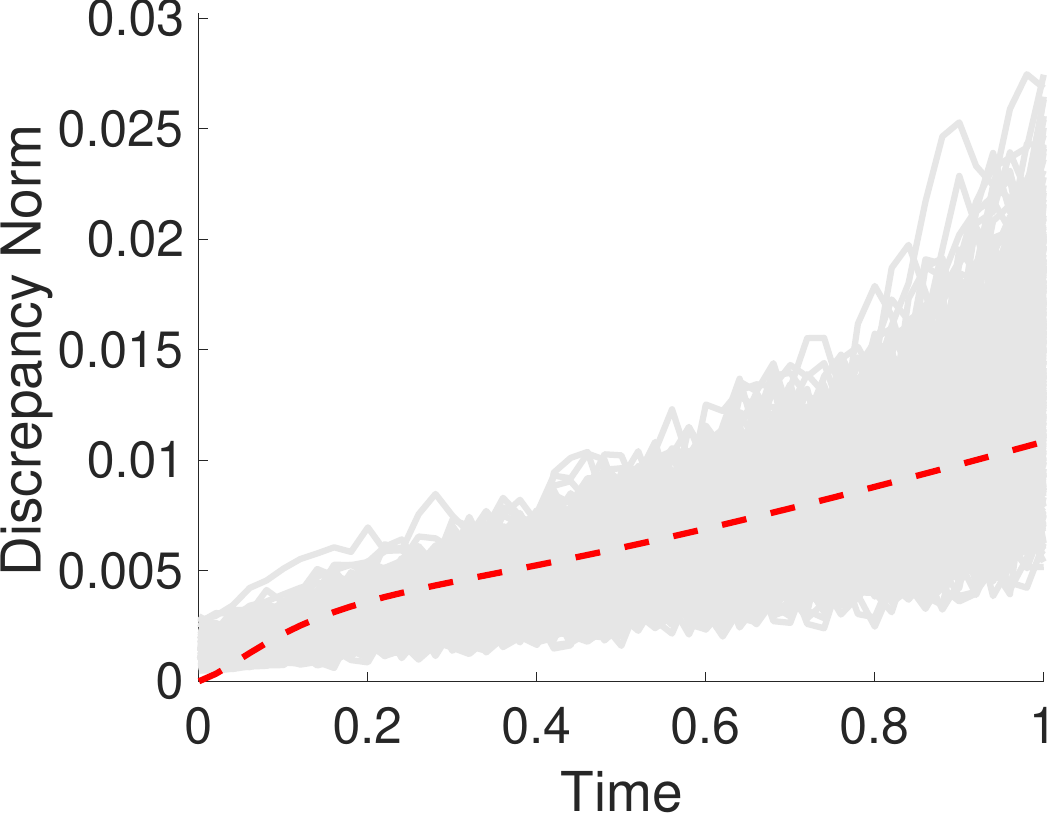}
    \caption{Spatial norm of prior discrepancy samples as a function of time. The left and right panels correspond to the two cases where $\vec{\alpha}_t = \vec{e}$ and $\vec{\alpha}_t$ is computed using the approach described in Subsection~\ref{ssec:time_weighting}, respectively. Each solid gray curve corresponds to a prior discrepancy sample and the broken red curve (appearing in both left and right panels) corresponds to the spatial norm of the discrepancy data $\vec{d}_1$ as a function of time.}
  \label{fig:trans_var}
\end{figure}

To highlight the role of the variance weighting hyper-parameter $\vec{\alpha}_t$, we consider two cases. The first is when $\vec{\alpha}_t = \vec{e}$ is taken as the default vector of $1$'s, thus not weighting the time nodes. In the second case, we initialize $\vec{\alpha}_t$ using the approach described in Subsection~\ref{ssec:time_weighting}. Figure~\ref{fig:trans_var} displays the norm (on the spatial domain $\Omega=[0,1]$) of prior discrepancy samples evaluated at $\z=\ztilde$, as a function of time. That is, we generate samples $\d(\ztilde, \T_\t \vec{\omega}_\t) \in \R^{n_s n_t}$ and compute the spatial norm at each time step. The left and right panels correspond to the two cases where $\vec{\alpha}_t = \vec{e}$ and $\vec{\alpha}_t$ is computed using the approach described in Subsection~\ref{ssec:time_weighting}, respectively, with the spatial norm of $\vec{d}_1$, as a function of time, overlaid on both plots. 
  
  Figure~\ref{fig:trans_var} demonstrates the role of $\vec{\alpha}_t$ weighting time so that the prior discrepancy has a small magnitude at time $t=0$ and increases as a function of time, which reflects the effect of compounding errors accumulating over the time steps. In comparison, if $\vec{\alpha}_t=\vec{e}$, we observe magnitudes that are too large in the early time step (all samples have discrepancy norms greater than the observed data for $t \in (0,0.2)$, but then subsequently most samples have magnitudes less than the observed data at time $t=1$. 
  
  \subsection{Stationary 2D problem}  \label{ssec:num_2d}
  In this subsection, we demonstrate our framework to visualize prior samples in higher dimensions, which is crucial to verify or adjust the hyper-parameters in problems of practical interest. To this end, we consider a stationary 2D analog the previous advection-diffusion examples. Specifically, the high-fidelity model is 
\begin{align*}
& -\Delta u+ \vec{v} \cdot \nabla u = z \qquad & \text{on } \Omega \\
& \nabla u \cdot \vec{n} = 0 & \text{on } \partial \Omega_n \\
& u= 0 & \text{on } \partial \Omega_d,
\end{align*}
where $\Omega=(0,1)^2$, $ \partial \Omega_n=\{0\} \times (0,1) \cup \{1\} \times (0,1) \cup (0,1) \times \{1\}$ is the Neumann boundary, and $\partial \Omega_d = (0,1) \times \{0\}$ is the Dirichlet boundary, $\vec{n}$ is the normal vector to the boundary, the optimization variable is the source function $z:\Omega \to \R$, and the velocity field is given by $\vec{v}=(5,5)$. The low-fidelity model is given by the diffusion equation 
\begin{align*}
& -\Delta u = z \qquad & \text{on } \Omega \\
& \nabla u \cdot \vec{n} = 0 & \text{on } \partial \Omega_n \\
& u= 0 & \text{on } \partial \Omega_d.
\end{align*}

We use a single high-fidelity evaluation to form a dataset $\{\z_1,\vec{d}_1\}$ and initialize the prior hyper-parameters $\alpha_\u,\beta_\u,\alpha_\z$, and $\beta_\z$. Following the framework presented in Section~\ref{ssec:viz_framework}, we generate a dataset 
$$\{ \Delta \z_k, \d(\ztilde,\t_i), \d(\ztilde + \Delta \z_k,\t_i) \}_{i=1,k=1}^{Q,P}.$$
Analysis of $\d(\ztilde,\t_i)$ verifies that $\alpha_\u$ and $\beta_\u$ were initialized with appropriate values. We omit details of this $\{\alpha_\u,\beta_\u\}$ analysis and focus instead on analysis of $\{\alpha_\z,\beta_\z\}$. 

\begin{figure}[!ht]
\centering
  \includegraphics[width=0.9\textwidth]{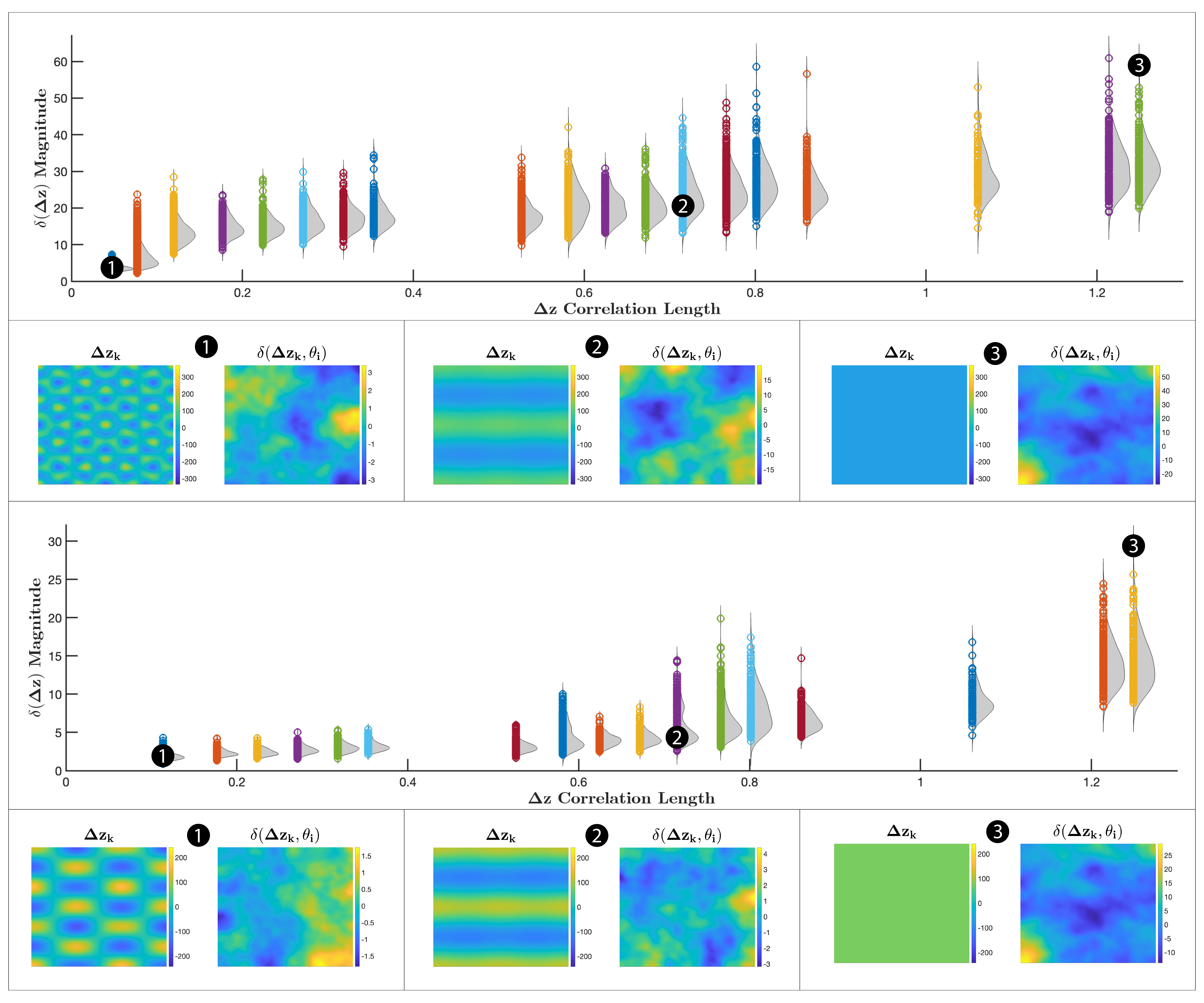}
    \caption{Depiction of an interactive visualization framework to examine the prior model discrepancy. The top scatter plot depicts all samples $\{ \Delta \z_k,\d(\ztilde+\Delta \z_k,\t_i) -\d(\ztilde,\t_i)\}_{i=1,k=1}^{Q,P}$, where the horizontal axis corresponds to the empirically estimated correlation length of $\Delta \z_k$ (with values binned for ease of visualization) and the vertical axis corresponds to the maximum absolute magnitude of $\d(\Delta \z_k,\t_i) =\d(\ztilde+\Delta \z_k,\t_i) -\d(\ztilde,\t_i)$. The horizontal axis labels indicate the average correlation length within that bin. Due to the large number of points within each bin (along a vertical column), we use a density plot to indicate the distribution of the magnitude of $\d(\Delta \z_k,\t_i)$. Upon clicking on a point, two surface plots appear corresponding to $\d(\Delta \z_k,\t_i)$ and $\Delta \z_k$. Below each scatter plot, we show these plots for three selections, with numbered circles indicating which point in the scatter plot the surfaces correspond to. The bottom scatter plot corresponds to recomputing the samples after adjusting the hyper-parameters.}
  \label{fig:interactive_viz_snapshot}
\end{figure}

Figure~\ref{fig:interactive_viz_snapshot} showcases how we use the interactive visualization to help adjust and verify $\{\alpha_{z},\beta_{z}\}$ through two sets of hyper-parameter instances. Each instance consists of an overview plot and three samples selected from the overview plot for detailed inspections. The overview plot shows a projection of all samples $\{ \Delta \z_k, \d(\Delta \z_k,\t_i)\}_{i=1,k=1}^{Q,P}$ on to a 2D space, where $\d(\Delta \z_k,\t_i) =\d(\ztilde+\Delta \z_k,\t_i) -\d(\ztilde,\t_i)$. The 2D projection is characterized by the maximum magnitude of $\d(\Delta z_k,\t_i)$ and the empirical estimation of the correlation length of $\Delta \z_k$. The projection clusters the samples into bins of correlation lengths.

It is difficult to interpret the distribution of the samples within each correlation length bin since many samples are projected into the same bin. We augment the scatter plot with a density plot (per bin) to better illustrate the distribution. The correlation length is a measure of $\Delta z_k$'s smoothness, so the leftmost samples in the scatter plot correspond to the high-frequency (most non-smooth) $\Delta \z_k$'s. We select three samples from the leftmost, middle, and rightmost bins. The plots of $\Delta z_k$ show decreasing frequency (increasing smoothness) from left to right. Since the model discrepancy prior is smoothing, the high-frequency $\Delta z_k$'s lead to lower magnitude changes in the discrepancy $\d$. The hyper-parameter $\beta_\z$ controls the relationship between the frequency of $\Delta z_k$ and the magnitude of $\d(\Delta z_k,\t)$. Taking a large $\beta_\z$ will result in $\d$ only being sensitive to the lowest frequency $\Delta z_k$'s, thus restricting the prior to discrepancies that are nearly constant (as functions of $\z$). On the other hand, taking a small $\beta_\z$ will result in $\d$ exhibiting non-trivial sensitivity to high-frequency $\Delta z_k$'s. 

By examining the density centers of each correlation length bin, we observe a relatively large maximum magnitude of $\d(\Delta z_k,\t_i)$ in the small correlation length bins. This observation is based on knowledge of the physical processes in the high- and low-fidelity models. Specifically, the process of diffusion dampens high-frequency sources, and consequently $\d(\Delta z_k,\t_i)$ should be small for high-frequency $\Delta z_k$'s. The magnitude observed in small correlation length bins indicates that $\beta_\z$ is too small. 

The misspecification of $\beta_\z$ was because the data used in the initialization algorithm was a spatially localized source. The empirical correlation length estimation assumes that the data is sampled from a Gaussian random field. However, spatially localized fields are not modeled well via Gaussian random fields, and consequently, the estimate of $\beta_\z$ was too small. The prior is adjusted by increasing $\beta_\z$ and repeating the visualization to confirm that the prior samples have characteristics that are commensurate with what is expected from the physical principles. 

The samples of $\delta(\Delta z_k,\theta_i)$ are about double the magnitude that we expect. Such an assessment is based on physical intuition about what is lacking in the low-fidelity model to prescribe an appropriate discrepancy magnitude. The misspecification of $\alpha_\z$ is due to the many approximations that are used to initialize it. We divide $\alpha_\z$ by $4.0$ to correct the discrepancy magnitude. 

The second hyper-parameter instance in Figure~\ref{fig:interactive_viz_snapshot} shows the result after adjusting $\{\alpha_{\z},\beta_{\z}\}$. As in the case before adjustment, we select three samples (from the leftmost, middle, and rightmost correlation length bins) for inspection. Observe that the magnitude of the lowest-frequency $\Delta z_k$ (the rightmost bin) is approximately halved, as a result of changing $\alpha_\z$. Furthermore, the magnitude of the discrepancy decays more rapidly as a function of the correlation length, as a result of changing $\beta_\z$. For instance, around correlation length $0.3$, the median magnitude was previously greater than $15$, but after the update it is less than $5$. Changing $\alpha_\z$ accounts for dividing the magnitude by $2$, while changing $\beta_\z$ accounts for the remaining decrease as high-frequency modes are damped by a larger $\beta_\z$. 

The combination of the overview plot and interactive inspection relieved us from the tedious process of trying to assess the quality of hyper-parameters by inspecting randomly chosen samples from $\{ \Delta \z_k,\d(\ztilde+\Delta \z_k,\t_i) -\d(\ztilde,\t_i)\}_{i=1,k=1}^{Q, P}$. Specifically, in addition to the quick and intuitive overview, the interactions on the overview plot allow us to quickly inspect outlying samples, such as the leftmost and the rightmost samples, to infer bounds on the samples. If those samples on the boundary of the projected space show unexpected outcomes, we can quickly identify the cause and make necessary adjustments.

  \section{Conclusion} \label{sec:conclusion}
Hyper-differential sensitivity analysis with respect to model discrepancy has considerable potential to enable approximate optimization with quantified uncertainties for many applications where currently only low-fidelity optimization is possible. However, to provide useful uncertainty estimates, it is crucial to specify the model discrepancy prior distribution appropriately. As was demonstrated in Subsection~\ref{ssec:stationary_1D_example}, misspecification of the model discrepancy prior may lead to an optimal solution posterior distribution that is too narrow, thus giving unwarranted confidence in the solution. On the other hand, if the model discrepancy prior distribution is not sufficiently informative, then the posterior will have large uncertainty in directions not informed by the high-fidelity data. Consequently, propagating this posterior discrepancy through the optimization problem will yield large uncertainty in the optimization solution such that it fails to be useful.

Specifying the prior model discrepancy is challenging for large-scale and coupled systems due to the overwhelming volume of information that must be processed in order to formulate the prior appropriately. For instance, in a coupled system with five transient states in three spatial dimensions (e.g., a fluid flow model with three velocity components, a pressure component, and a temperature component), there are over 30 prior hyper-parameters that must be understood through visualization of many transient vector fields. This article provides a pragmatic approach to specify the prior through a combination of algorithmic initialization of the prior hyper-parameters and a visualization framework to understand their influence on the prior discrepancy. Importantly, our approach is scalable in both the computational complexity of the algorithmic initialization and in the visual complexity of the prior discrepancy analysis. Our mathematical analysis and numerical demonstrations coupled with visualization provide insight to guide hyper-parameter specification that would otherwise be intractable due to the volume of data generated from prior sampling.

Traditional Bayesian approaches to hyper-parameter selection may be considered as alternatives to our approach. For instance, choosing hyper-parameters that maximize the data likelihood. However, there are several important problem characteristics that motivate our approach. We focus on problems where the number of high-fidelity model evaluations $N$ is small. In the case where $N=1$, there is no information in the data to inform $\alpha_\z$ using data likelihood maximization. For $N>1$ but much less than the dimension of $\z$, there will be many local extrema in the data likelihood maximization as a result of the hyper-parameters having confounding effects. On the other hand, for any $N \ge 1$, our initialization approach leverages the physical interpretation of the hyper-parameters to sequentially select the hyper-parameters. The only solution variability comes from estimating expected values. Hybrid approaches can also be explored. For example, our proposed method can be combined with data likelihood maximization to (i) reduce the number of local minima by fixing a subset of hyper-parameters, (ii) inform the prior distribution of the hyper-parameters, or (iii) serve as an initialization for the data likelihood optimization algorithm.

This article focused on Gaussian priors, which are reasonable given the characteristics of the mathematical object being estimated (a model discrepancy operator). However, we have not discussed enforcement of physics-based laws. For instance, if both the high- and low-fidelity models conserve mass, then the model discrepancy should also conserve mass. For physical laws that can be expressed as linear operators, the discrepancy may be parameterized as $\vec{C} \d(\z,\t)$, where $\vec{C}$ is the linear operator enforcing the conservation law. Consequently, the proposed approach of prior hyper-parameter specification and discrepancy calibration may still be possible. However, many constraints may require the use of nonlinear operators or inequality constraints, which prohibit much of our analysis. Future work should explore which constraints can be imposed on the prior.
 
 Lastly, we note that this work may provide a foundation for uncertainty quantification beyond our HDSA-MD framework. In particular, our approach to prior modeling may be useful in operator learning. Calibration of $\d$ is an example of learning an affine operator. Future work may explore extensions to nonlinear operator learning such as uncertainty quantification on a linearization of the learned operator, for instance, the output layer of a Neural Operator.

\acknowledgements
 
  This work was supported by the U.S. Department of Energy,
Office of Science, Office of Advanced Scientific Computing Research Field Work
Proposal Number 23-02526 (Ab-initio Visualization for Innovative Science). Li, Ouermi, and Johnson were partially supported by the Intel OneAPI CoE. Sandia National Laboratories is a multi-mission laboratory managed and operated by National Technology \& Engineering Solutions of Sandia, LLC (NTESS), a wholly owned subsidiary of Honeywell International Inc., for the U.S. Department of Energy’s National Nuclear Security Administration (DOE/NNSA) under contract DE-NA0003525. This written work is authored by an employee of NTESS. The employee, not NTESS, owns the right, title and interest in and to the written work and is responsible for its contents. Any subjective views or opinions that might be expressed in the written work do not necessarily represent the views of the U.S. Government. The publisher acknowledges that the U.S. Government retains a non-exclusive, paid-up, irrevocable, world-wide license to publish or reproduce the published form of this written work or allow others to do so, for U.S. Government purposes. The DOE will provide public access to results of federally sponsored research in accordance with the DOE Public Access Plan.
 
 All data used in this article was generated by sampling from Gaussian random fields and solving partial differential equations. The source code is not publicly available, but data may be shared upon reasonable request to the corresponding author.
 
\appendix

\section{Proof of Theorem 1} \label{proof_thm:fit_solutions}

\begin{proof}
Let
\begin{eqnarray*}
\mathbf{A}_\ell =
\left( \begin{array}{cc}
\I_{n_u} & \I_{n_u} \otimes \z_\ell^T \M_\z
\end{array} \right) \in \R^{n_u \times n_\theta}
\end{eqnarray*}
so that $\d(\z_\ell,\t) = \A_\ell \t$, and
\begin{eqnarray*}
\mathbf{A} =
\left( \begin{array}{c}
\A_1 \\
\A_2 \\
\vdots \\
\A_N
\end{array} \right) \in \R^{n_u N \times n_\theta} 
\qquad 
\text{and}
\qquad 
\mathbf{d} =
\left( \begin{array}{c}
\vec{d}_1 \\
\vec{d}_2 \\
\vdots \\
\vec{d}_N
\end{array} \right) \in \R^{n_u N} .
\end{eqnarray*}
Then $\d(\z_\ell,\t)=\vec{d}_\ell$, $\ell=1,2,\dots,N$, is equivalent to $\A \t = \vec{d}$. Since $\{\z_\ell\}_{\ell=1}^N$ are linearly independent and $\M_\z$ is a nonsingular matrix, $\{ \M_\z \z_\ell\}_{\ell=1}^N$ are linearly independent, which implies that the rows of $\A$ are linearly independent and hence $\text{rank}(\A) = n_u N$. Furthermore, $\vec{d} \in \R^{n_u N}$ and hence $\vec{d}$ is in the range space of $\A$. The result follows with $\t^\star=\A^\dagger \vec{d}$, where $\A^\dagger$ is the pseudo-inverse of $\A$, and $\Theta_{\text{interp}}$ is the null space of $\A$.
\end{proof}

\section{Proof of Theorem 2} \label{proof_thm:trace_class_cov}

 \begin{proof}
 Observe that 
 \begin{align*}
\text{Tr}_{\M_\t}(\W_\t^{-1}) = (1+\ztilde^T \M_\z \W_\z^{-1} \M_\z \ztilde + \text{Tr}_{\M_\z}(\W_\z^{-1})) \text{Tr}_{\M_\u}(\W_\u^{-1}).
\end{align*}
Furthermore, $\ztilde^T \M_\z \W_\z^{-1} \M_\z \ztilde = \vert \vert \E_\z^{-1} \M_\z \ztilde \vert \vert_{\M_\z}^2$, and hence
 \begin{align*}
\lim\limits_{\U_h \to \U, \Z_h \to \Z} \ztilde^T \M_\z \W_\z^{-1} \M_\z \ztilde = \lim\limits_{ \Z_h \to \Z}\vert \vert \E_\z^{-1} \M_\z \ztilde \vert \vert_{\M_\z}^2 = \| \mathcal E_z^{-1} \tilde{z} \|_{\Z}^2,
\end{align*}
where $\mathcal E_z$ is the operator discretized by $\E_\z$. Since $\lim\limits_{\U_h \to \U} \text{Tr}_{\M_\u}(\W_\u^{-1})$ and $\lim\limits_{\Z_h \to \Z} \text{Tr}_{\M_\z}(\W_\z^{-1})$ exist and are finite, it follows that $\lim\limits_{\U_h \to \U, \Z_h \to \Z} \text{Tr}_{\M_\t}(\W_\t^{-1})$ exists and is finite.
 \end{proof}

\section{Proof of Theorem 3} \label{proof_thm:exp_delta_diff}

 \begin{proof}
 Computing the squared norm of $\d(\ztilde + \Delta \z,\vec{T}_\t \vec{\omega}_\t) - \d(\ztilde ,\vec{T}_\t \vec{\omega}_\t) = \T_\u \L(\vec{\omega}_\t) \T_\z^T \M_\z \Delta \z $ and rearranging the multiplications gives
\begin{align*}
\| \d(\ztilde + \Delta \z,\vec{T}_\t \vec{\omega}_\t) - \d(\ztilde ,\vec{T}_\t \vec{\omega}_\t) \|_{\M_\u}^2 = \|  (\T_\u \otimes  \Delta \z^T  \M_\z \T_\z) \vec{\ell}(\vec{\omega}_\t) \|_{\M_\u}^2,
\end{align*}
where $\vec{\ell}:\R^{n_\theta} \to \R^{n_u n_z}$ maps vectors of dimension $n_\theta=n_u (n_z+1)$ to their $n_u+1$ through $n_\theta$ components, i.e., $\vec{\ell}(\vec{\omega}_\t)$ corresponding to vectorizing the entries of $\L(\vec{\omega}_\t)$.

Hence, $\| \d(\ztilde + \Delta \z,\L_\t \vec{\omega}_\t) - \d(\ztilde ,\L_\t \vec{\omega}_\t) \|_{\M_\u}^2$ is distributed as a weighted sum of chi-squared random variables, where the weights are defined by the eigenvalues of 
\begin{align*}
(\T_\u \otimes  \Delta \z^T  \M_\z \T_\z)^T \M_\u (\T_\u \otimes  \Delta \z^T  \M_\z \T_\z) = \T_\u^T \M_\u \T_\u \otimes \T_\z^T \M_\z \Delta \z \Delta \z^T \M_\z \T_\z.
\end{align*}
If follows that the expected value of the discrepancy difference norm squared is
\begin{align*}
\mathbb{E}_{\vec{\omega_\t}} \left[ \| \d(\ztilde + \Delta \z,\L_\t \vec{\omega}_\t) - \d(\ztilde ,\L_\t \vec{\omega}_\t) \|_{\M_\u}^2 \right] & =  \text{Tr}_2( \T_\u^T \M_\u \T_\u \otimes \T_\z^T \M_\z \Delta \z \Delta \z^T \M_\z \T_\z) \\
& = \text{Tr}_{\M_\u}(\W_\u^{-1}) \left( \Delta \z^T \M_\z \W_\z^{-1} \M_\z \Delta \z \right),
\end{align*}
where the latter equality comes from writing the trace of a Kronecker product as the product of the traces, and commuting the order of matrix multiplication within the trace operation. Recalling that $\mathbb{E}_{\vec{\omega_\t}} \left[ \| \d(\ztilde ,\L_\t \vec{\omega}_\t) \|_{\M_\u}^2 \right] = \text{Tr}_{\M_\u}(\W_\u^{-1})$ and that $\W_\z^{-1} = \alpha_\z \E_\z^{-1} \M_\z \E_\z^{-1}$, we have
\begin{align*}
\frac{\mathbb{E}_{\vec{\omega_\t}} \left[ \| \d(\ztilde + \Delta \z,\L_\t \vec{\omega}_\t) - \d(\ztilde ,\L_\t \vec{\omega}_\t) \|_{\M_\u}^2 \right] }{\mathbb{E}_{\vec{\omega_\t}} \left[ \| \d(\ztilde ,\L_\t \vec{\omega}_\t) \|_{\M_\u}^2 \right] } = \alpha_\z \Delta \z^T \M_z \E_\z^{-1} \M_\z \E_\z^{-1} \M_\z \Delta \z .
\end{align*}
 \end{proof}
 
\section{Proof of Theorem 4} \label{proof_thm:exp_delta_z}

\begin{proof}
Samples of $\Delta \z$ are given by
\begin{align*}
\Delta \z =\| \ztilde \|_{\M_\z} \frac{ \E_\z^{-1} \M_\z^{\frac{1}{2}} \vec{\omega}_\z }{ \| \E_\z^{-1} \M_\z^{\frac{1}{2}} \vec{\omega}_\z \|_{\M_\z} },
\end{align*}
where $\vec{\omega}_\z \in \R^{n_z}$ is a standard normal Gaussian random vector. It follows that 
\begin{align*}
\Delta \z^T \M_z \E_\z^{-1} \M_\z \E_\z^{-1} \M_\z \Delta \z  = \| \ztilde \|_{\M_\z}^2 \frac{ \vec{\omega}_\z^T \M_\z^{\frac{1}{2}} \E_\z^{-1} \M_z \E_\z^{-1} \M_\z \E_\z^{-1} \M_\z \E_\z^{-1} \M_\z^{\frac{1}{2}} \vec{\omega}_\z} { \vec{\omega}_\z^T \M_\z^{\frac{1}{2}} \E_\z^{-1} \M_z \E_\z^{-1} \M_\z^{\frac{1}{2}} \vec{\omega}_\z} .
\end{align*}

Leveraging the spectral decomposition $\E_\z = \M_\z \V \vec{\Lambda} \V^T \M_\z$ and the associated decomposition of its inverse $\E_\z^{-1} = \V \vec{\Lambda}^{-1} \V^T$. It follows that
\begin{align*}
\E_\z^{-1} \M_z \E_\z^{-1} = \V \vec{\Lambda}^{-2} \V^T \qquad \text{and} \qquad  \E_\z^{-1} \M_z \E_\z^{-1} \M_\z \E_\z^{-1} \M_\z \E_\z^{-1} = \V \vec{\Lambda}^{-4} \V^T.
\end{align*}
Hence,
\begin{align*}
\Delta \z^T \M_z \E_\z^{-1} \M_\z \E_\z^{-1} \M_\z \Delta \z  = \| \ztilde \|_{\M_\z}^2  \frac{ \vec{\omega}_\z^T \M_\z^{\frac{1}{2}} \V \vec{\Lambda}^{-4} \V^T \M_\z^{\frac{1}{2}} \vec{\omega}_\z} { \vec{\omega}_\z^T \M_\z^{\frac{1}{2}} \V \vec{\Lambda}^{-2} \V^T \M_\z^{\frac{1}{2}} \vec{\omega}_\z} .
\end{align*}
Since $\M_\z^\frac{1}{2} \V$ is an orthogonal matrix, $\vec{\nu}_\z = \M_\z^{\frac{1}{2}} \V \vec{\omega}_\z$ is a standard normal random vector. Accordingly, 
\begin{align*}
\mathbb E_{\Delta \z} \left[ \Delta \z^T \M_z \E_\z^{-1} \M_\z \E_\z^{-1} \M_\z \Delta \z \right] = \| \ztilde \|_{\M_\z}^2 \mathbb E_{\vec{\nu}_\z} \left[ \frac{ \vec{\nu}_\z^T \vec{\Lambda}^{-4} \vec{\nu}_\z} { \vec{\nu}_\z^T  \vec{\Lambda}^{-2} \vec{\nu}_\z} \right].
\end{align*}
\end{proof}

\section{Empirical correlation length estimation} \label{appendix_corr_length}

 This appendix summarizes how the correlation length of a Gaussian random field is estimated from data. For generality, we assume that the data is denoted by $\vec{f} \in \R^m$, which corresponds to evaluations of a function $f:\Omega \to \R$, where $\Omega \subset \R^s$ is a domain with either $s=1,2,$ or $3$. This formulation is applicable for the state and optimization variables with dependence on time or space. Furthermore, we assume that nodal values are provided along with $\vec{f}$, i.e., the nodal values are contained in $\mathbf{X} \in \R^{m \times s}$, so that the $i^{th}$ component of $\vec{f}$ corresponds to evaluating the function $f$ for the $i^{th}$ row of $\mathbf{X}$.
 
 Since we are calibrating hyperparamters for a spatially homogenous Gaussian prior, we treat $\vec{f}$ as a realization of a Gaussian random field. We seek to compute the correlation length $\kappa \in (0,\infty)$, which is defined as the distance for which the correlation between evaluations of $f$ at points $\vec{x},\vec{y}$, such that $\| \vec{x} - \vec{y} \|_2 = \kappa$, have a correlation of $0.1$. 
 
 To compute the correlation length from data, we use evaluations of $f$ at various inputs $\vec{x},\vec{y}$ to estimate the correlation for various distances $\| \vec{x} - \vec{y} \|_2$ and estimate $\kappa$ as the distance for which we achieve a correlation estimate of approximately $0.1$. To simplify this computation, we use the nodal data $\mathbf{X}$ and function evaluation data $\vec{f}$ to construct an interpolating function $\hat{f}$ such that $\hat{f} \approx f$. We assume that the nodal data is sufficiently fine that interpolation error is negligible relative to others errors arising from assumptions and sampling.
 
 Algorithm~\ref{alg:correlation_length_estimation} outlines our approach to estimate the correlation length. We present the approach in dimension $s=1$ on the domain $\Omega=(0,1)$ for simplicity, but the algorithm generalizes to higher dimensions.
 
 \begin{algorithm}[ht!!]
	\caption{Correlation length estimation}
	\begin{algorithmic} [1] \label{alg:correlation_length_estimation}
		\STATE Input: $\vec{x}, \vec{f} \in \R^m$
		\STATE Compute the mean $\overline{f} = \frac{1}{m} \sum\limits_{i=1}^m f_i$
		\STATE Compute the variance $\sigma^2(f) = \frac{1}{m-1} \sum\limits_{i=1}^m (f_i-\overline{f})^2$
		\STATE Construct an interpolating function $\hat{f}$
		\STATE Set a distance increment $\Delta \kappa \in (0,1)$
		\STATE Initialize the correlation $\rho=1$
		\STATE Initialize the correlation length $\kappa=0$
		\WHILE{$\rho > 0.1$}
		\STATE Increment $\kappa = \kappa + \Delta  \kappa$
		\STATE Compute $\vec{f}^{\kappa} \in \R^m$ by evaluating $\hat{f}$ at $x_i+\Delta  \kappa$, $i=1,2,\dots,m$
		\STATE Compute the covariance $c_\kappa = \frac{1}{m-1} \sum\limits_{i=1}^m (f_i-\overline{f}) (f_i^{\kappa}-\overline{f})$
		\STATE Set the correlation $\rho = \frac{c_\kappa}{\sigma^2(f)}$
		\ENDWHILE
		\STATE Return: Correlation length $\kappa$
		\STATE Note: For some $i$ in lines 10-12, $x_i \pm \Delta  \kappa \notin \Omega$, and hence evaluating $\hat{f}$ is not appropriate. These datapoints are omitted and the normalization factor $m-1$ in Line 12 is adjusted accordingly.
	\end{algorithmic}
\end{algorithm}
 
 \section{Determining perturbation directions} \label{appendix_z_pert}
To determine the directions $\Delta \z$ that maximize the prior discrepancy variation $\d(\ztilde+\Delta \z,\vec{T}_\t \vec{\omega}_\t)-\d(\ztilde,\vec{T}_\t \vec{\omega}_\t)$, observe that~\eqref{eqn:prior_discrepancy_samples} implies that
 \begin{align*}
\d(\ztilde+\Delta \z,\vec{T}_\t \vec{\omega}_\t)-\d(\ztilde,\vec{T}_\t \vec{\omega}_\t) = \T_\u \vec{L}(\vec{\omega}_\t) \T_\z^T \M_\z \Delta \z .
 \end{align*}
We seek the perturbation directions $\Delta \z$ that maximize $\| \T_\z^T \M_\z \Delta \z \|_2^2$. Recalling that $\T_\z = \E_\z^{-1} \M_\z^\frac{1}{2}$, we have
\begin{align*}
\| \T_\z^T \M_\z \Delta \z \|_2^2= \Delta \z^T \M_\z \E_\z^{-1} \M_\z  \E_\z^{-1} \M_\z \Delta \z.
\end{align*}
Since $\E_\z^{-1} = \V \vec{\Lambda}^{-1} \V^T$ and $\V^T \M_\z \V = \vec{I}$,
\begin{align*}
\| \T_\z^T \M_\z \Delta \z \|_2^2= \Delta \z^T \M_\z \V \vec{\Lambda}^{-2} \V^T \M_\z \Delta \z.
\end{align*}
Consequently, the maximum perturbation directions are given by the eigenvectors corresponding to the smallest eigenvalues of $\E_\z$, or equivalently, the largest eigenvalues of $\E_\z^{-1}$.

	\bibliographystyle{IJ4UQ_Bibliography_Style}

	\bibliography{dasco}

\begin{thebibliography}{10}

\bibitem{hart_bvw_cmame}
Hart, J. and {van Bloemen Waanders}, B., Hyper-differential sensitivity
  analysis with respect to model discrepancy: Optimal solution updating, {\em
  Computer Methods in Applied Mechanics and Engineering}, 412, 2023.

\bibitem{hart_bvw_mods}
Hart, J. and {van Bloemen Waanders}, B., Hyper-differential sensitivity
  analysis with respect to model discrepancy: Optimal solution posterior
  sampling, {\em Foundations of Data Science}, 7(1):99--133, 2025.

\bibitem{multifidelity_review_peherstorfer}
Peherstorfer, B., Willcox, K., and Gunzburger, M., Survey of multifidelity
  methods in uncertainty propagation, inference, and optimization, {\em SIAM
  Review}, 60(3):550--591, 2018.

\bibitem{multifidelity_quasinewton_bryson}
Bryson, D.E. and Rumpfkeil, M.P., Multifidelity {Q}uasi-{N}ewton method for
  design optimization, {\em AIAA JOURNAL}, 56(10), 2018.

\bibitem{trmm_lewis}
Alexandrov, N., {Jr.}, J.D., Lewis, R., and Torczon, V., A trust-region
  framework for managing the use of approxima- tion models in optimization,
  {\em Structural Optimization}, 15:16--23, 1998.

\bibitem{biegler_rom_opt}
Agarwal, A. and Biegler, L., A trust-region framework for constrained
  optimization using reduced order modeling, {\em Optimization and
  Engineering}, 14:3--35, 2013.

\bibitem{willcox_multifi_opt_2012}
March, A. and Willcox, K., Provably convergent multifidelity optimization
  algorithm not requiring high-fidelity derivatives, {\em AIAA Journal},
  50(5):1079--1089, 2012.

\bibitem{ohagan2001}
Kennedy, M.C. and O'Hagan, A., Bayesian calibration of computer models, {\em
  Journal of the Royal Statistical Society}, 63(3):425--464, 2001.

\bibitem{ohagan2002}
Oakley, J. and O'Hagan, A., Bayesian inference for the uncertainty distribution
  of computer model outputs, {\em Biometrika}, 89(4):769--784, 2002.

\bibitem{ohagan2014}
Brynjarsd{\'o}ttir, J. and O'Hagan, A., Learning about physical parameters: the
  importance of model discrepancy, {\em Inverse Problems}, 30(11), 2014.

\bibitem{Higdon_2008}
Higdon, D., Gattiker, J., Williams, B., and Rightley, M., Computer model
  calibration using high-dimensional output, {\em Journal of the American
  Statistical Association}, 103(482):570--583, 2008.

\bibitem{Arendt_2012}
Arendt, P.D., Apley, D.W., and Chen, W., Quantification of model uncertainty:
  Calibration, model discrepancy, and identifiability, {\em Journal of
  Mechanical Design}, 134, 2012.

\bibitem{Maupin}
Maupin, K.A. and Swiler, L.P., Model discrepancy calibration across
  experimental settings, {\em Reliability Engineering \& System Safety}, 200,
  2020.

\bibitem{Ling_2014}
Ling, Y., Mullins, J., and Mahadevan, S., Selection of model discrepancy priors
  in bayesian calibration, {\em Journal of Computational Physics},
  276:665--680, 2014.

\bibitem{Gardner_2021}
Gardner, P., Rogers, T., Lord, C., and Barthorpe, R., Learning model
  discrepancy: A gaussian process and sampling-based approach, {\em Mechanical
  Systems and Signal Processing}, 152(107381), 2021.

\bibitem{Vogel_99}
Vogel, C.R., Sparse matrix computations arising in distributed parameter
  identification, {\em SIAM J. Matrix Anal. Appl.}, pp. 1027--1037, 1999.

\bibitem{Archer_01}
Ascher, U.M. and Haber, E., Grid refinement and scaling for distributed
  parameter estimation problems, {\em Inverse Problems}, 17:571--590, 2001.

\bibitem{Haber_01}
Haber, E. and Ascher, U.M., Preconditioned all-at-once methods for large,
  sparse parameter estimation problems, {\em Inverse Problems}, 17:1847--1864,
  2001.

\bibitem{Vogel_02}
Vogel, C.R., {\em Computational Methods for Inverse Problems}, SIAM Frontiers
  in Applied Mathematics Series, 2002.

\bibitem{Biegler_03}
L.~T.~Biegler, O.~Ghattas, M.H. and van Bloemen~Waanders, B. (Eds.), {\em
  Large-Scale PDE-Constrained Optimization}, Vol.~30, Springer-Verlag Lecture
  Notes in Computational Science and Engineering, 2003.

\bibitem{Biros_05}
Biros, G. and Ghattas, O., Parallel {L}agrange-{N}ewton-{K}rylov-{S}chur
  methods for {PDE}-constrained optimization. {Parts I-II}, {\em SIAM J. Sci.
  Comput.}, 27:687-- 738, 2005.

\bibitem{Laird_05}
Laird, C.D., Biegler, L.T., van Bloemen~Waanders, B., and Bartlett, R.A., Time
  dependent contaminant source determination for municipal water networks using
  large scale optimization, {\em ASCE J. Water Res. Mgt. Plan.}, pp. 125--134,
  2005.

\bibitem{Hintermuller_05}
Hintermuller, M. and Vicente, L.N., Space mapping for optimal control of
  partial differential equations, {\em SIAM J. Opt.}, 15:1002--1025, 2005.

\bibitem{Hazra_06}
Hazra, S.B. and Schulz, V., Simultaneous pseudo-timestepping for aerodynamic
  shape optimization problems with state constraints, {\em SIAM J. Sci.
  Comput.}, 28:1078--1099, 2006.

\bibitem{Biegler_07}
Biegler, L.T., Ghattas, O., Heinkenschloss, M., Keyes, D., and van
  Bloemen~Waanders, B. (Eds.), {\em Real-Time PDE-Constrained Optimization},
  Vol.~3, SIAM Computational Science and Engineering, 2007.

\bibitem{Borzi_07}
Borzi, A., High-order discretization and multigrid solution of elliptic
  nonlinear constrained optimal control problems, {\em J. Comp. Applied Math},
  200:67--85, 2007.

\bibitem{Hinze_09}
Hinze, M., Pinnau, R., Ulbrich, M., and Ulbrich, S., {\em Optimization with PDE
  Constraints}, Springer, 2009.

\bibitem{Biegler_11}
Biegler, L., Biros, G., Ghattas, O., Heinkenschloss, M., Keyes, D., Mallick,
  B., Marzouk, Y., Tenorio, L., van Bloemen~Waanders, B., and Willcox, K.
  (Eds.), {\em Large-Scale Inverse Problems and Quantification of Uncertainty},
  John Wiley and Sons, 2011.

\bibitem{frontier_in_pdeco}
Antil, H., Kouri, D.P., Lacasse, M., and Ridzal, D. (Eds.), {\em Frontiers in
  PDE-Constrained Optimization}, Springer, 2018.

\bibitem{shapiro_SIAM_review}
Bonnans, J.F. and Shapiro, A., Optimization problems with perturbations: A
  guided tour, {\em SIAM Review}, 40(2):228--264, 1998.

\bibitem{fiacco}
Fiacco, A.V. and Ghaemi, A., Sensitivity analysis of a nonlinear structural
  design problem, {\em Computers and Operations Research}, 9(1):29--55, 1982.

\bibitem{Griesse_part_1}
Griesse, R., Parametric sensitivity analysis in optimal control of a reaction
  diffusion system. {I}. solution differentiability, {\em Numerical Functional
  Analysis and Optimization}, 25(1-2):93--117, 2004.

\bibitem{Griesse_part_2}
Griesse, R., Parametric sensitivity analysis in optimal control of a
  reaction-diffusion system -- part {II}: practical methods and examples, {\em
  Optimization Methods and Software}, 19(2):217--242, 2004.

\bibitem{griesse2}
Brandes, K. and Griesse, R., Quantitative stability analysis of optimal
  solutions in {PDE}-constrained optimization, {\em Journal of Computational
  and Applied Mathematics}, 206:908--926, 2007.

\bibitem{ghattas_infinite_dim_bayes_1}
Bui-{T}hanh, T., Ghattas, O., Martin, J., and Stadler, G., A computational
  framework for infinite-dimensional {B}ayesian inverse problems. {P}art {I}:
  The linearized case, with applications to global seismic inversion, {\em SIAM
  Journal on Scientific Computing}, 35(6):A2494--A2523, 2013.

\bibitem{ghattas_infinite_dim_bayes_2}
Petra, N., Martin, J., Stadler, G., and Ghattas, O., A computational framework
  for infinite-dimensional {B}ayesian inverse problems. {P}art {II}: Stochastic
  newton mcmc with application to ice sheet flow inverse problems, {\em SIAM
  Journal on Scientific Computing}, 36(4):A1525--A1555, 2014.

\bibitem{stuart_inv_prob}
Stuart, A.M., Inverse problems: a {B}ayesian perspective, {\em Acta Numerica},
  2010.

\bibitem{stuart_learning_lin_op}
de~Hoop, M.V., Kovachki, N.B., Nelsen, N.H., and Stuart, A.M., Convergence
  rates for learning linear operators from noisy data, {\em SIAM/ASA Journal on
  Uncertainty Quantification}, 11(2), 2023.

\bibitem{lindgren11}
Lindgren, F., Rue, H., and Lindstr{\"o}m, J., An explicit link between
  {G}aussian fields and {G}aussian {M}arkov random fields: The stochastic
  partial differential equation approach, {\em Journal of the Royal Statistical
  Society Series B}, 73:423--498, 2011.

\bibitem{matern_cov_villa}
Villa, U. and O'{L}eary-Roseberry, T., A note on the relationship between
  {PDE}-based precision operators and matern covariances, {\em
  arXiv:2407.00471}, 2024.

\end{thebibliography}
	\end{document}